\documentclass[a4paper,12pt]{amsart}

\usepackage[margin=1in]{geometry} 

\usepackage{amssymb,amsmath,amsthm,mathrsfs,enumerate,multicol,esint}
\usepackage{hyperref}
\usepackage{color}
\usepackage{makecell}
\usepackage{tikz}
\usetikzlibrary{arrows,calc}

\allowdisplaybreaks

\title[Hermite Calder\'on--Zygmund operators on endpoint spaces]{Calder\'on--Zygmund operators and endpoint spaces for Hermite expansions}

\author{The Anh Bui}
\address{School of Mathematical and Physical Sciences, 
Macquarie University, 
NSW 2109, Australia}
\email{the.bui@mq.edu.au}

\author{Fu Ken Ly}
\address{School of Mathematics and Statistics, The Learning Hub, The University of Sydney, 
NSW 2006, Australia}
\email{ken.ly@sydney.edu.au}


\newcommand{\RR}{\mathbb{R}} %Euclidean space
\newcommand{\NN}{\mathbb{N}} % Natural numbers
\newcommand{\ZZ}{\mathbb{Z}} %integers
\newcommand{\sz}{\mathscr{S}} % schwartz
\newcommand{\MM}{\mathcal{M}}

\newcommand{\N}{\mathcal{N}}
\newcommand{\K}{\mathcal{K}}
\newcommand{\Lip}{\Lambda}
\newcommand{\vc}{\infty}
\newcommand{\f}{\frac}
\newcommand{\lesi}{\lesssim}

\newcommand{\LL}{\mathcal{L}}
\newcommand{\QQ}{\mathbb{Q}} %projection
\newcommand{\ip}[1]{\langle #1 \rangle} % inner product
\newcommand{\wt}[1]{\widetilde{#1}} %shortcut for widetilde
\newcommand{\diff}{\triangle}
\newcommand{\ve}{\varepsilon}
\newcommand{\vp}{\varphi}

\newcommand{\CN}{\mathcal{C}} %cancellation class
\newcommand{\SM}{\mathcal{S}} %smoothness class - hermite
\newcommand{\floor}[1]{\lfloor #1 \rfloor}
\newcommand{\ceil}[1]{\lceil #1\rceil}
\newcommand{\A}[1]{A^{(#1)}}
\newcommand{\AH}{\mathcal{A}}
\newcommand{\Axy}{\big(\mathcal{A}_i^{(x)}-\mathcal{A}_i^{(y)}\big)}
\newcommand{\vph}{\varphi}
\newcommand{\cro}{\varrho} %critical radius
\newcommand{\rr}{\rho} %hormander class index for frequency variable

\newcommand{\bmo}{\text{\upshape bmo}}

\newcommand{\HCZO}{\text{\upshape HCZO}}

\newcommand{\dd}{\delta}

\newcommand{\om}{\omega}
\newcommand{\VV}{\mathcal{V}}
\newcommand{\UU}{\mathcal{U}}
\newcommand{\I}{\mathcal{I}}
\newcommand{\half}{\frac{1}{2}}
\newcommand{\R}{\mathcal{R}}

\newcommand{\D}{\partial}
\newcommand{\uu}{\upsilon}

\newcommand{\pt}{e^{-t\LL}}
\newcommand{\F}{\mathcal{F}}

\DeclareMathOperator{\supp}{supp\,} %support of a function
\DeclareMathOperator{\proj}{Proj} % projection operator
\theoremstyle{plain}
\newtheorem{Theorem}{Theorem}[section]
\newtheorem{Lemma}[Theorem]{Lemma}
\newtheorem{Proposition}[Theorem]{Proposition}

\newtheorem{Definition}[Theorem]{Definition}
\theoremstyle{definition}
\newtheorem{Remark}[Theorem]{Remark}
\theoremstyle{remark}

\numberwithin{equation}{section}%allows numbering of equations based on section

\def\barint{\kern4pt
\raise3.4pt\hbox{\vrule height.8pt width5pt}%
\kern-9pt % -(4pt + 5pt)
\int}
\def\XXint#1#2#3{{\setbox0=\hbox{$#1{#2#3}{\int}$}
     \vcenter{\hbox{$#2#3$}}\kern-.5\wd0}}

\begin{document}

\begin{abstract}
Let $\LL-\Delta +|x|^2$ be the Hermite operator on $\mathbb R^n$, and $T$ be a Calderon--Zygmund  type operator that is modelled on certain singular integrals related to $\LL$.
We establish necessary and sufficient conditions for $T$ to be bounded on various function spaces including the Hardy spaces and the  Lipschitz spaces associated to $\mathcal L$. We then apply our results to study the boundedness of the Riesz transforms and pseudo-multipliers associated to $\mathcal L$.
 \end{abstract}

\maketitle

%\tableofcontents

\section{Introduction}
The purpose of this work is to shed further light on the molecular and cancellation theory for operators and endpoint spaces (which include the Hardy and Lipschitz spaces) in the context of Hermite expansions. 

It is well known that a Calder\'on--Zygmund operator $T$ is bounded on the Hardy space $H^p(\RR^n)$ for the range $\f{n}{n+1}<p\le 1$ if and only if $T^*(1)=0$. For $p<\f{n}{n+1}$, one requires additional regularity on $T$, as well as higher order cancellations in the form of 
\begin{align}\label{eq:T1-zero}
T^*(x^\gamma) =0, \qquad\text{for}\quad |\gamma|\le \floor{n(1/p-1)}.
\end{align}
For the other endpoints, which include (homogeneous) Lipschitz spaces $\dot{\Lip}_s(\RR^n)$ and various spaces of mean oscillation, similar conditions are required on $T$ itself, such as $T(1)=0$. Some references include \cite{AM,FTW,HL,HHL,Lemarie,LW,Meyer,Torres}.

Unfortunately, there exist examples of Calder\'on--Zygmund operators for which $T(1)$ or $T^*(1)$ are non-vanishing. One class of examples includes the usual pseudo-differential operators with symbols from the H\"ormander classes $S^0_{\rr,\dd}$, $0\le \dd\le\rr<1$, which are known to preserve the local Hardy spaces $h^p(\RR^n)$ of Goldberg \cite{Go}, but not $H^p(\RR^n)$. 
Another class of examples includes various singular integrals associated to operators with potential  such as the Hermite operator $=-\Delta+|x|^2$ and, more generally, Schr\"odinger operators $-\Delta+V$ with certain potentials $V$; again, such operators are well behaved on endpoint spaces derived from the underlying  operator, but not necessarily on those in the classical scales.

This raises the question of what are suitable cancellation conditions in place of \eqref{eq:T1-zero} for situations that lie beyond the classical settings. 
For $h^p(\RR^n)$, there has been some quite recent progress in this line of inquiry and the picture that has emerged is a kind of Lipschitz or Campanato-type condition. A typical result here states that if a Calder\'on--Zygmund operator $T$ has sufficient additional decay at infinity on its kernel, then it is bounded on $h^p(\RR^n)$ for $\f{n}{n+1}<p\le 1$ if and only if 
\begin{align}\label{eq:T1-lip}
T^*(1)\in \Lip_{n(1/p-1)}(\RR^n).
\end{align}
For $p$ below $\f{n}{n+1}$ one has analagous higher moment conditions; see \cite{BL22, DHZ, DLPV21, DLPV22, Ko} for further details.

Likewise, there has been related developments for operators with potentials. These include critiera for so-called Hermite-- or Schr\'odinger--Calder\'on--Zygmund operators to be bounded on various endpoint spaces including their associated Hardy and BMO type spaces (see \cite{BCFST, BHQ, BLL18,MSTZ}). However, to date, these results have dealt only with the range $\f{n}{n+1}<p\le 1$ (for Hardy spaces) and $0\le s<1$ (for Lipschitz spaces). It is therefore interesting to ask whether one can develop a suitable cancellation theory in such settings for small $p$ and large smoothness index $s$ in accordance with the inhomogeneous regime of $h^p(\RR^n)$ and $\Lip_s(\RR^n)$. 

\bigskip
The aim of the present article is to give a positive answer to the above question in the context of Hermite expansions. In particular, we are able to establish a theory for the full range $0<p\le 1$ and $s\ge 0$ and we  also provide some examples that illustrate how our criteria may be applied.

In order to describe our main results, let us introduce some definitions and notations relevant to the Hermite context. 
For $n=1$ the Hermite function of degree $k\in\NN_0=\NN\cup\{0\}$ is 
\begin{align*}
	h_k(t)=(2^k k! \sqrt{\pi})^{-1/2}H_k(t)e^{-t^2/2} \qquad \forall t\in\RR{},
\end{align*}
where $H_k(t)=(-1)^ke^{t^2}\partial_t^k(e^{-t^2})$
is the $k$th Hermite polynomial. 
For $n\ge 2$, the $n$-dimensional Hermite functions $h_\xi$ are defined over the multi-indices $\xi\in \NN_0$ by 
\begin{align*}
	h_\xi(x)=\prod_{j=1}^n h_{\xi_j}(x_j)\qquad \forall x\in\RR^n.
\end{align*}
Such functions are eigenfunctions of the Hermite operator
$$\LL=-\Delta+|x|^2$$
 in the sense that $\LL(h_\xi)=(2|\xi|+n)h_\xi$ for every $\xi\in\NN_0^n$;  in addition, they form an orthonormal basis for $L^2(\RR^n)$. It is also worth noticing here that $\LL$ has a non-zero spectral gap.

A useful tool for analysis in the Hermite setting is the so-called `critical radius' function $\cro:\RR^n\to \RR$  given by
\begin{align*}
\cro(x):=\f{1}{1+|x|},\qquad x\in\RR^n.
\end{align*}
Observe that for any $x,y\in\RR^n$ we have $\cro(y)\sim\cro(x)$ whenever $|x-y|<\cro(x)$.

We now introduce the following notions of (higher-order) Hermite--Calder\'on--Zygmund operators, which generalise previous definitions as found in \cite{BCFST,MSTZ}.

\begin{Definition}[Hermite--Calder\'on--Zygmund operators]\label{def:HCZO}
Suppose that $T$ is a continuous linear operator from $\mathcal{D}(\RR^n)$ to $\mathcal{D}'(\RR^n)$ with an associated kernel $K(x,y)$. 
Let $M\in\NN_0$ and $\ve \in (0,1]$. Then we say $T$ is a $(M,\ve)$-Hermite--Calder\'on--Zygmund operator (denoted $T\in \HCZO(M,\ve)$) if all the following conditions hold.
\begin{enumerate}[\upshape(i)]
\item $T$ is bounded on $L^2(\RR^n)$.
\item For $x\ne y$ we have
\begin{align*}
|K(x,y)|\lesssim  |x-y|^{-n}\Big(1+\f{|x-y|}{\cro(y)}\Big)^{-M-\ve}.
\end{align*}
\item For $|x-y|>2|x-x'|$,
\begin{align*}
|\partial_1^\gamma K(x,y)-\partial_1^\gamma K(x',y)|\lesssim  \f{|x-x'|^{\ve}}{|x-y|^{n+M+\ve}}, \qquad |\gamma|= M.
\end{align*}
\item For $|x-y|>2|y-y'|$,
\begin{align*}
|\partial_2^\gamma K(x,y)-\partial_2^\gamma K(x,y')|\lesssim  \f{|y-y'|^{\ve}}{|x-y|^{n+M+\ve}}, \qquad |\gamma|= M.
\end{align*}
\end{enumerate}
If $T$ satisfies (i),(ii) and (iii) only, then we say $T\in \HCZO_{(1)}(M,\ve)$; if $T$ satisfies (i), (ii) and (iv) only, then we say $T\in \HCZO_{(2)}(M,\ve)$. 
\end{Definition}
Here and in the sequel, for a multi-index $\gamma\in\NN_0^n$, we denote by $\partial^\gamma_i =\partial^{\gamma_1}_i\partial^{\gamma_2}_i\cdots\partial^{\gamma_n}_i$ to mean the partial derivatives with respect to the $i$th variable.
It is clear that a Hermite-Calderon--Zygmund operator is a Calder\'on--Zygmund operator with stronger decay at infinity. For existing work on $(0,\ve)$-order operators and some relevant examples one can consult \cite{BCFST,BLL18,Ky,MSTZ,ST05a,ST03,ST05b,T93} and the references therein.

\bigskip
We are now ready to present the main results of this article. Our first  result gives necessary and sufficient conditions for a Hermite--Calder\'on--Zygmund operator (with regularity in the second variable) to be bounded on the Hermite--Hardy spaces $h^p_\LL(\RR^n)$ for all $0<p\le 1$ (for a definition of these spaces, see Definition \ref{def:hp-spaces}). Observe that, in analogy with \eqref{eq:T1-lip}, the conditions are given in terms of the Hermite--Lipschitz spaces $\Lip^s_\LL(\RR^n)$, $s\ge0$ (see Definition \ref{def:Lip-spaces}). 
\begin{Theorem}[Hermite Calder\'on--Zygmund operators on Hardy spaces]\label{thm:hardy}
Let $\om\ge0$, $\ve\in(0,1]$ and $T\in\HCZO_{(2)}\big(\floor{\om},\ve\big)$. Consider the following condition 
\begin{align}\label{eq:hardycond}
\sup_{x_0\in\RR^n} \cro(x_0)^{\om-|\alpha|} \big\Vert T^*[(\cdot-x_0)^\alpha \chi]\big\Vert_{\Lip^{\om}_\LL} <\infty, \qquad \forall\;|\alpha|\le \floor{\om}
\end{align}
where $\chi\in C^\infty_0(B(x_0,2\cro(x_0))$ with $\chi=1$ on $B(x_0,\cro(x_0))$. 
\begin{enumerate}[\upshape(a)]
\item
If \eqref{eq:hardycond} holds for some $\om>0$ with $\om^*\ne0$, then $T$ extends to a bounded operator on the Hardy space $h^p_\LL(\RR^n)$ for every  $\f{n}{n+\floor{\om}+\om^*\land\ve}<p\le 1$.
\item Conversely, if $T$ is bounded on the Hardy space $h^p_\LL(\RR^n)$ for some $0<p\le 1$, then \eqref{eq:hardycond} holds for $\om=n(\f{1}{p}-1)$. 
\end{enumerate}
\end{Theorem}
\noindent Note here that for $\om^*:=\om-\floor{\om}$ for $\om\in\RR$.

Our second main result establishes necessary and sufficient conditions for a Hermite--Calder\'on--Zygmund operator (with regularity in the first variable) to be bounded on the Hermite--Lipschitz spaces.
\begin{Theorem}[Hermite Calder\'on--Zygmund operators on Lipschitz spaces]\label{thm:lip}
Let $\om\ge0$, $\ve\in(0,1]$ and $T\in\HCZO_{(1)}\big(\floor{\om},\ve\big)$. Consider the following condition 
\begin{align}\label{eq:lipcond}
\sup_{x_0\in\RR^n} \cro(x_0)^{\om-|\alpha|} \big\Vert T[(\cdot-x_0)^\alpha \chi]\big\Vert_{\Lip^{\om}_\LL} <\infty, \qquad \forall\;|\alpha|\le \floor{\om}
\end{align}
where $\chi\in C^\infty_0(B(x_0,2\cro(x_0))$ with $\chi=1$ on $B(x_0,\cro(x_0))$. 
\begin{enumerate}[\upshape(a)]
\item
If \eqref{eq:lipcond} holds for some $\om>0$ with $\om^*\ne0$, then $T$ extends to a bounded operator on the Lipschitz space $\Lip^s_\LL(\RR^n)$ for every  $0\le s<\floor{\om}+\om^*\land\ve$.
\item Conversely, if $T$ is bounded on the Lipschitz space $\Lip^s_\LL(\RR^n)$ for some $s\ge 0$, then \eqref{eq:lipcond} holds for $\om=s$. 
\end{enumerate}
\end{Theorem}
We wish to emphasise that our results above extend the indices previously considered in the literature which were $\f{n}{n+1}< p\le 1$ and $0\le s<1$. It is also worth pointing out that part (a) of Theorem \ref{thm:hardy} can be reformulated in the following way. Let $p\in (0,1]$, $(\f{n}{p})^*<\ve\le 1$ and suppose that $T\in \HCZO_{(2)}\big(\floor{n(\f{1}{p}-1)},\ve\big)$. If $T$ satisfies \eqref{eq:hardycond} with $\om = \floor{n(\f{1}{p}-1)}+\theta$ for some $(\f{n}{p})^*<\theta\le 1$ then $T$ is bounded on $h^p_\LL(\RR^n)$. 

Let us offer a few remarks on our proofs. The two crucial tools in our approach are duality (Theorem \ref{thm-duality}) and a new molecular characterisation for Hermite--Hardy spaces (Proposition \ref{prop: molecules}). One of the key impediments with working with small $p$ in the inhomogeneous or Schr\"odinger (and thus Hermite) contexts is a suitable molecular theory for the range $p<\f{n}{n+1}$ and, while the atomic theory for the Hermite--Hardy spaces has been available since \cite{Dz}, a satisfying molecular theory has thus far been out of reach. 
In this article we establish a new molecular characterisation for  Hermite--Hardy spaces where each molecule $m$, whose mass is concentrated around some ball $B(x_0,r)$, satisfies an `approximate cancellation' of the form
$$ \Big|\int (x-x_0)^\gamma m(x)\,dx\Big|\lesi r^{n-n/p+\om}\cro(x_0)^{|\alpha|-\om}.$$ 
for all $|\gamma|\le \om$ (see Definition \ref{def:molecule}). The antecedents of these molecules are inspired by recent developments in \cite{BL22, DLPV21, LN21} (see also \cite{BLL18,DHZ,Ko}) where molecules with similar moment conditions have been employed for $h^p(\RR^n)$ and for other function spaces in the Hermite context.

In the final section of this paper, we show how our main results may be applied to some operators that have attracted recent interest, namely, the so-called Hermite pseudo-multipliers and the higher-order Hermite Riesz transforms (Section \ref{sec:applications}). We show that these are Hermite--Calder\'on--Zygmund operators in the sense of Definition \ref{def:HCZO} and apply Theorems \ref{thm:hardy} and \ref{thm:lip} to recover some of the known results in the literature. One particular case however, is new: the action of pseudo-multipliers on $\bmo_\LL$, the space of bounded mean oscillation for Hermite operators (which is the space $\Lip^0_\LL$).

The paper is organised as follows. In Section \ref{sec:spaces} we introduce and establish some key results  concerning the Hermite--Hardy and Hermite--Lipschtiz spaces including a new molecular characterisation (Proposition \ref{prop: molecules}) and duality (Theorem \ref{thm-duality}). Our main results are then proved in Section \ref{sec:proofs}. Finally, in Section \ref{sec:applications} we give applications to Hermite pseudo-multipliers and higher-order Riesz transforms. Some of the proofs in this final section are fairly technical and we place these---along with other important facts and identities in the Hermite context---in the Appendices.

\bigskip
{\bf Notation:} 
We set$\NN_0:=\NN\cup\{0\}$; $\ZZ^-$ will denote the negative integers, that is, $\ZZ^-=\ZZ\backslash\NN_0$, and $\ZZ^+$ will denote the positive integers. For $s\in\RR$ we define $\floor{s}$ to mean the greatest integer not exceeding $s$, and also define $s^*:=s-\floor{s}$. 
The usual bracket notation $\ip{x}=1+|x|$ will also be employed, and the critical radius function can be thus defined as $\cro(x):=\ip{x}^{-1}$. For a locally integrable function $f$ and measureable set $E\subset\RR^n$, we define $\displaystyle\fint_E f :=\f{1}{|E|}\int f$ to mean the average of $f$ over $E$. 
We denote by $B(x,r)$ the ball centred $x\in\RR^n$ and radius $r>0$. By a 'ball $B$' we mean $B
=B(x_B,r_B)$ with some fixed centre $x_B$ and radius $r_B$. Once a ball has been fixed from the context, we will also use the shorthand $\cro_B:=\cro(x_B)$.
For a ball $B$, the set $U_j(B)$ will denote the annulus $2^j B\backslash 2^{-j-1}B$ when $j\ge 1$, and denote $B$ when $j=0$.

\section{Hermite endpoint spaces}\label{sec:spaces}
In this section we define the endpoint spaces under consideration, namely the Hermite--Hardy and Hermite--Lipschitz spaces, and then establish some important results concerning these spaces that will be needed in the proofs of the main results.

To begin we outline some needed estimates on the Hermite heat semigroup and its kernel estimates. 
The Hermite heat semigroup $\{e^{-t\LL}\}_{t>0}$ is defined by
$$ e^{-t\LL}f(x)=\sum_{\xi\in \NN_0^n} e^{-t(2|\xi|+n)} \ip{f,h_\xi}h_\xi(x), \qquad f\in L^2(\RR^n),$$
with kernel given by
$$ \pt(x,y) = \sum_{\xi\in \NN_0^n} e^{-t(2|\xi|+n)} h_\xi(x)h_\xi(y).$$
For further details see \cite{BD,ST03,T93}.

In the sequel we will be making use of the following estimates.
\begin{Lemma}[Heat kernel estimates]\label{lem:HK}
Let $k\in\NN_0$. 
\begin{enumerate}[\upshape(a)]
\item
There exist constants $C_k, c>0$ such that for any $N>0$ and $x,y\in\RR^n$,
\begin{align*}
\big|\LL^k\pt(x,y)\big| \le C_k \f{e^{-\f{|x-y|^2}{ct}}}{t^{n/2+k}}\Big(1+\f{\sqrt{t}}{\cro(x)}+\f{\sqrt{t}}{\cro(y)}\Big)^{-N}.
\end{align*}
\item
For each $\gamma\in \NN_0^n$, there exist $C_\gamma,c>0$ such that for any $N>0$ and $x,y\in\RR^n$,
\begin{align*}
\big|\partial_y^\gamma \pt(x,y)\big| \le C_\gamma \f{e^{-\f{|x-y|^2}{ct}}}{t^{(n+|\gamma|)/2}}\Big(1+\f{\sqrt{t}}{\cro(x)}+\f{\sqrt{t}}{\cro(y)}\Big)^{-N}.
\end{align*}
\end{enumerate}
\end{Lemma}
\begin{proof}
The proof of (a) can be found in \cite[Proposition 4]{DGMTZ} and \cite[Lemma 2.1]{DT21}. See also Lemmas 2.1 and 2.2 of \cite{BD}. 
The proof of (b) can be found in \cite[Proposition 3.3]{BHH}.
\end{proof}

We will also utilise some estimates involving certain orthogonal projectors associated with the Hermite functions. For each $N\in\NN_0$  we  define the orthogonal projection of $f$ onto $\bigoplus^N_{k=0} \text{span}\{h_\xi: |\xi|=k\}$ by
\begin{align*}
	\QQ_N f = \sum_{k=0}^N \sum_{|\xi|=k} \ip{f, h_\xi} h_\xi
	\qquad\text{with kernel}\qquad
	\QQ_N(x,y)=\sum_{k=0}^N \sum_{|\xi|=k}h_\xi(x)h_\xi(y).
\end{align*}
Then we have the following bounds on $\QQ_N$ (see \cite{LN21,Ly22,PX}):
 for any $N\in\NN_0$ and $\mu\ge 0$, we have
\begin{align}\label{eq:QQ est}
\QQ_{4^j+N}(x,x)\lesssim 2^{jn}\Big(1+\f{1+|x|}{2^j}\Big)^{-\mu} \qquad \forall \; j\in\NN_0, \;x\in\RR^n.
\end{align}

\subsection{Hardy spaces and a new molecular characterization}
In this section we define Hardy spaces associated with the Hermite operator and give a new notion of molecule for all $0<p\le 1$. 

Let $\MM_\LL$ be the maximal operator associated to the Hermite heat semigroup given by
$$ \MM_\LL f(x):=\sup_{t>0}|e^{-t\LL}f(x)|.$$

\begin{Definition}[Hermite Hardy spaces]\label{def:hp-spaces}
For $0<p\le 1$ we define the Hermite Hardy space $h^p_\LL(\RR^n)$ by
$$h^p_\LL(\RR^n) := \big\{ f\in \mathscr{S}': \MM_\LL f\in L^p(\RR^n)\big\}$$
and set $\Vert f\Vert_{h^p_\LL}:=\Vert \MM_\LL f\Vert_{L^p}$. 
\end{Definition}
It is well known that these spaces admit the following atomic decomposition (see \cite{Dz}).
\begin{Definition}[Atoms]
Let $0<p\le 1<q\le \infty$ and $M\in\NN_0$. We say that the function $a$ is a  $(p,q,M)$-atom associated to a ball $B$ if 
\begin{enumerate}[\upshape(i)]
\item $r_B\le \f{1}{2}\cro_B$;
\item $\supp a\subset B$;
\item $\Vert a\Vert_{L^q}\le |B|^{\f{1}{q}-\f{1}{p}}$;
\item If $r_B<\f{1}{8}\cro_B$ then $\displaystyle \int x^\alpha a(x)\,dx=0$ for every $|\alpha|\le M$.
\end{enumerate}
\end{Definition}

Then for $0<p\le 1$ and $M\ge \floor{n(\f{1}{p}-1)}$ we have that $f\in h^p_\LL(\RR^n)$ if and only if $f$ can be represented as $f=\sum_j c_j a_j$ where each $a_j$ is a $(p,q,M)$ atom and $c_j$ are scalars satisfying $\sum_j|c_j|^p<\infty$. In addition, $\Vert f\Vert_{h^p_\LL}\sim \inf \big\{(\sum_j |c_j|^p)^{1/p}\big\}$ where the infimum is taken over all possible representations $f=\sum_jc_ja_j$.   See for example \cite[Theorem 1.12]{Dz} and \cite[Theorem 1.4]{BHH}.

We now introduce the following  new notion of molecules for $h^p_\LL$.  Molecules for $\f{n}{n+1}<p\le 1$ were previously given in \cite{BLL18}.
\begin{Definition}[Molecules]\label{def:molecule}
Let $0<p\le 1<q\le \infty$,  $\delta>0$ and $\om\ge 0$. We say that the function $m$ is a  $(p,q,\delta,\om)$-molecule associated to a ball $B$ if 
\begin{enumerate}[\upshape(i)]
\item $r_B\le \f{1}{2}\cro_B$;
\item $\Vert m\Vert_{L^q(U_j(B))}\le 2^{-j\delta}|2^jB|^{\f{1}{q}-\f{1}{p}}$ for every $j\in\NN_0$;
\item If $r_B<\f{1}{8}\cro_B$ then $\displaystyle \Big|\int (x-x_B)^\alpha m(x)\,dx\Big| \le |B|^{1-\f{1}{p}} r_B^{|\alpha|} \Big(\f{r_B}{\cro_B}\Big)^{\om-|\alpha|}$ for every $|\alpha|\le \floor{\om}$.
\end{enumerate}
\end{Definition}

Then we have the following
\begin{Proposition}[Molecular characterisation]\label{prop: molecules}
Let $0<p<1<q\le \infty$, $\om\ge0$ and $\delta>\max\{0,\floor{\om}-n(\f{1}{p}-1)\}$.  
If $m$ is a $(p,q,\delta,\om)$-molecule  then $m\in h^p_\LL(\RR^n)$ for all $\f{n}{n+\om}<p\le 1$. \end{Proposition}

Our proof of Proposition \ref{prop: molecules} will follow the approach in \cite{BL22}, but adapted to the Hermite context.
We shall require the following lemma.
\begin{Lemma}\label{lem:AE}
Let  $s>0$ and $0<p<1$. Suppose that $b$ is a function supported in a ball $B$ with $r_B\le\cro_B$ such that 
\begin{align}\label{eq:AE0}
\Vert b\Vert_{L^q}\le |B|^{\f{1}{q}-\f{1}{p}}
\end{align}
for some $q\ge 1$. Assume $b$ also satisfies one of the following conditions:
\begin{enumerate}[\upshape (a)]
\item There exists $\lambda\in (0,1]$ such that $\lambda \cro_B\le r_B \le \cro_B$; or
\item For every $|\alpha|\le \floor{s}$ we have $\displaystyle \int b(x)(x-x_B)^\alpha dx=0$; or
\item For some $|\alpha|\le \floor{s}$ we have
$$ \int b(x)(x-x_B)^\gamma dx=0, \qquad \text{for}\quad\gamma\ne \alpha,\quad |\gamma|\le \floor{s} $$
and for $\gamma =\alpha$,
$$ \Big|\int b(x)(x-x_B)^\alpha dx\Big|\le |B|^{1-\f{1}{p}} \Big(\f{r_B}{\cro_B}\Big)^{s-|\alpha|}r_B^{|\alpha|}.$$
\end{enumerate}
Then we have 
\begin{align}\label{eq:AE1}
|e^{-t\LL}b(x)|\lesi \f{r_B^s}{|x-x_B|^{n+s}} |B|^{1-\f{1}{p}},\qquad x\in\RR^n\backslash 4B,
\end{align}
and as a consequence, for any $\f{n}{n+s} < \wt{p} \le 1$,
\begin{align}\label{eq:AE2}
\Vert \MM_\LL(b)\Vert_{L^{\wt{p}}}\lesi |B|^{1/\wt{p}-1/p}.
\end{align}
\end{Lemma}
\begin{Remark}\label{rem:AE}
Given  $M\in\NN_0$ and any $(p,q,M)$-atom $a$, we have for any $\f{n}{n+M+1}<\wt{p}\le 1$,
$$ \Vert a\Vert_{h^{\wt{p}}_\LL} \lesi |B|^{1/\wt{p}-1/p}.$$
One can see this by applying Lemma \ref{lem:AE}  with $s=M+1$ and $\lambda=1/8$.

\end{Remark}
\begin{proof}[Proof of Lemma \ref{lem:AE}]
We first estimate \eqref{eq:AE1}. Under assumption (a) we have $\Vert b\Vert_{L^1}\le |B|^{1-1/p}$ from \eqref{eq:AE0} and H\"older's inequality. Note also that for $x\in \RR^n\backslash 4B$ and $y\in B$ we have $|x-y|\sim |x-x_B|$. Since $r_B\le \cro_B$ we also have $\cro(y)\sim \cro_B$. Then by Lemma \ref{lem:HK} we have
\begin{align*}
\big|e^{-t\LL}b(x)\big|
&\lesi \int_B \f{e^{-|x-y|^2/ct}}{t^{n/2}}\Big(1+\f{\sqrt{t}}{\cro(x)}+\f{\sqrt{t}}{\cro(y)}\Big)^{-N} |b(y)|\,dy \\
&\lesi \f{e^{-|x-x_B|^2/ct}}{t^{n/2}}\Big(1+\f{\sqrt{t}}{\cro(x)}+\f{\sqrt{t}}{\cro_B}\Big)^{-N} \Vert b\Vert_{L^1} \\
&\lesi \f{\cro_B^N}{|x-x_B|^{n+s}} t^{\f{s-N}{2}} |B|^{1-1/p}
\end{align*}
By setting $N=s$ and observing that $\cro_B\le \lambda^{-1} r_B$ we obtain \eqref{eq:AE1}.

We now turn to assumptions (b) and (c).  Let us write
\begin{align*}
e^{-t\LL}b(x) 
&=\int_B \big[ \pt(x,y)-\sum_{|\gamma|\le \floor{s}} \f{1}{\gamma !} \partial_2^\gamma \pt(x,x_B) (y-x_B)^\gamma \big] b(y)\,dy \\
&\qquad+ \sum_{|\gamma|\le  \floor{s}} \f{1}{\gamma !} \partial_2^\gamma \pt(x,x_B) \int_B (y-x_B)^\gamma b(y)\,dy
\qquad=: I +II
\end{align*}
For term $I$ we have, by Lemma \ref{lem:HK} and  Taylor's theorem, for some $\wt{y}$ on the line segment between $y$ and $x_B$,
\begin{align*}
|I| \le \sum_{|\gamma|=\floor{s}+1} \f{1}{\gamma!}\int_B \big|\partial_2^\gamma \pt(x,\wt{y})\big| |y-x_B|^{|\gamma|} |b(y)|\,dy 
\lesi \int_B \f{e^{-|\wt{y}-x_B|^2/ct}}{\sqrt{t}^{n+\floor{s}+1}} |y-x_B|^{\floor{s}+1}|b(y)|\,dy. 
\end{align*}
Since $\wt{y}\in B$ and $x\in \RR^n\backslash 4B$ then $|\wt{y}-x_B|\sim |x-x_B|$. Thus
\begin{align*}
|I| 
\lesi \int_B \f{e^{-|x-x_B|^2/ct}}{\sqrt{t}^{n+\floor{s}+1}} |y-x_B|^{\floor{s}+1}|b(y)|\,dy 
\lesi \f{r_B^{\floor{s}+1}}{|x-x_B|^{n+\floor{s}+1}} \Vert b\Vert_{L^1}
\le \f{r_B^{\floor{s}+1}}{|x-x_B|^{n+\floor{s}+1}}  |B|^{1-\f{1}{p}}.
\end{align*}
Then, since $|x-x_B|\ge 3 r_B$ and $\floor{s}+1-s\ge 0$, we have 
\begin{align*}
\f{r_B^{\floor{s}+1}}{|x-x_B|^{n+\floor{s}+1}}
= \f{r_B^{\floor{s}+1}}{|x-x_B|^{n+s}} \f{1}{|x-x_B|^{\floor{s}+1-s}}
\lesi 
\f{r_B^s}{|x-x_B|^{n+s}}.
\end{align*}
Thus we arrive at the required estimate in \eqref{eq:AE1} for term $I$.
We turn to term $II$. Observe that assumption (b) implies $II=0$ immediately, completing the estimate \eqref{eq:AE1} for case (b).

For case (c), all the integrals in $II$ vanish except for $\gamma=\alpha$. Thus we have by Lemma \ref{lem:HK},
\begin{align*}
|II|
&=\f{1}{\alpha !} \big|\partial_2^\alpha \pt(x,x_B)\big| \Big|\int (y-x_B)^\alpha b(y)\,dy\Big|\\
&\lesi \f{e^{-|x-x_B|^2/ct}}{\sqrt{t}^{n+|\alpha|}} \Big(1+\f{\sqrt{t}}{\cro(x)}+\f{\sqrt{t}}{\cro_B}\Big)^{-N} \Big|\int (y-x_B)^\alpha b(y)\,dy\Big| \\
&\le \f{e^{-|x-x_B|^2/ct}}{\sqrt{t}^{n+|\alpha|}} \Big(\f{\cro_B}{\sqrt{t}}\Big)^N |B|^{1-\f{1}{p}}\Big(\f{r_B}{\cro_B}\Big)^{s-|\alpha|}r_B^{|\alpha|}.
\end{align*}
By choosing $N=s-|\alpha|$ and then applying the calculation
\begin{align*}
e^{-|x-x_B|^2/ct} \lesi \Big(\f{\sqrt{t}}{|x-x_B|}\Big)^{n+s},
\end{align*}
we arrive at the required estimate in \eqref{eq:AE1} for term $II$. 

Finally we consider \eqref{eq:AE2}. We write 
$$\Vert \MM_\LL (b)\Vert_{L^{\wt{p}}} \le \Vert \MM_\LL (b)\Vert_{L^{\wt{p}}(4B)} +\Vert \MM_\LL (b)\Vert_{L^{\wt{p}}(\RR^n\backslash 4B)}. $$
Then by the H\"older's inequality, the $L^q$ boundedness of $\MM_\LL$ and assumption \eqref{eq:AE0},
\begin{align*}
\Vert \MM_\LL (b)\Vert_{L^{\wt{p}}(4B)}^{\wt{p}}
\le \Vert \MM_\LL (b)\Vert_{L^q}^{\wt{p}} |4B|^{1-\wt{p}/q}
\lesi \Vert b\Vert_{L^q}^{\wt{p}} |B|^{1-\wt{p}/q} 
\le |B|^{1-\wt{p}/p}
\end{align*}
Next, from \eqref{eq:AE1}  and the fact that $\wt{p}>\f{n}{n+s}$, we have
\begin{align*}
\Vert \MM_\LL (b)\Vert_{L^{\wt{p}}(4B)}^{\wt{p}}
\lesi |B|^{\wt{p}(1-1/p)} r_B^{s\wt{p}} \int_{\RR^n\backslash 4B} |x-x_B|^{-(n+s)\wt{p}}\,dx
\lesi |B|^{1-\wt{p}/p}.
\end{align*}
Together the latter two estimate gives \eqref{eq:AE2}, which concludes the proof of the lemma.
\end{proof}

\begin{proof}[Proof of Proposition \ref{prop: molecules}]
Let $m$ be a molecule as specified in the Proposition, and associated with some ball $B$. Our aim is to demonstrate that for some constant $C>0$ one has
\begin{align}\label{eq:mol0}\Vert \MM_\LL (m) \Vert_{L^p}\le C.\end{align}
For each $j\ge 0$ we set $ \chi_j:=\chi_{U_j(B)}$ and $ m_j:=m\chi_j$
and write
\begin{align}\label{eq:mol0.1}
m =\sum_{j\ge 0} m \chi_{U_j(B)} = \sum_{j\ge 0}m_j.
\end{align}
Let us first consider the case when $\f{1}{8}\cro_B\le r_B\le \f{1}{2}\cro_B$. For each $j\ge 0$ we apply part (a) of Lemma \ref{lem:AE} to $b=2^{j\delta} m_j$ with ball $2^j B$ and $\wt{p}=p$.
This, along with \eqref{eq:mol0.1}, gives
\begin{align*}
\Vert \MM_\LL(m)\Vert_{L^p}\lesi \sum_{j\ge 0} \Vert \MM_\LL(m_j)\Vert_{L^p}=\sum_{j\ge 0}2^{-j\delta}\Vert \MM_\LL(b)\Vert_{L^p}\lesi\sum_{j\ge 0} 2^{-j\delta} <\infty
\end{align*}
and thus $m\in h^p_\LL$ for any $0<p\le 1$.

We now consider the case $r_B<\f{1}{8}\cro_B$, which will constitute the bulk of the remaining proof.
Here we adapt the technique of \cite{TW} to refine the decomposition \eqref{eq:mol0.1}  into a sum of scaled atoms. 

To begin, take $\VV_j$ to be the span of the polynomials $\big\{(x-x_B)^\alpha\big\}_{|\alpha|\le \floor{\om}}$ on $U_j(B)$ and $\UU_j$ to be the corresponding inner product space given by
$$ \ip{f,g}_j:=\fint_{U_j(B)} f(x)g(x)\,dx.$$
Let $\{\uu_{j,\alpha}\}_{|\alpha|\le \floor{\om}}$ be an orthonormal basis for $\VV_j$ obtained via the Gram--Schmidt process applied to $\big\{(x-x_B)^\alpha\big\}_{|\alpha|\le \floor{\om}}$
which, through homogeneity and uniqueness of the process, gives
\begin{align}\label{eq:mol1}
\uu_{j,\alpha}(x)=\sum_{|\beta|\le \floor{\om}} \lambda_{\alpha,\beta}^j (x-x_B)^\beta,
\end{align}
where for each $|\alpha|, |\beta|\le \floor{\om}$ we have
\begin{align}\label{eq:mol2}
|\uu_{j,\alpha}(x)|\le C\qquad\text{and}\qquad |\lambda_{\alpha,\beta}^j |\lesi (2^j r_B)^{-|\alpha|}. 
\end{align}
Let $\{\nu_{j,\alpha}\}_{|\alpha|\le \floor{\om}}$ be the dual basis of $\big\{(x-x_B)^\alpha\big\}_{|\alpha|\le \floor{\om}}$ in $\VV_j$; 
that is, it is the unique collection of polynomials such that
\begin{align}\label{eq:mol3} \ip{\nu_{j,\alpha}, (\cdot-x_B)^\beta}_j=\delta_{\alpha,\beta}, \qquad  |\alpha|, |\beta|\le \floor{\om}. \end{align}
satisfying (by \eqref{eq:mol2})
\begin{align}\label{eq:mol5}
\Vert\nu_{j,\alpha}\Vert_\infty \lesi (2^jr_B)^{-|\alpha|},\qquad \forall \; |\alpha|\le \floor{\om}.
\end{align}
Now let $P_j:=\proj_{\VV_j} (m _j)$ be the orthogonal projection of $m_j$ onto $\VV_j$. Then  we have
\begin{align}\label{eq:mol6}
P_j = \sum_{|\alpha|\le \floor{\om}} \ip{m_j, \uu_{j,\alpha}}_j \uu_{j,\alpha} = \sum_{|\alpha|\le \floor{\om}} \ip{m_j,(\cdot-x_B)^\alpha}_j \nu_{j,\alpha}.
\end{align}
Define also the numbers
\begin{align}\label{eq:mol7}
\N_{j,\alpha}:= \left
\lbrace 
	\begin{array}{ll}
	\sum_{k\ge j} |U_k(B)| \ip{m_k,(\cdot-x_B)^\alpha}_k \qquad &\text{for}\quad j\ge 1,\\
			\int m(x)(x-x_B)^\alpha\,dx \qquad &\text{for}\quad j=0.
	\end{array}
\right.
\end{align}
and observe that $|U_j(B)| \ip{m_k,(\cdot-x_B)^\alpha}_k = \N_{j,\alpha}-\N_{j+1,\alpha}$. Then using \eqref{eq:mol6}-\eqref{eq:mol7} one has
\begin{align*}
\sum_{j\ge 0}P_j
=\sum_{|\alpha|\le \floor{\om}} \sum_{j\ge 0}\ip{m_j,(\cdot-x_B)^\alpha}_j \nu_{j,\alpha}
=\sum_{|\alpha|\le \floor{\om}} \sum_{j\ge 0}\big(\N_{j,\alpha}-\N_{j+1,\alpha}\big)\f{\nu_{j,\alpha}}{|U_j(B)|},
\end{align*}
which we may then use, along with summation by parts, to obtain following decomposition of $m$,
\begin{align*}
m=\sum_{j\ge 0} a_j +\sum_{|\alpha|\le \floor{\om}} \sum_{j\ge 0} a_{j,\alpha} + \sum_{|\alpha|\le \floor{\om}} a_\alpha,
\end{align*}
where
\begin{align*}
a_j :=m_j-P_j, &&
a_{j,\alpha} :=\N_{j+1,\alpha}\Big(\f{\nu_{j+1,\alpha}}{|U_{j+1}(B)|}-\f{\nu_{j,\alpha}}{|U_{j}(B)|}\Big), &&
a_\alpha:=\f{\nu_{0,\alpha}}{|B|}\N_{0,\alpha}.
\end{align*}

Let us now outline the important properties of the functions in the above decomposition. 
For $a_j$ we observe that for $|\alpha|\le \floor{\om}$,
\begin{align}\label{eq:mol8}
\supp a_j \subset 2^jB, &&
 \int a_j(x)(x-x_B)^\alpha dx=0, &&
  \Vert a_j\Vert_{L^q}\le C_1 2^{-j\delta}|2^jB|^{1/q-1/p}.
\end{align}
The support property is immediate, and the second property in \eqref{eq:mol8} is due to the fact that $P_j$ is the orthogonal projection onto $\VV_j$. For the third property in \eqref{eq:mol8} one notes that by \eqref{eq:mol2} and H\"older's inequality
\begin{align*}
|P_j|
\le \sum_{|\alpha|\le \floor{\om}} \big|\ip{m_j,\uu_{j,\alpha}}_j\big| |\uu_{j,\alpha}|
\lesi \Vert m_j \Vert_{L^q} |U_j(B)|^{-1/q}.
\end{align*}
Then $\Vert P_j\Vert_{L^q}\lesi \Vert m_j\Vert_{L^q}$, concluding the estimate.

Next, for $a_{j,\alpha}$ we have, for $|\beta|\le \floor{\om}$,
\begin{align}\label{eq:mol9}
\supp a_{j,\alpha} \subset 2^{j+1}B, &&
 \int a_{j,\alpha}(x)(x-x_B)^\beta dx=0, &&
 \Vert a_{j,\alpha}\Vert_{L^q}\le C_2 2^{-j\delta}|2^jB|^{1/q-1/p}.
\end{align}
Again the support is immediate. For the  second property in \eqref{eq:mol9}) we see that 
\begin{align*}
\int a_{j,\alpha}(x)(x-x_B)^\beta dx 
&= \N_{j+1,\alpha} \big[ \ip{\nu_{j+1,\alpha},(\cdot-x_B)^\beta}_{j+1} - \ip{\nu_{j,\alpha},(\cdot-x_B)^\beta}_j\big], 
\end{align*}
and each inner product now vanishes because of the orthogonality of the dual basis \eqref{eq:mol3}.  For the third property in \eqref{eq:mol9},  we  have by H\"older's inequality,
\begin{align*}
|\N_{j,\alpha}|
\le \sum_{k\ge j} \Vert m\Vert_{L^q(U_k(B))} \big\Vert (\cdot-x_B)^\alpha \big\Vert_{L^{q'}(U_k(B))}
\lesi \sum_{k\ge j} 2^{-k\delta} (2^kr_B)^{|\alpha|} |2^k B|^{1-1/p},
\end{align*}
and consequently,
\begin{align}\label{eq:mol10} 
|\N_{j,\alpha}|
\lesi 2^{-j\delta} (2^jr_B)^{|\alpha|} |2^jB|^{1-1/p} \sum_{k\ge j} 2^{-(k-j)[\delta-|\alpha|+n(1/p-1)]}
\lesi 2^{-j\delta} (2^jr_B)^{|\alpha|} |2^jB|^{1-1/p}.
\end{align}
Note that our hypothesis on $\delta$ ensures the sum converges, since it holds that $\delta-|\alpha| -n(1/p-1)\ge \delta-\floor{\om}-n(1/p-1)>0.$

It follows  readily then that
\begin{align*}
\Vert a_{j,\alpha}\Vert_{L^q}
\le |\N_{j+1,\alpha}| \Big(\f{\Vert\nu_{j+1,\alpha}\Vert_{L^q}}{|U_{j+1}(B)|}+\f{\Vert\nu_{j,\alpha}\Vert_{L^q}}{|U_{j}(B)|}\Big)
\lesi |\N_{j+1,\alpha}| (2^jr_B)^{-|\alpha|}|2^jB|^{1/q-1} 
\end{align*}
which yields the third estimate in \eqref{eq:mol9}, in view of \eqref{eq:mol10}.

Finally, for $a_\alpha$ we have, for every $1\le q\le \infty$,
\begin{align}
&\supp a_{\alpha} \subset B, \label{eq:mol11a} \\
 &\int a_{\alpha}(x)(x-x_B)^\beta dx=\left
\lbrace 
	\begin{array}{ll}
	0 \qquad &\text{if}\quad \beta\ne\alpha\\
	\N_{0,\alpha} \qquad &\text{if}\quad \beta=\alpha,
	\end{array}
\right. \label{eq:mol11b} \\
 &\Big|\int a_{\alpha}(x)(x-x_B)^\beta dx\Big|\le \left
\lbrace 
	\begin{array}{ll}
	0 \qquad &\text{if}\quad \beta\ne\alpha\\
	|B|^{1-\f{1}{p}} \Big(\f{r_B}{\cro_B}\Big)^{\om-|\alpha|}r_B^{|\alpha|} \qquad &\text{if}\quad \beta=\alpha,
	\end{array}
\right. \label{eq:mol11c} \\
&\Vert a_{\alpha}\Vert_{L^q}\le C_3 |B|^{1/q-1/p}. \label{eq:mol11d} 
\end{align}
The support \eqref{eq:mol11a} follows since $\supp a_\alpha=\supp \nu_{0,\alpha} \subset B$. We also have 
$$ \int a_\alpha(x)(x-x_B)^\beta dx = \N_{0,\alpha}\fint_B \nu_{0,\alpha}(x)(x-x_B)^\beta dx = \N_{0,\alpha} \ip{\nu_{0,\alpha},(\cdot-x)^\beta}_0. $$
and the orthogonality property \eqref{eq:mol11b} now follows from the dual basis property \eqref{eq:mol3}. Inequality \eqref{eq:mol11c} follows from orthogonality of $a_\alpha$ and the estimate $|\N_{0,\alpha}|$ using Definition \ref{def:molecule} (iv). For \eqref{eq:mol11d} we have by \eqref{eq:mol5} and Definition \ref{def:molecule} (iv)
\begin{align*}
|a_\alpha(x)| 
\le \f{|\nu_{0,\alpha}|}{|B|}|\N_{0,\alpha}| 
\lesi \f{r_B^{-|\alpha|}}{|B|} |B|^{1-1/p} r_B^{|\alpha|} \Big(\f{r_B}{\cro_B}\Big)^{\om-|\alpha|} 
= |B|^{-1/p} \Big(\f{r_B}{\cro_B}\Big)^{\om-|\alpha|}.
\end{align*}
Then for any $q\ge 1$ we have, since $r_B\le \f{1}{2}\cro_B$, 
$$ \Vert a_\alpha\Vert_{L^q} \le \Vert a_\alpha\Vert_{L^\infty} |B|^{1/q} \lesi |B|^{1/q-1/p} \Big(\f{r_B}{\cro_B}\Big)^{\om-|\alpha|}\le |B|^{1/q-1/p}$$
as required.

Returning now to the estimate \eqref{eq:mol0} we apply Lemma \ref{lem:AE} to scaled multiples of $a_j, a_{j,\alpha}$ and $a_\alpha$ with $\wt{p}=p$. In particular on applying Lemma \ref{lem:AE} (b) to $b=C_12^{j\delta}a_j$ with ball $2^jB$ and to $b=C_2 2^{j\delta} a_{j,\alpha}$ with ball $2^{j+1}B$, we obtain
$$ \Vert \MM_\LL(a_j)\Vert_{L^p}+\Vert \MM_\LL(a_{j,\alpha})\Vert_{L^p}\lesi 2^{-j\delta}.$$
In a similar fashion, Lemma \ref{lem:AE} (c) applied to $b=C_3 a_\alpha$ with ball $B$ gives
$$ \Vert \MM_\LL(a_\alpha)\Vert_{L^p}\lesi 1.$$
Gathering these estimates together we see that
\begin{align*}
\Vert \MM_\LL(m)\Vert_{L^p}
\le \sum_{j\ge 0}\Vert \MM_\LL(a_j)\Vert_{L^p}
+\sum_{|\alpha|\le \floor{\om}}\sum_{j\ge 0} \Vert \MM_\LL(a_{j,\alpha})\Vert_{L^p}
+\sum_{|\alpha|\le \floor{\om}} \Vert \MM_\LL(a_\alpha)\Vert_{L^p}
<\infty,
\end{align*}
which gives \eqref{eq:mol0}, allowing us to conclude that $m\in h^p_\LL$.
\end{proof}

\subsection{Lipschitz spaces and duality}
In this section we define the Lipschitz spaces associated to the Hermite operator, and prove that they are the dual of the Hardy spaces.

We say that $\vph$ is an \emph{admissible function} if  $\varphi\in C^\infty(\RR_+)$ and
\begin{align*}
	\supp\varphi\subset \big[\tfrac{1}{4},1\big] \ \ \ \text{and} \ \ \  |\varphi|>c>0 \quad\text{on}\quad \big[2^{-7/4},2^{1/4}\big]
\end{align*}
for some $c>0$. Given an admissible function $\vph$, we set $\varphi_j(\lambda):=\varphi(2^{-j}\lambda)$ if $j\in \ZZ$  and call the resulting collection $\{\varphi_j\}_{j\in\ZZ}$ an \emph{admissible system}.

Since the Hermite functions $h_\xi$ with $\xi\in \NN_0^n$ are members of $\sz(\RR^n)$, then for an admissible system $\{\varphi_j\}_{j\in\NN_0}$ we may define the operators $\varphi_j(\sqrt{\LL})$ on $\sz'(\RR^n)$ by
$$ \varphi_j(\sqrt{\LL})f(x)=\sum_{\xi\in\NN_0^n}\vph_j(\xi)  \ip{f,h_\xi}h_\xi(x)\qquad \forall f\in\sz'(\RR^n), x\in \RR^n,$$
where  $ \ip{f,\phi}=f(\phi)$ for $f\in\sz'(\RR^n)$ and with the understanding that $\vph_j(\xi):=\vph_j(\sqrt{2|\xi|+n})$.  
Note that since the Hermite operator $\LL$ has a spectral gap, then  $\varphi_j(\sqrt{\LL})f \equiv 0$ for  every $j\in \ZZ^-$ and $f\in \sz'(\RR^n)$.

Denote by $\{\I_j\}_{j\in\NN_0}$ the following subsets of $\NN_0^n$: 
$\I_j=\big\{\xi\in\NN_0^n: \f{1}{2}4^{j-2}-\frac{n}{2}\le |\xi|\le  \f{1}{2}4^{j}-\frac{n}{2} \big\}$ for $j\in \NN$.
In view of the support of $\varphi_j,$ it follows that
$$ \varphi_j(\sqrt{\LL})f(x)=\sum_{\xi\in \I_j}\vph_j(\xi) \ip{f,h_\xi}h_\xi(x).$$
with  kernels given by
$$ \vp_j(\sqrt{\LL})(x,y) 
= \sum_{\xi\in\I_j} \vph_j(\xi)  h_\mu(x) h_\mu(y).
$$
To define our spaces we also need a certain class of functions with admissible growth; we say that a locally integrable function belongs to the class $\F(\RR^n)$ if  there exists $\gamma>0$ such that
\[
\int_{\RR^n} \f{|f(x)|^2}{(1+|x|)^{n+\gamma}}dx<\vc.
\] 
We may now define the Lipschitz spaces associated with Hermite expansions as follows. 
\begin{Definition}[Hermite Lipschitz spaces]\label{def:Lip-spaces}
For $s>0$ we say that a function $f\in \mathcal F(\RR^n)$ belongs to the Hermite Lipschitz space $\Lip_\LL^s(\RR^n)$ if 
\begin{equation}\label{eq-equivalentnorms for Lip spaces}
\Vert f\Vert_{\Lip_\LL^s} 
:= \sup_{j\in \ZZ} 2^{js}\big\Vert\vph_j(\sqrt{\LL})f\big\Vert_{L^\infty}<\infty.
\end{equation}
\end{Definition}
\noindent For  $s=0$ we define $\Lip^0_\LL$ to be the space of bounded mean oscillation  associated with the Hermite operator introduced in \cite{DGMTZ}. That is, 
$$ \Lip_\LL^0(\RR^n):=\bmo_\LL(\RR^n) = \big\{f\in L^1_{loc}(\RR^n): \Vert f\Vert_{\bmo_\LL}<\infty\big\}$$
where  $\Vert f\Vert_{\bmo_\LL}$ is the infimum of all $C>0$ such that 
$$ \sup_B \fint_B|f(x)-f_B|\,dx \le C \quad\text{and}\quad \sup_{r_B\ge \cro(x_B)} \fint_B|f(x)|\,dx \le C.$$

It is important to note that since $\varphi_j(\LL)f=0$ for  $j\in \ZZ^-$ and $f\in \sz'(\RR^n)$, then \eqref{eq-equivalentnorms for Lip spaces} is equivalent to
\begin{align}\label{eq:Lip space}
\Vert f\Vert_{\Lip_\LL^s} 
= \sup_{j\ge 0} 2^{js}\big\Vert\vph_j(\sqrt{\LL})f\big\Vert_{L^\infty}<\infty.
\end{align}
From \eqref{eq:Lip space} one can see that $\Lip_\LL^s = B^{\infty,\infty}_s(\LL)$ for $s>0$, where $B^{p,q}_s(\LL)$ are the Besov spaces associated to the Hermite operator (see \cite{BD,LN21,BLL21,PX}).
With this in mind we have the following characterization: for every $\ell \in \ZZ^+$ with $\ell>s/2$, 
\begin{equation}\label{eq1-equivalent norms}
	\|f\|_{\Lip_\LL^{s}} \sim \sup_{t>0} t^{-s/2}\|(t\LL)^\ell e^{-t\LL}f\|_{L^\vc}.
\end{equation}
For a proof of  \eqref{eq1-equivalent norms} see \cite[Corollary 3.8]{BBD}, and 
for other characterizations and related notions see \cite{DT19,DT20,DT21}.

The main result of this section shows that the Lipschtiz spaces as defined above are the dual of the Hermite--Hardy spaces. 

\begin{Theorem}[Duality]\label{thm-duality}
For $0<p\le 1$, we have 
\begin{align*}
(h^p_\LL)^* = \left\lbrace
\begin{array}{lll}
&\bmo_\LL, &\qquad p=1 \smallskip\\ 
&\Lip_\LL^{n(\f{1}{p}-1)}, &\qquad 0<p<1.
\end{array}\right.
\end{align*}	
\end{Theorem}
Note that the case $p=1$ was proved in \cite{DGMTZ}, and thus we only give here the proof of  $p<1$. 
In order to prove this result, we need an alternative characterization for the Lipschitz space $\Lip_\LL^{s}$ that is related to the classical Campanato spaces. Such spaces are defined as follows. 
\begin{Definition}[Hermite--Campanato spaces]
	Let $\alpha>0, N\ge 1+\floor{\f{n\alpha}{2}}$. We say that $f\in \mathcal \F(\RR^n)$ belongs to $\mathfrak L_{\mathcal L}(\alpha,N)$ if 
	\[
	\|f\|_{\mathfrak L_{\mathcal L}(\alpha,N)}:=\sup_{B: \ {\rm balls}} \f{1}{|B|^\alpha}\Big(\fint_B |(I-e^{r_B^2\mathcal L})^N f(x)|^2 dx\Big)^{1/2}<\vc,
	\]
	where the supremum is taken over all balls $B$ in $\RR^n$.
\end{Definition}

In the following, we see that the  Hermite--Lipschitz and Hermite--Campanato spaces  are equivalent for suitable indices.
\begin{Proposition}\label{prop:campanato-lip}
	For $s>0$ and $N\ge 1+\floor{\f{s}{2}}$, we have $\Lip_\LL^{s}=\mathfrak L_{\mathcal L}\big(\f{s}{n},N\big)$
	with equivalent norms.
\end{Proposition}	
\begin{proof}[Proof of Proposition \ref{prop:campanato-lip}]
We note that our proof will make use of the characterisation given in \eqref{eq:Lip space}. 
	We first prove that 
\begin{align}\label{eq:lip to campanato}
	\|f\|_{\mathfrak L_{\mathcal L}(\f{s}{n},N)}\lesi \|f\|_{\Lip_\LL^{s}}.
\end{align}

    To begin, we write, for each ball $B$,
    \[
    \begin{aligned}
    	\fint_B |(I-e^{r_B^2\mathcal L})^N f(x)|^2 dx&=\fint_B \Big|\int_{[0,r_B^2]^N}\LL^N e^{-|\vec u|\LL}f(x) d\vec{u}\Big|^2 dx,
    \end{aligned}
    \]
    where $d \vec u = du_1\ldots du_N$ and $|\vec u| = u_1 + \ldots + u_N.$

 Now,  since $N>s/2$, we may use \eqref{eq1-equivalent norms} along with the geometric-harmonic mean inequality to obtain
    \[
    |\LL^N e^{-|\vec u|\LL}f(x)|\lesi |\vec u|^{s/2-N}\|f\|_{\Lip_\LL^{s}} \lesi \prod_{i=1}^N u_i^{\f{s}{2N}-1}\|f\|_{\Lip_\LL^{s}}.
    \]
Then it follows readily that
    \[
    \begin{aligned}
    	\fint_B |(I-e^{r_B^2\mathcal L})^N f(x)|^2 dx
	\lesi \fint_B \Big|\int_{[0,r_B^2]^N} \prod_{i=1}^N u_i^{\f{s}{2N}-1}\|f\|_{\Lip_\LL^{s}} d\vec{u}\,\Big|^2 dx
    	\lesi r_B^{2s}\|f\|_{\Lip_\LL^{s}}^2,
    \end{aligned}
    \]
    which implies \eqref{eq:lip to campanato}.

We now prove the reverse inequality, 
\begin{align}\label{eq:lip to campanato}
    \|f\|_{\Lip_\LL^{s}}\lesi \|f\|_{\mathfrak L_{\mathcal L}(\f{s}{n},N)}.
\end{align}
Suppose that $f\in \mathfrak L_{\mathcal L}\big(\f{s}{n},N\big)$ with $N\ge 1 +\floor{\f{s}{2}}$. Then we have
 for each $\ell \in\ZZ^+$,
    \[
    \begin{aligned}
    	t^{-s/2}(t\LL)^\ell e^{-t\LL}f &= c\int_0^\vc t^{-s/2}(t\LL)^\ell e^{-t\LL}(u\LL)^\ell e^{-u\LL}(I-e^{-u\LL})^N f \f{du}{u}\\
    	&=c\int_0^\vc t^{-s/2}(ut)^\ell\LL^{2\ell} e^{-(t+u)\LL} (I-e^{-u\LL})^N f \f{du}{u}.
    \end{aligned}
    \]
This in turn implies, from the bounds on the heat kernel of $(t+u)^{2\ell}\LL^{2\ell} e^{-(t+u)\LL}$ (see Lemma \ref{lem:HK}), that
    \[
    \begin{aligned}
    	|t^{-s/2}&(t\LL)^\ell e^{-t\LL}f(x)|\\
    	&\lesi \int_0^\vc t^{-s/2}\f{(ut)^\ell}{(t+u)^{2\ell}} \int_{\RR^n} \f{e^{-\f{|x-y|^2}{c(t+u)}}}{(t+u)^{n/2}} |(I-e^{-u\LL})^N f(y)|dy \f{du}{u}\\
    	&\lesi \sum_{k=0}^\vc \int_0^\vc t^{-s/2}\f{(ut)^\ell}{(t+u)^{2\ell}} \int_{U_k(B(x,\sqrt u))} \f{e^{-\f{|x-y|^2}{c(t+u)}}}{(t+u)^{n/2}} |(I-e^{-u\LL})^N f(y)|dy \f{du}{u}.
    \end{aligned}
    \]

    Now, using the Cauchy--Schwarz inequality, we have for each $k\in \NN_0$,
    \[
    \begin{aligned}
    	\int_{U_k(B(x,\sqrt u))} &\f{e^{-\f{|x-y|^2}{c(t+u)}}}{(t+u)^{n/2}} |(I-e^{-u\LL})^N f(y)|dy\\
    	&\le \Big(\int_{U_k(B(x,\sqrt u))} \f{e^{-2\f{|x-y|^2}{c(t+u)}}}{(t+u)^{n/2}}dy\Big)^{\half}
    	\times \Big(\int_{U_k(B(x,\sqrt u))}|(I-e^{-u\LL})^N f(y)|^2dy\Big)^{\half}.
    \end{aligned}
    \]
    It is easy to see that the first factor in the previous inequality satisfies
    \[
    \begin{aligned}
    	\Big(\int_{U_k(B(x,\sqrt u))}\f{e^{-2\f{|x-y|^2}{c(t+u)}}}{(t+u)^{n/2}}dy\Big)^{\half}&\lesi \f{|U_k(B(x,\sqrt u))|^{\half}}{(t+u)^{n/2}}e^{-\f{2^{2k}u}{c'(t+u)}}
    	\sim \Big(\f{2^k\sqrt u}{t+u}\Big)^{\f{n}{2}}e^{-\f{2^{2k}u}{c'(t+u)}}
    	    	\sim \f{e^{-\f{2^{2k}u}{2c'(t+u)}}}{(t+u)^{n/4}}.
    \end{aligned}
    \]
For the second factor, we decompose the annulus $U_k(B(x,\sqrt u))$ via
    \[
    U_k(B(x,\sqrt u))\subset \bigcup_{j=1}^{K_0} B_{kj},
    \]
    where  $B_{kj}$ are balls with $r_{B_{kj}}=\sqrt u$ for all $k,j$, and $K_0\sim 2^{kn}$ is a universal constant. 
    Then, since $f\in \mathfrak{L}_\LL\big(\f{s}{n},N\big)$ with $N\ge 1+\floor{\f{s}{2}}$ we see that
    \begin{align*}
    \int_{B_{kj}}|(I-e^{-u\LL})^N f(y)|^2dy \lesi |B_{kj}|^{1+2s/n}\|f\|^2_{\mathfrak L_{\mathcal L}(\f{s}{n},N)}.
    \end{align*}
    It follows then that
    \[
    \begin{aligned}
    	\int_{U_k(B(x,\sqrt u))}|(I-e^{-u\LL})^N f(y)|^2dy
    	\lesi \sum_{j=1}^{K_0} |B_{kj}|^{1+2s/n}\|f\|^2_{\mathfrak L_{\mathcal L}(\f{s}{n},N)}
    	\sim 2^{kn} u^{n/2+s}\|f\|^2_{\mathfrak L_{\mathcal L}(\f{s}{n},N)}.
    \end{aligned}
    \]
   Combining the above estimates we arrive at the following: for and $\ell \in \ZZ^+$ and a.e. $x\in \RR^n$,
    \[
    \begin{aligned}
    	|t^{-s/2}&(t\LL)^\ell e^{-t\LL}f(x)|
    	\lesi \sum_{k=0}^\vc 2^{kn/2}\|f\|_{\mathfrak L_{\mathcal L}(\f{s}{n},N)}\int_0^\vc t^{-s/2}\f{(ut)^\ell}{(t+u)^{2\ell}}\f{e^{-\f{2^{2k}u}{2c'(t+u)}}}{(t+u)^{n/4}} u^{n/4+s/2} \f{du}{u}.
    \end{aligned}
    \]
 Now take  $\ell>s/2$; a simple calculation shows that 
    \[
    \begin{aligned}
    	\int_0^\vc t^{-\f{s}{2}}\f{(ut)^\ell}{(t+u)^{2\ell}}\f{e^{-\f{2^{2k}u}{2c'(t+u)}}}{(t+u)^{\f{n}{4}}} u^{\f{n}{4}+\f{s}{2}} \f{du}{u}
	\lesi 2^{-(\f{n}{2}+s)k} t^{\ell-\f{s}{2}} \int_0^\infty u^{\ell-1} (t+u)^{\f{s}{2}-2\ell} du
    	 \lesi 2^{-(\f{n}{2}+s)k}.
    \end{aligned}
    \]
    Therefore, for any $\ell>s/2$,
    \[
    \begin{aligned}
    	|t^{-s/2}(t\LL)^\ell e^{-t\LL}f(x)|
    	\lesi \sum_{k=0}^\vc 2^{-sk}\|f\|_{\mathfrak L_{\mathcal L}(\f{s}{n},m)}
    	\lesi \|f\|_{\mathfrak L_{\mathcal L}(\f{s}{n},m)}.
    \end{aligned}
    \]
    This, along with \eqref{eq1-equivalent norms}, implies \eqref{eq:lip to campanato},
   completing our proof of Proposition \ref{prop:campanato-lip}. 
\end{proof}

\begin{proof}[Proof of Theorem \ref{thm-duality}:]  
	As mentioned earlier, we only need to establish the case $p<1$. For such $p$, we have from Theorem 1.3 in \cite{SY} and Theorem 4.1 in \cite{Yan}, that
	\[
	(h^p_\LL)^* = \mathfrak L_{\mathcal L}\big(\tfrac{1}{p}-1,N\big),
	\]
	for every $N\ge 1+\big\lfloor\f{n}{2}(\f{1}{p}-1)\big\rfloor$.
	This, in conjunction with Proposition \ref{prop:campanato-lip}, yields Theorem \ref{thm-duality}.
\end{proof}

\bigskip

Finally, we close this section with a lemma that describes a type of function from $\Lip_\LL^s$ that will be important in  the proofs of our main results. 
\begin{Lemma}\label{lem:lip eg}
For any $x_0\in\RR^n$ and $\alpha\in\NN_0^n$ we set 
$$ g_{x_0,\alpha}(x):=(x-x_0)^\alpha \chi(x)$$
where $\chi\in C^\infty_0(\RR^n)$ with $\chi\equiv 1$ on $B_0:=B(x_0,\cro(x_0))$, $\chi\equiv 0$ on $(2B_0)^c$ and $\Vert \partial^\gamma \chi\Vert_{L^\infty}\lesi \cro(x_0)^{-|\gamma|}$ for $\gamma\in\NN_0^n$. 

Then $g_{x_0,\alpha}\in \Lip_\LL^s$  for any $s\ge 0$ and there exists $C_{n,s,\alpha}$ depending on $s$, $n$ and $\alpha$ such that
$$ \Vert g_{x_0,\alpha}\Vert_{\Lip_\LL^s}\le C_{n,s,\alpha} \,\cro(x_0)^{|\alpha|-s}.$$
\end{Lemma}
\begin{proof}[Proof of Lemma \ref{lem:lip eg}]
We first show that $g_{x_0,\alpha}\in \Lip_\LL^0(\RR^n)\equiv \bmo_\LL(\RR^n)$. Observe that for any ball $B$, 
\begin{align*}
\fint_B |g_{x_0,\alpha}|\le \sup_{x\in 2B_0} |(x-x_0)^\alpha| \fint_B |\chi|  \le C_\alpha \,\cro(x_0)^{|\alpha|}.\end{align*}
It follows readily that 
\begin{align*}
\Vert g_{x_0,\alpha} \Vert_{\bmo_\LL} \le 3\sup_{B\subset\RR^n} \fint_B |g_{x_0,\alpha}| \le 3C_\alpha  \,\cro(x_0)^{|\alpha|}
\end{align*}
as required. 

We next consider $s>0$. For each $j\ge 0$ we have, by the Cauchy--Schwarz inequality,
\begin{align*}
\big|\vph_j(\sqrt{\LL})g_{x_0,\alpha}(x)\big|
&=\Big| \sum_{\xi\in\I_j} h_\xi(x)\int \chi(y)(y-x_0)^\alpha h_\xi(y)\,dy\Big|\\
&\le \Big(\sum_{\xi\in\I_j} h_\xi(x)^2\Big)^{\half} \Big(\sum_{\xi\in \I_j} \Big|\int \chi(y)(y-x_0)^\alpha h_\xi(y)\,dy\Big|^2\Big)^\half.
\end{align*}
for every $x\in\RR^n$. Now firstly in view of \eqref{eq:QQ est} we observe that
\begin{align}\label{eq:htoQQ} \Big(\sum_{\xi\in\I_j} h_\xi(x)^2\Big)^{\half} \le \QQ_{4^j} (x,x)^\half \lesi 2^{jn/2},\end{align}
and secondly from \eqref{eq:CN0} we have, for any $N\ge 0$,
$$ \Big(\sum_{\xi\in \I_j} \Big|\int \chi(y)(y-x_0)^\alpha h_\xi(y)\,dy\Big|^2\Big)^\half 
\lesi \cro(x_0)^{|\alpha|+n-2N}2^{-j(2N-n/2)}.
$$
Combining the preceding three calculations yields
\begin{align*}
\big|\vph_j(\sqrt{\LL})g_{x_0,\alpha}(x)\big| 
\lesi \cro(x_0)^{|\alpha|+n-2N}2^{-j(2N-n)}, \qquad x\in\RR^n.
\end{align*}
Now any $s>0$, if $1+|x_0|\ge 2^j$ then we choose $N=0$; otherwise we set $N=\floor{\f{n+s}{2}}+1$. In either case we obtain
\begin{align*}
2^{js}\big\Vert\vph_j(\sqrt{\LL})g_{x_0,\alpha}\big\Vert_{L^\infty}
\lesi \cro(x_0)^{|\alpha|-s}\Big(\f{1+|x_0|}{2^j}\Big)^{2N-n-s}
\lesi \cro(x_0)^{|\alpha|-s},
\end{align*}
where the constant depends on $N = N(n,s)$. Since our estimate above is independent of $j$ then we have
\begin{align*}
\Vert g_{x_0,\alpha}\Vert_{\Lip_\LL^s} 
= \sup_{j\ge 0} 2^{js}\big\Vert\vph_j(\sqrt{\LL})g_{x_0,\alpha}\big\Vert_{L^\infty}
\lesi \cro(x_0)^{|\alpha|-s},
\end{align*}
which completes the proof of Lemma \ref{lem:lip eg}.
\end{proof}

\section{Proofs of the main results}\label{sec:proofs}

\subsection{Proof of Theorem \ref{thm:hardy}}
{\it Part (a):}
Fix $p \in (\f{n}{n+\floor{\om}+\om^*\land\ve}, 1]$. We shall prove that $T$ maps  $(p,2,\floor{\om})$-atoms into multiples of  $(p,2,\delta,\floor{\om}+\mu)$-molecules with 
\begin{align*}
\delta = \floor{\om}+\ve -n\big(\tfrac{1}{p}-1\big)>0\qquad\text{and}\qquad
 \mu=\min\{\om^*,\ve\}.
\end{align*}
Let $a$ be a $(p,2,\floor{\om})$-atom associated to some ball $B$. We first check that $Ta$ satisfies Definition \ref{def:molecule} (ii). For $j=0,1,2$, by the $L^2$ boundedness of $T$ we have
$$ \Vert Ta \Vert_{L^2(U_j(B))} \lesi \Vert a\Vert_{L^2} \le |B|^{1/2-1/p}.$$
For $j\ge 3$ we consider two cases.

\underline{Case 1:} $r_B\ge \f{1}{8}\cro_B$. 
From Minkowski's inequality and  Definition \ref{def:HCZO} (ii) we have
\begin{align*}
\Vert Ta\Vert_{L^2(U_j(B))}
&\le \int_B \Big(\int_{U_j(B)} |K(x,y)|^2 dx\Big)^{\half}|a(y)|\,dy\\
&\lesi \int_B \Big(\int_{U_j(B)} |x-y|^{-2n} \big(1+\tfrac{|x-y|}{\cro(y)}\big)^{-2(\floor{\om}+\ve)} dx\Big)^{\half}|a(y)|\,dy
\end{align*}
Since $\f{1}{8}\cro_B\le r_B\le \f{1}{2}\cro_B$  then for $y\in B$ we have $\cro(y)\sim \cro_B\sim r_B$.
Observe also that for $j\ge 3$ the fact that $x\in U_j(B)$ and $y\in B$ implies $|x-y|\gtrsim 2^j r_B$. These two facts give 
$$ |x-y|^{-n} \big(1+\tfrac{|x-y|}{\cro(y)}\big)^{-(\floor{\om}+\ve)} 
\lesi 2^{-j(n+\floor{\om}+\ve)} r_B^{-n}\Big(\f{\cro(y)}{r_B}\Big)^{\floor{\om}+\ve}\sim 2^{-j(\floor{\om}+\ve)}|2^jB|^{-1}.$$
Inserting this estimate into the previous calculation gives
\begin{align*}
\Vert Ta\Vert_{L^2(U_j(B))}
\lesi 2^{-j(\floor{\om}+\ve)}|2^jB|^{-\half}\Vert a\Vert_{L^1} 
\lesi 2^{-j(\floor{\om}+\ve - n(1/p-1)} |2^j B|^{\half-\f{1}{p}}
\le 2^{-j\delta} |2^j B|^{\half-\f{1}{p}}.
\end{align*}

\underline{Case 2:} $r_B<\f{1}{8}\cro_B$. 
The cancellation of $a$, along with Taylor's theorem, the mean value theorem, and  Minkowski's inequality, gives
\begin{align*}
\Vert Ta \Vert_{L^2(U_j(B))}
&=\Big\{\int_{U_j(B)} \Big(\int_B\big[K(x,y)-\sum_{|\gamma|\le \floor{\om}}\f{1}{\gamma!} \partial_2^\gamma K(x,x_B) (y-x_B)^\gamma\big] a(y)\,dy\Big)^2dx\Big\}^{\half} \\
&=\Big\{\int_{U_j(B)} \Big(\sum_{|\gamma|=\floor{\om}}\f{1}{\gamma!}\int_B \big[\partial_2^\gamma K(x,\wt{y})-\partial_2^\gamma K(x,x_B)\big](y-x_B)^\gamma a(y)\,dy\Big\}^{\half} \\
&\le \sum_{|\gamma|=\floor{\om}} \f{1}{\gamma!} \int_B\Big(\int_{U_j(B)} \big|\partial_2^\gamma K(x,\wt{y})-\partial_2^\gamma K(x,x_B)\big|^2 |y-x_B|^{2|\gamma|}dx\Big)^{\half}|a(y)|\,dx
\end{align*}
where $\wt{y}$ is some point on the line segment joining $y$ and $x_B$. Since $x\in U_j(B)$ and $\wt{y}\in B$, then for $j\ge 3$ we have
\begin{align*}
|x-x_B|>2r_B\ge 2|\wt{y}-x_B| 
\end{align*}
and so we may apply the estimate in Definition \ref{def:HCZO} (iv) to obtain
\begin{align*}
\Vert Ta\Vert_{L^2(U_j(B))}
\lesi \int_B \Big(\int_{U_j(B)} \f{|y-x_B|^{2(\floor{\om}+\ve)}}{|x-x_B|^{2(n+\floor{\om}+\ve)}} \,dx\Big)^{\half}|a(y)|\,dy 
&\lesi r_B^{-n} 2^{-j(n+\floor{\om}+\ve)}|2^j B|^{\f{1}{2}}\Vert a\Vert_{L^1}.
\end{align*}
From the fact that $\delta=\floor{\om} +\ve -n(\f{1}{p}-1)>0$, it follows that
\begin{align*}
\Vert Ta\Vert_{L^2(U_j(B))}
\lesi 2^{-j(\floor{\om}+\ve)} |2^jB|^{-\f{1}{2}}|B|^{1-\f{1}{p}}
\le 2^{-j\delta}|2^jB|^{1/2-1/p}.
\end{align*}
With this latter estimate obtained, we have thus shown that $Ta$ satisfies Definition \ref{def:molecule}~(ii).

We now show that $Ta$ satisfies Definition \ref{def:molecule} (iii). Note that is this is only required when $r_B<\f{1}{8}\cro_B$. Accordingly, let $\chi_B\in C^\infty_0(\RR^n)$ with $\chi_B\equiv 1$ on $B(x_B,\cro_B)$ and $\chi_B\equiv 0$ on $B(x_B, 2\cro_B)^c$. Then for each $|\alpha|\le \floor{\om}$,
\begin{align*}
\Big|\int (x-x_B)^\alpha Ta(x)\,dx\Big|
\le \Big|\ip{Ta,(\cdot-x_B)^\alpha \chi_B}\Big|+\Big|\ip{Ta,(\cdot-x_B)^\alpha (1-\chi_B)}\Big|
=:I+II.
\end{align*}
Now by duality (Theorem \ref{thm-duality}), Lemma \ref{lem:AE} (b) (see Remark \ref{rem:AE}) with $\wt{p}=\f{n}{n+\om}$, and our hypothesis \eqref{eq:hardycond}, we have
\begin{align*}
I
=\big|\ip{a,T^*\big[(\cdot-x_B)^\alpha \chi_B\big]}\big|
\le \Vert a\Vert_{h_\LL^{\f{n}{n+\om}}} \big\Vert T^*\big[(\cdot-x_B)^\alpha \chi_B\big\Vert_{\Lip_\LL^{\om}} 
\lesi |B|^{1+\f{\om}{n}-\f{1}{p}} \cro_B^{|\alpha|-\om}.
\end{align*}
Since $r_B\le \f{1}{2}\cro_B$ and $\om^*\ge \mu$ then we conclude
\begin{align*}
I\lesi |B|^{1-\f{1}{p}}r_B^{|\alpha|} \Big(\f{r_B}{\cro_B}\Big)^{\om-|\alpha|} \lesi |B|^{1-\f{1}{p}}r_B^{|\alpha|} \Big(\f{r_B}{\cro_B}\Big)^{\floor{\om}+\mu-|\alpha|}.
\end{align*}

For term $II$ we use the cancellation of $a$, Taylor's expansion and the Mean Value Theorem to write 
\begin{align*}
II
&=\Big|\int (1-\chi_B)(x-x_B)^\alpha \int_B\big[K(x,y)-\sum_{|\gamma|\le \floor{\om}}\f{1}{\gamma!} \partial_2^\gamma K(x,x_B) (y-x_B)^\gamma\big] a(y)\,dy\,dx\Big| \\
&\le \int\limits_{B(x_B,\cro_B)^c} |x-x_B|^{|\alpha|} \sum_{|\gamma|=\floor{\om}}\f{1}{\gamma !}\int_{U_j(B)} \big|\partial_2^\gamma K(x,\wt{y})-\partial_2^\gamma K(x,x_B)\big| |y-x_B|^{|\gamma|} |a(y)|\,dy\,dx
\end{align*}
where for each $y\in B$, $\wt{y}$ is a point on the line segment joining $y$ and $x_B$. Since $y\in B$ with $r_B<\f{1}{8}\cro_B$ and $x\in B(x_B,\cro_B)^c$ then 
$$ |x-x_B|\ge \cro_B>8r_B\ge 2|y-x_B|\ge 2|\wt{y}-x_B|,$$
and so we may use the estimate in Definition \ref{def:HCZO} (iv) again to obtain
\begin{align*}
II
&\lesi \int\limits_{B(x_B,\cro_B)^c} |x-x_B|^{|\alpha|} \int_B \f{|\wt{y}-x_B|^\ve}{|x-x_B|^{n+\floor{\om}+\ve}}|y-x_B|^{\floor{\om}}|a(y)|\,dy\,dx \\
&\le r_B^{\floor{\om}+\ve}\,\Vert a\Vert_{L^1}\int\limits_{B(x_B,\cro_B)^c} \f{dx}{|x-x_B|^{n+\floor{\om}+\ve-|\alpha|}}.
\end{align*}
Since $\floor{\om}\ge |\alpha|$ and $\ve >0$ then the integral is bounded by a constant independent of $x_B$. Thus we have
\begin{align*}
II
\lesi |B|^{1-\f{1}{p}}r_B^{\floor{\om}+\ve} \cro_B^{|\alpha|-\floor{\om}-\ve}
=|B|^{1-\f{1}{p}}r_B^{|\alpha|} \Big(\f{r_B}{\cro_B}\Big)^{\floor{\om}+\ve-|\alpha|},
\end{align*}
and since $r_B\le \f{1}{2}\cro_B$ and $\ve\ge \mu$ we conclude
\begin{align*}
II
\lesi |B|^{1-\f{1}{p}}r_B^{|\alpha|} \Big(\f{r_B}{\cro_B}\Big)^{\floor{\om}+\mu-|\alpha|}.
\end{align*}

Thus combining the estimates for $I$ and $II$, and taking into account that $r_B<\f{1}{8}\cro_B$, we have 
\begin{align*}
\Big|\int (x-x_B)^\alpha Ta(x)\,dx\Big|
\lesi |B|^{1-\f{1}{p}}r_B^{|\alpha|} \Big(\f{r_B}{\cro_B}\Big)^{\floor{\om}+\mu-|\alpha|}
\end{align*}
and our proof of (iii) is complete on recalling that $\mu=\min\{\om^*,\ve\}$.

{\it Part (b):}
Since $T$ is bounded on $h^p_\LL(\RR^n)$, then by duality $T^*$ is bounded on $\Lip_\LL^{n(\f{1}{p}-1)}$. Recall by Lemma \ref{lem:lip eg}, that for any $x_0\in\RR^n$ and $\chi$ as defined in the Lemma  the function $g_{x_0,\alpha}:=(x-x_0)^\alpha \chi$ satisfies $g_{x_0,\alpha}\in \Lip_\LL^{n(\f{1}{p}-1)}$  with 
$$ \Vert g_{x_0,\alpha} \Vert_{\Lip_\LL^{n(\f{1}{p}-1)}}\lesi \cro(x_0)^{|\alpha|-n(\f{1}{p}-1)}$$
for any $\alpha\in\NN_0^n$. Combining the above facts we thus have
\begin{align*}
\big\Vert T^*[(\cdot-x_0)^\alpha \chi]\big\Vert_{\Lip^{n(\f{1}{p}-1)}_\LL} 
= \Vert Tg_{x_0,\alpha}\Vert_{\Lip^{n(\f{1}{p}-1)}_\LL} 
\lesi \Vert g_{x_0,\alpha}\Vert_{\Lip^{n(\f{1}{p}-1)}_\LL} 
\lesi \cro(x_0)^{|\alpha|-n(\f{1}{p}-1)}
\end{align*}
which shows that \eqref{eq:hardycond} holds for $\om=n(\f{1}{p}-1)$.

\subsection{Proof of Theorem \ref{thm:lip}}
{\it Part (a).} Firstly observe that if $T\in \HCZO_{(1)}\big(\floor{\om},\ve\big)$ then $T^*\in \HCZO_{(2)}\big(\floor{\om},\ve\big)$. Thus we may apply Theorem \ref{thm:hardy} to $T^*$ to see that $T^*$ is bounded on $h^{\f{n}{n+s}}_\LL(\RR^n)$ for every $s$ satisfying
$$ \f{n}{n+\floor{\om}+\om^*\land\ve} < \f{n}{n+s}\le 1.$$
By duality we then have that $T$ is bounded on $\Lip^s_\LL(\RR^n)$ for every $0\le s<\floor{\om}+\om^*\land \ve$.

{Part (b).} 
Let $g_{x_0,\alpha}$ and $\chi$ be as defined in Lemma \ref{lem:lip eg} for $|\alpha|\le \floor{\om}$. Then the boundedness of $T$ on $\Lip_\LL^s(\RR^n)$ and Lemma \ref{lem:lip eg} implies
$$ \big\Vert T\big[(\cdot-x_0)^\alpha \chi\big]\big\Vert_{\Lip^s_\LL} 
=\Vert Tg_{x_0,\alpha}\Vert_{\Lip^s_\LL}
\lesi \Vert g_{x_0,\alpha}\Vert_{\Lip^s_\LL}
\lesi \cro(x_0)^{|\alpha|-s}.$$
This gives \eqref{eq:lipcond} as required.

\section{Applications}\label{sec:applications}
In this section we give applications of our results to pseudo-multipliers and hermite Riesz transforms.

\subsection{Pseudo-multipliers}\label{sec:pdo}
Let $\sigma:\RR^n\times\NN_0^n\to \mathbb{C}$ be a bounded function. We define the Hermite pseudo-multiplier associated with $\sigma$ by
$$ \sigma(x,\LL) f(x):=\sum_{\xi\in\NN_0^n} \sigma(x,\xi)\ip{f,h_\xi}h_\xi(x)$$
with kernel 
$$K^\sigma(x,y)=\sum_{\xi\in\NN_0^n} \sigma(x,\xi)h_\xi(x)h_\xi(y).$$
Let $\K,\N\in \NN_0\cup\{\infty\}$. We say that the function $\sigma:\RR^n\times \NN_0^n\to \mathbb{C}$ belongs to the class $\SM^{m,\K,\N}_{\rr,\dd}$ if $\sigma(\cdot,\xi)\in C^\N(\RR^n)$ and 
\begin{align*}
|\partial_x^\nu \diff_\xi^\kappa \sigma(x,\xi)| \lesi \ip{\xi}^{\f{m}{2}+\f{\dd}{2}|\nu|-\rr|\kappa|}, \qquad |\kappa|\le\K,\; |\nu|\le \N
\end{align*}
Observe that $\SM^{m,\K,\N}_{\rr,\dd}\subseteq\SM^{m',\K',\N'}_{\rr',\dd'}$ if one of the following holds: $m\le m', \rr'\le \rr, \dd\le\dd', \K'\le\K, \N'\le \N$. That is, the classes $\SM^{m,\K,\N}_{\rr,\dd}$ are increasing with $m,\dd,\K,\N$ and decreasing with $\rr$. 

The operators $\sigma(x,\mathcal L)$ can be viewed as pseudo-differential
	operators associated to the Hermite operator $\mathcal L$. Such operators also occurs as Weyl transforms of radial
	functions on $\mathbb C^n$. For further details we refer to  \cite{Epp, BT, LN21} and the references therein.

\begin{Theorem}[Pseudo-multipliers are Hermite--Calder\'on--Zygmund operators]\label{thm:pdo-hczo}
Fix $M\in\NN_0$ and let $\sigma:\RR^n\times\NN_0^n\to\mathbb{C}$ be a symbol. 
\begin{enumerate}[\upshape(a)]
\item Suppose that $\sigma\in\SM^{0,n+M+1,\N}_{1,\dd}$ for some $0\le \dd<1$ and  with $\N=2\floor{\f{n+M+1}{2(1-\dd)}}$. Then $\sigma(\cdot,\LL)\in \HCZO(M,\ve)$ for any $0<\ve<1$. 
\item Suppose that $\sigma\in \SM^{0,n+M+1,1}_{1,1}$, along with the estimates
\begin{align}\label{eq:conditionCN}
\sup_x\sup_\xi \Big(\fint_{B(x,\cro(x))} 
\big|\cro(y)^{|\gamma|} \partial_y^\gamma \sigma(y,\xi)\big|^2\,dy\Big)^\half <\infty, \qquad |\gamma|\le 2\floor{(n+M)/2}+2.
\end{align} 
 Then $\sigma(\cdot,\LL)\in \HCZO_{(2)}(M,\ve)$ for any $0<\ve<1$. 
\end{enumerate}
\end{Theorem}
We show that our main results allow us to recover the following from \cite{LN21} (see \cite[Corollary 5.9]{LN21}).

\begin{Theorem}[Pseudo-multipliers on Hardy spaces]\label{thm:pdo-hardy}
Let $0<p\le 1$ and $M=\floor{n(\f{1}{p}-~1)}$. Assume that $\sigma:\RR^n\times\NN_0^n\to\mathbb{C}$ satisfies any of the conditions in Theorem \ref{thm:pdo-hardy}. 
Then $\sigma(\cdot,\LL)$ extends to a bounded operator on $h^p_\LL(\RR^n)$.
\end{Theorem}
For comparison we observe that assuming $L^2(\RR^n)$ boundedness, \cite[Theorem 1.4]{BT} shows that $\sigma(\cdot,\LL)$ is bounded on $L^p(\RR^n)$ and is weak $(1,1)$ if $\sigma\in \SM^{0,n+1,1}_{1,0}$. In contrast we show that if $\sigma\in \SM^{0,n+1,n+1}_{1,0}$ then $\sigma(\cdot,\LL)$ maps $h^1_\LL\to h^1_\LL$, and if $\sigma\in \SM^{0,\floor{\f{n}{p}}+1,\floor{\f{n}{p}}+1}_{1,0}$ then $\sigma(\cdot,\LL)$ maps $h^p_\LL\to  h^p_\LL$. The increased number of derivatives in $\N$ is required because we are proving something a little stronger than $L^p$ estimates. 

The reader may observe from the proof below that we actually prove if $M\in \NN_0$, then under either assumption (a) or (b) from Theorem \ref{thm:pdo-hardy}, we have $\sigma(\cdot,\LL)\in \HCZO_{(2)}(M,\ve)$ for any $0<\ve <1$.

For Lipschitz spaces, the pseudo-multipliers considered are known to hold from \cite{LN21} for $s>0$. However the result there does not include  $\bmo_\LL$ (the case $s=0$). The following result addresses this gap and is thus new for $\bmo_\LL$.
\begin{Theorem}[Pseudo-multipliers on Lipschitz spaces]\label{thm:pdo lip}
Assume that $\sigma$ satisfies one of the following conditions.
\begin{enumerate}[\upshape (a)]
\item $\sigma\in \SM^{0,n+1,\N}_{1,\dd}$ for $0\le \dd<1$ and $\N\ge 2\ceil{\f{n+1}{2(1-\dd)}}$.
\item $\sigma\in \SM^{0,n+1,1}_{1,1}$ and \eqref{eq:conditionCN} holds with $M=0$. 
\end{enumerate}
Then $\sigma(\cdot,\LL)$ extends to a bounded operator on $\Lip^s_\LL(\RR^n)$ for $0\le s<1$. 
\end{Theorem}

Note that for the special case $\dd=0$ from the previous two theorems, we have that if $\sigma\in\SM^{0,n+1,n+1}_{1,0}$ then $\sigma(\cdot,\LL)$ is bounded on $h^p_\LL(\RR^n)$ for all $\f{n}{n+1}<p\le $, on $\Lip^s_\LL(\RR^n)$ for all $0\le s<1$, and on all $L^p(\RR^n)$ for $1<p<\infty$.

Before giving the proof of Theorem \ref{thm:pdo-hardy}, we obtain the some required estimates for the kernel of $\sigma(\cdot,\LL)$ which are interesting in their own right. 
Choose an admissible function $\vph$ satisfying
\begin{align}\label{eq:partition}
\sum_{k\ge 0} \vph(2^{-j}\lambda)=1, \qquad \lambda\ge \tfrac{1}{2}.
\end{align}
Then we write
$$ \sigma(\cdot,\LL) = \sum_{j\ge0} \sigma(\cdot,\LL)  \vph_j(\sqrt{\LL}) = \sum_{j\ge 0} \sigma_j(\cdot,\LL),$$
where the operators $\sigma_j(\cdot,\LL) $ are integral operators with kernels given by
$$K_j^\sigma(x,y)=\int K^\sigma(x,z)\vph_j(\sqrt{\LL})(z,y)\,dz.$$
Then we have the following estimates for $K_j(x,y)$.
\begin{Lemma}\label{lem:ddKj}
Let $m\in\RR$, $\N\ge 0$, $\K\ge 1$, $\rr\in[0,1]$ and  $\delta\in[0,1]$. 
Suppose that $\sigma\in \SM^{m,\K,\N}_{\rr,\delta}$. 
Let $\{\vp_j\}_j$ be an admissible system. Then for any $\mu>0$, $|\gamma|\le \N$, $ N\le \K$ and $\eta\in\NN_0^n$, we  have
\begin{align}\label{eq:Kjab}
|x-y|^{N}\big|\partial_x^\gamma\partial_y^\eta K_j^\sigma(x,y)\big|
\lesssim 2^{j(n+m+|\gamma|+|\eta|+N(1-2\rr))}  \big(1+\tfrac{2^{-j}}{\cro(x)}+ \tfrac{2^{-j}}{\cro(y)}\big)^{-\mu}
\end{align}
for all $j\in \NN_0$.
\end{Lemma}
The proof of Lemma \ref{lem:ddKj} is fairly technical and will be given in Appendix \ref{app:technical}.

\subsubsection{Proof of Theorem \ref{thm:pdo-hczo}}

We first note that the  $L^2(\RR^n)$ boundedness of $\sigma(\cdot,\LL)$ in either case follows directly from \cite[Corollary 5.4]{LN21}. Next, the size estimate in Definition \ref{def:HCZO} (ii) (in both cases) is obtained by applying the following result with $\K=n+M+1$. 
\begin{Proposition}\label{prop:kernelV1}
Suppose that  $\mu\ge 0$, $\N\ge 0$ and   $\K\ge 1$. Let $\sigma\in \SM^{0,\K,\N}_{1,1}$ and $K^\sigma(x,y)$ be the kernel of $\sigma(\cdot,\LL)$. Assume further that $\eta\in\NN_0^n$, $0\le|\gamma|\le \N$ and
$$\K> n+|\gamma|+|\eta|+\mu.$$ 
Then the following holds.
\begin{align}\label{eq:kernelV1.1}
\big|\partial_x^\gamma \partial_y^\beta K^\sigma(x,y)\big| \lesssim  |x-y|^{-n-|\gamma|-|\eta|}\Big(1+\tfrac{|x-y|}{\cro(x)}+\tfrac{|x-y|}{\cro(y)}\Big)^{-\mu}
\end{align}
\end{Proposition}
\begin{proof}[Proof of Proposition \ref{prop:kernelV1}]
Set $n_{\gamma,\beta}:=n+|\gamma|+|\eta|>0$. 
We have
\begin{align*}
\big|\partial_y^\gamma \partial_x^\eta K^\sigma(x,y)\big| 
\le \sum_{j\ge0} \big|\partial_y^\gamma \partial_x^\eta K_j^\sigma(x,y)\big| 
=\sum_{2^j \le |x-y|^{-1}} \dots+\sum_{2^j>|x-y|^{-1}} \dots
=:I+II.
\end{align*}

Then by Lemma \ref{lem:ddKj} with $N=0$, we have for any $\mu\ge 0$
\begin{align*}
I
\lesi \sum_{j=0}^{\floor{-\log_2(|x-y|)}} 2^{jn_{\gamma,\eta}} \big(1+\tfrac{2^{-j}}{\cro(x)}+\tfrac{2^{-j}}{\cro(y)}\big)^{-\mu}  
\le \big(1+\tfrac{|x-y|}{\cro(x)}+\tfrac{|x-y|}{\cro(y)}\big)^{-\mu}   \sum_{j=0}^{\floor{-\log_2(|x-y|)}} 2^{jn_{\gamma,\eta}},
\end{align*}
where in the last step we used that $2^{-j}\ge |x-y|$. Since $n_{\gamma,\eta}>0$ then we can compute the sum
\begin{align*}
 \sum_{j=0}^{\floor{-\log_2(|x-y|)}} 2^{jn_{\gamma,\eta}}
 \lesi 2^{n_{\gamma,\eta}(\floor{-\log_2(|x-y|)})}
 \sim |x-y|^{-n_{\gamma,\eta}},
\end{align*}
which allows us to arrive at
\begin{align*}
I \lesi\big(1+\tfrac{|x-y|}{\cro(x)}+\tfrac{|x-y|}{\cro(y)}\big)^{-\mu} |x-y|^{-n_{\gamma,\eta}}.
\end{align*}
For term $II$ we observe that if $2^j>|x-y|^{-1}$ then $2^j|x-y|>1$ and so
$$ 1+\tfrac{|x-y|}{\cro(x)}+\tfrac{|x-y|}{\cro(y)} \le 2^j|x-y|\big(1+\tfrac{2^{-j}}{\cro(x)}+\tfrac{2^{-j}}{\cro(x)}\big)$$
and hence 
$$\big(1+\tfrac{2^{-j}}{\cro(x)}+\tfrac{2^{-j}}{\cro(y)}\big)^{-\mu} \lesssim 2^{j\mu}|x-y|^\mu \big(1+\tfrac{|x-y|}{\cro(x)}+\tfrac{|x-y|}{\cro(y)}\big)^{-\mu}.$$
Then by Lemma \ref{lem:ddKj} with $N=\K$ we have
\begin{align*}
II 
&\lesi |x-y|^{-\K}\sum_{j=\floor{-\log_2(|x-y|)}+1}^\infty 2^{j(n_{\gamma,\beta}-\K)}\big(1+\tfrac{2^{-j}}{\cro(x)}+\tfrac{2^{-j}}{\cro(y)}\big)^{-\mu}\\
&\lesi |x-y|^{\mu-\K}\big(1+\tfrac{|x-y|}{\cro(x)}+\tfrac{|x-y|}{\cro(y)}\big)^{-\mu}\sum_{j=\floor{-\log_2(|x-y|)}+1}^\infty 2^{j(n_{\gamma,\eta}+\mu-\K)}.
\end{align*}
Now since by assumption $\K>n_{\gamma,\beta}+\mu$ then the series converges and can be controlled by
\begin{align*}
\sum_{j=\floor{-\log_2(|x-y|)}+1}^\infty 2^{j(n_{\gamma,\eta}+\mu-\K)}
\lesi 2^{\floor{-\log_2(|x-y|)}(n_{\gamma,\eta}+\mu-\K)}
\sim |x-y|^{-(n_{\gamma,\eta}+\mu-\K)}.
\end{align*}
Inserting this estimate into the calculation for $II$ we obtain
\begin{align*}
II \lesi  \big(1+\tfrac{|x-y|}{\cro(x)}+\tfrac{|x-y|}{\cro(y)}\big)^{-\mu} |x-y|^{-n_{\gamma,\eta}}.
\end{align*}
Combining our estimates for $I$ and $II$ we arrive at \eqref{eq:kernelV1.1} as required. This completes the proof of Proposition \ref{prop:kernelV1}.
\end{proof}

Let us now turn to the regularity estimates in Definition \ref{def:HCZO} (iii)-(iv). These are contained the following estimates. For part (a) we apply Proposition \ref{prop:czk} (b) to obtain both Definition \ref{def:HCZO} (iii) and (iv), noting that $2\floor{(n+M+1)/(2(1-\dd))}\ge M+1$. For part (b) we apply Proposition \ref{prop:czk} (a) to obtain Definition \ref{def:HCZO} (iv).
\begin{Proposition}\label{prop:czk}
Let  $M\in\NN_0$ and $\sigma:\RR^n\times\NN_0^n\to \mathbb{C}$ be a symbol.Consider the following estimates on $K^\sigma(x,y)$, the associated kernel of $\sigma(\cdot,\LL)$: for any $0<\ve <1$ there exists $C_\ve>0$ such that
\begin{align}\label{eq:czk2}
\big|\partial_2^\eta K^\sigma(x,y)-\partial_2^\eta K^\sigma(x,y')\big|\lesi \f{|y-y'|^\ve}{|x-y|^{n+M+\ve}},\qquad |\eta|=M,
\end{align}
for every $|x-y|>2|y-y'|$, and
\begin{align}\label{eq:czk1}
\big|\partial_1^\gamma K^\sigma(x,y)-\partial_1^\gamma K^\sigma(x',y)\big|\lesi \f{|x-x'|^\ve}{|x-y|^{n+M+\ve}},\qquad |\gamma|=M,
\end{align}
for every $|x-y|>2|x-x'|$.
\begin{enumerate}[\upshape(a)]
\item
If $\sigma\in\SM_{1,1}^{0,n+M+1,0}$ then \eqref{eq:czk2} holds.
\item If $\sigma\in\SM_{1,1}^{0,n+M+1,M+1}$, then both \eqref{eq:czk2} and \eqref{eq:czk1} holds.
\end{enumerate}
\end{Proposition}
\begin{proof}[Proof of Proposition \ref{prop:czk}]
Applying Lemma \ref{lem:ddKj} with $\K=n+M+1$, $m=0$, $\rr=1$ and  $\N\ge 0$ gives
\begin{align}\label{eq:czk3}
|x-y|^N \big|\partial_1^\gamma \partial_2^\eta K_j^\sigma(x,y)\big| \lesi 2^{j(n+|\gamma|+|\eta|-N)} 
\end{align}
for $\eta\in \NN_0^n$, $|\gamma|\le \N$ and $N\le n+M+1$. 

We first consider part (a). Let  $N =n+M+\ve$. 
Then taking  $|\gamma|=0$ in \eqref{eq:czk3} we have
\begin{align}\label{eq:czk4}
|x-y|^{n+M+\ve} \big|\partial_2^\eta K_j^\sigma(x,y)\big|
\lesi 2^{-j\ve}, \quad |\eta|=M,
\end{align}
and
\begin{align}
\label{eq:czk5}
|x-y|^{n+M+\ve} \big|\partial_2^\eta K_j^\sigma(x,y)\big|
\lesi 2^{j(1-\ve)}, \quad |\eta|=M+1.
\end{align}
Let $|\eta|=M$. Then we write
\begin{align*}
&|x-y|^{n+M+\ve}\big| \partial_2^\eta K^\sigma(x,y)-\partial_2^\eta K^\sigma(x,y')\big| \\
&\qquad\qquad\le \sum_{j\ge 0} |x-y|^{n+M+\ve}\big| \partial_2^\eta K_j^\sigma(x,y)-\partial_2^\eta K_j^\sigma(x,y')\big| 
=I+II
\end{align*}
where
\begin{align*}
I&:=\sum_{2^j>|y-y'|^{-1}} |x-y|^{n+M+\ve}\big| \partial_2^\eta K_j^\sigma(x,y)-\partial_2^\eta K_j^\sigma(x,y')\big|\\
II&:=\sum_{2^j>|y-y'|^{-1}} |x-y|^{n+M+\ve}\big| \partial_2^\eta K_j^\sigma(x,y)-\partial_2^\eta K_j^\sigma(x,y')\big|.
\end{align*}

By \eqref{eq:czk4} we have 
\begin{align*}
I
&=\sum_{j=\floor{-\log_2 |y-y'|}+1}^\infty |x-y|^{n+M+\ve}\big| \partial_2^\eta K_j^\sigma(x,y)-\partial_2^\eta K_j^\sigma(x,y')\big| \\
&\le \sum_{j=\floor{-\log_2 |y-y'|}+1}^\infty |x-y|^{n+M+\ve}\Big(\big| \partial_2^\eta K_j^\sigma(x,y)\big|+\big|\partial_2^\eta K_j^\sigma(x,y')\big|\Big) \\
&\le \sum_{j=\floor{-\log_2 |y-y'|}+1}^\infty 2^{-j\ve}\Big(1+\f{|x-y|^{n+M+\ve}}{|x-y'|^{n+M+\ve}} \Big).
\end{align*}
We note that since $|x-y|>2|y-y'|$ then $|x-y|$ is comparable to $|x-y'|$. Thus we can conclude our estimate of term $I$ by 
\begin{align*}
I\lesi\sum_{j=\floor{-\log_2 |y-y'|}+1}^\infty  2^{-j\ve}
\lesi 2^{-\ve\floor{-\log_2|y-y'|}} 
\sim |y-y'|^\ve.
\end{align*}
For the second term $II$ we use the mean value theorem and estimate \eqref{eq:czk5} to obtain
\begin{align*}
II
&= \sum_{j=0}^{\floor{-\log_2|y-y'|}}|x-y|^{n+M+\ve}\big| \partial_2^\eta K_j^\sigma(x,y)-\partial_2^\eta K_j^\sigma(x,y')\big|  \\
&\le \sum_{j=0}^{\floor{-\log_2|y-y'|}}|x-y|^{n+M+\ve} |y-y'|\big| \nabla_2\partial_2^\eta K_j^\sigma(x,\wt{y})\big| \\
&\lesi \sum_{j=0}^{\floor{-\log_2|y-y'|}} 2^{j(1-\ve)}|y-y'|\f{|x-y|^{n+M+\ve} }{|x-\wt{y}|^{n+M+\ve}},
\end{align*}
where $\wt{y}$ is some point on the line segment joining $y$ and $y'$. Now $|x-y|>2|y-y'|$ also implies that $|x-y|$ is comparable to $|x-\wt{y}|$. Thus we have
\begin{align*}
II
\lesi |y-y'|\sum_{j=0}^{\floor{-\log_2|y-y'|}} 2^{j(1-\ve)} 
\lesi |y-y'| 2^{(1-\ve)\floor{-\log_2|y-y'|}}
\sim |y-y'|^\ve.
\end{align*}
This yields \eqref{eq:czk2} and completes the proof of part (a).

For part (b) we observe that $\SM^{0,n+M+1,M+1}_{1,1}\subset\SM^{0,n+M+1,0}_{1,1}$ so that \eqref{eq:czk2} holds automatically. It remains to show \eqref{eq:czk1}, which can be obtained along the same lines as \eqref{eq:czk2}. However, in place of \eqref{eq:czk4}-\eqref{eq:czk5}, here one uses instead
\begin{align*}
|x-y|^{n+M+\ve}\big|\partial_1^\gamma K_j^\sigma(x,y)\big| 
\lesi 2^{-j\ve}, \quad |\gamma|=M, 
\end{align*}
and
\begin{align*}
|x-y|^{n+M+\ve}\big|\partial_1^\gamma K_j^\sigma(x,y)\big| 
\lesi 2^{j(1-\ve)}, \quad |\gamma|=M+1.
\end{align*}
which can be derived by setting $|\eta|=0$ and taking $|\gamma|$ as $M$ or $M+1$ respectively in \eqref{eq:czk3}.
We shall omit the rest of the details. 
\end{proof}
This completes the proof of Proposition \ref{prop:czk}, and hence also that of Theorem \ref{thm:pdo-hczo}.

\subsubsection{Proof of Theorem \ref{thm:pdo-hardy}}

Let $p$ and $M$ be as given in the assumptions and choose $0<\ve <1$ such that $\f{n}{n+M+\ve}<p$. Our aim is to apply Theorem \ref{thm:hardy} to $T=\sigma(\cdot,\LL)$ and $\om=M+\ve$.

Firstly, observe from Theorem \ref{thm:pdo-hczo} that under either condition we have $\sigma(\cdot,\LL)\in \HCZO_{(2)}(M,\ve)$. Our main task here is to show that the cancellation condition \eqref{eq:hardycond} holds with $\om=M+\ve$. We set 
\begin{align}\label{eq:pdo hardy1} g_{x_0,\beta}(x):=(x-x_0)^\beta \chi(x).\end{align}
By duality and orthogonality we may write
\begin{align*}
\vph_j(\sqrt{\LL})\sigma(x,\LL) ^* g_{x_0,\beta}(x) = \sum_{\xi\in\I_j} \vph_j(\xi)h_\xi(x)\int\sigma(y,\xi)h_\xi(y)g_{x_0,\beta}(y)\,dy.
\end{align*}
Then by Cauchy--Schwarz and estimate \eqref{eq:htoQQ} we have
\begin{align*}
\big|\vph_j(\sqrt{\LL})\sigma(x,\LL) ^* g_{x_0,\beta}(x)\big|
&\lesi \Big(\sum_{\xi\in\I_j}h_\xi(x)^2\Big)^{\f{1}{2}}\Big(\sum_{\xi\in\I_j}\Big|\int \sigma(y,\xi)h_\xi(y) g_{x_0,\beta}(y)\,dy\Big|^2\Big)^{\f{1}{2}}\\
&\lesi 2^{jn/2}\Big(\sum_{\xi\in\I_j}\Big|\int \sigma(y,\xi)h_\xi(y) g_{x_0,\beta}(y)\,dy\Big|^2\Big)^{\f{1}{2}}.
\end{align*}
At this point we invoke the estimates in Lemma \ref{lem:CN}. Under assumption (a) we apply estimate \eqref{eq:CN1} from Lemma \ref{lem:CN} to get
\begin{align*}
\Big|\vph_j(\sqrt{\LL}) \sigma(x,\LL)^* g_{x_0,\beta}(x)\Big|
\lesi  \cro(x_0)^{|\beta|+n-2N}2^{-j(2N-n)} \max\{1,2^j\cro(x_0)\}^{2N\dd}
\end{align*}
Thus for any $\om>0$, 
\begin{align*}
2^{j\om}\Vert \vph_j(\sqrt{\LL}) \sigma(\cdot,\LL)^* g_{x_0,\beta}\Vert_{L^\infty}
&\lesi \cro(x_0)^{|\beta|-\om} \Big(\f{1+|x_0|}{2^j}\Big)^{2N-n-\om}\max\Big\{1,\f{2^j}{1+|x_0|}\Big\}^{2N\dd}\\
&=\cro(x_0)^{|\beta|-\om} \Big(\f{1+|x_0|}{2^j}\Big)^{2N(1-\dd)-n-\om}\max\Big\{\f{1+|x_0|}{2^j},1\Big\}^{2N\dd}.
\end{align*}
If $1+|x_0|\ge 2^j$ we choose $N=0$. Otherwise we choose $N =\ceil{\f{n+\om}{2(1-\dd)}}$. In either case we obtain 
\begin{align*}
2^{js}\Vert \vph_j(\sqrt{\LL}) \sigma(\cdot,\LL)^* g_{x_0,\beta}\Vert_{L^\infty}
\lesi \cro(x_0)^{|\beta|-s}.
\end{align*}
Thus setting $\om=M+\ve$ gives us the required result since 
$$ \N=2\ceil{\tfrac{n+M+1}{2(1-\dd)}} \ge 2\ceil{\tfrac{n+M+\ve}{2(1-\dd)}}, \qquad \forall\; \ve\in [0,1].$$
For assumption (b)  we use estimate \eqref{eq:CN2} of Lemma \ref{lem:CN} to get
\begin{align*}
\Big|\vph_j(\sqrt{\LL}) \sigma(x,\LL)^* g_{x_0,\beta}(x)\Big|
\lesi  \cro(x_0)^{|\beta|+n-2N}2^{-j(2N-n)}.
\end{align*}
Then by setting $N=0$ if $1+|x_0|\le 2^j$, and $N=\floor{\f{n+M}{2}}+1$ otherwise, we have
\begin{align*}
2^{j(M+\ve)}\Vert \vph_j(\sqrt{\LL}) \sigma(\cdot,\LL)^* g_{x_0,\beta}\Vert_{L^\infty}
\lesi \cro(x_0)^{|\beta|-M-\ve}  \Big(\f{1+|x_0|}{2^j}\Big)^{2N-n-M-\ve}
\lesi \cro(x_0)^{|\beta|-M-\ve}.
\end{align*}
This completes the proof of  \eqref{eq:hardycond} under both assumptions (a) and (b). 

We may now invoke Theorem \ref{thm:hardy} to complete the proof of  Theorem \ref{thm:pdo-hardy}.

\subsubsection{Proof of Theorem \ref{thm:pdo lip}}
Under either assumption (a) or (b) we have $\sigma(\cdot,\LL)\in \HCZO(0,\ve)$ for any $\ve\in(0,1)$. For part (a) this fact follows directly from Theorem \ref{thm:pdo-hczo} (a) for the case $M=0$. For part (b), we follow the proof of Theorem \ref{thm:pdo-hczo} (b) except here we may use Proposition \ref{prop:czk} (b) instead of (a).

Next we claim that the following estimate holds. For any $0<\om<1$,
\begin{align}\label{eq:pdo lip1}
\sup_{x_0\in\RR^n} \cro(x_0)^\om \Vert \sigma(\cdot,\LL) \chi \Vert_{\Lip^\om_\LL} <\infty.
\end{align}
with $\chi\in C^\infty_0(B(x_0,2\cro(x_0))$ with $\chi=1$ on $B(x_0,\cro(x_0))$.
Assuming \eqref{eq:pdo lip1},  then given any $0<s<1$ we may then invoke Theorem \ref{thm:lip}  for some $s<\om<1$ to conclude our result.

We now devote our attention to obtaining \eqref{eq:pdo lip1}, which will complete our proof. 
Firstly, arguing as in \cite[Theorem 4.3]{LN21} we have the following estimates for each $0\le~|\gamma|\le~1$,
\begin{align}\label{eq:pdo lip2a}
\big|\partial_x^\gamma \sigma(x,\LL)\vph_k(\sqrt{\LL})(x,y)\big| \lesi \f{2^{k(n+1)}}{(1+2^k|x-y|)^{n+1}}, \qquad k\in\NN_0.
\end{align}
Using \eqref{eq:pdo lip2a} and arguing as in \cite[Lemma 3.4]{LN21} we also have 
\begin{align}\label{eq:pdo lip2b}
\big|\vph_j(\sqrt{\LL})\sigma(x,\LL)\vph_k(\sqrt{\LL})(x,y)\big|
\lesi \f{2^{(j\land k)n-[(j-k)\vee 0]}}{(1+2^{j\land k}|x-y|)^{n+1}}, \qquad j,k\in\NN_0.
\end{align}
Now choose an admissible function $\vph$ satisfying \eqref{eq:partition}.
Observe that in view of the support of $\vph$ we have 
$$ \vph_k(\sqrt{\LL})\vph_\ell(\sqrt{\LL}) =0, \qquad \text{whenever}\quad |\ell-k|\ge 3.$$
As a consequence, we have 
\begin{align*}
\vph_k(\sqrt{\LL}) = \sum_{k\ge 0}\vph_k(\sqrt{\LL})\vph_\ell(\sqrt{\LL})
=\sum_{\ell=k-2}^{k+2} \vph_k(\sqrt{\LL})\vph_\ell(\sqrt{\LL}),
\end{align*}
and therefore, for each $j\in\NN_0$,
\begin{align*}
\vph_j(\sqrt{\LL})\sigma(\cdot,\LL)\chi
= \sum_{k\ge 0} \sum_{\ell=k-2}^{k+2} \vph_j(\sqrt{\LL})\sigma(\cdot,\LL)\vph_k(\sqrt{\LL})\vph_\ell(\sqrt{\LL})\chi.
\end{align*}
Then for each $x\in\RR^n$,
\begin{align*}
\big|\vph_j(\sqrt{\LL})\sigma(x,\LL)\chi(x)\big|
&\le \sum_{k\ge 0} \sum_{\ell=k-2}^{k+2} \int \big| \vph_j(\sqrt{\LL})\sigma(x,\LL)\vph_k(\sqrt{\LL})(x,y)\big| |\vph_\ell(\sqrt{\LL})\chi(y)|\,dy\\
&\le  \sum_{k\ge 0} \sum_{\ell=k-2}^{k+2} \Vert\vph_\ell(\sqrt{\LL})\chi\Vert_\infty \int \big| \vph_j(\sqrt{\LL})\sigma(x,\LL)\vph_k(\sqrt{\LL})(x,y)\big|\,dy.
\end{align*}
Now observe  that by Lemma \ref{lem:lip eg} we have  $\chi\in \Lip_\LL^\om(\RR^n)$ so that $\Vert \chi \Vert_{\Lip_\LL^\om}\lesi \cro(x_0)^{-\om}$. In particular it follows that
\begin{align*}
\Vert \vph_\ell(\sqrt{\LL}) \chi\Vert_{\infty} 
= 2^{-\ell\om} \big(2^{\ell\om}\Vert \vph_\ell(\sqrt{\LL})\chi\Vert_\infty\big)
\le 2^{-\ell\om} \Vert \chi \Vert_{\Lip_\LL^\om}
\lesi 2^{-\ell\om}\cro(x_0)^{-\om}
\sim 2^{-k\om}\cro(x_0)^{-\om}
\end{align*}
since $|\ell-k|\le 2$. Hence we arrive at
\begin{align*}
\big|\vph_j(\sqrt{\LL})\sigma(x,\LL)\chi(x)\big|
\lesi \cro(x_0)^{-\om}\sum_{k\ge 0} 2^{-k\om} \int \big| \vph_j(\sqrt{\LL})\sigma(x,\LL)\vph_k(\sqrt{\LL})(x,y)\big|\,dy.
\end{align*}

Now invoking \eqref{eq:pdo lip2b} we obtain
\begin{align*}
\big|\vph_j(\sqrt{\LL})\sigma(x,\LL)\chi(x)\big|
&\lesi \cro(x_0)^{-\om}\sum_{k\ge 0} 2^{-k\om} \int \big| \vph_j(\sqrt{\LL})\sigma(x,\LL)\vph_k(\sqrt{\LL})(x,y)\big|\,dy\\
&\lesi 2^{-j\om}\cro(x_0)^{-\om}\sum_{k\ge 0} 2^{(j-k)\om-[(j-k)\vee 0]}\int \f{2^{(j\land k)n}}{(1+2^{(j\land k)}|x-y|)^{n+1}}dy\\
&\lesi 2^{-j\om}\cro(x_0)^{-\om}\Big(\sum_{k<j} 2^{-(j-k)(1-\om)} +\sum_{k\ge j}2^{-(k-j)\om}\Big).
\end{align*}
Since $0<\om<1$ then we conclude
\begin{align*}
\big|\vph_j(\sqrt{\LL})\sigma(x,\LL)\chi(x)\big|
\lesi \cro(x_0)^{-\om}2^{-j\om}, \qquad x\in\RR^n.
\end{align*}
The preceding estimate implies 
$$ \cro(x_0)^{\om}2^{j\om}\big\Vert\vph_j(\sqrt{\LL})\sigma(\cdot,\LL)\chi\big\Vert_\infty <\infty, $$
and taking supremum over  $j\in\NN_0$we arrive at \eqref{eq:pdo lip1}.
This completes the proof of Theorem \ref{thm:pdo lip}.

\subsection{ Riesz transforms}\label{sec:RT}

For each $i\in\{1,\dots,n\}$ we define hermite derivatives by
\begin{align}\label{eq:hermitederiv}
A_i^* = \f{\partial}{\partial x_i}+x_i \qquad\text{and}\qquad A_i=-\f{\partial}{\partial x_i}+x_i.
\end{align}
These are sometimes called the \emph{annihilation} and \emph{creation} operators respectively. Moreover the operator $\LL$ can be factored as $\LL= \f{1}{2} \sum_i^n(A_iA_i^* + A_i^*A_i)$. Additional properties relevant to this paper is described in Appendix \ref{app:hermite}.

Let $\alpha\in\NN_0^n$ and set $\AH^\alpha=\AH_1^{\alpha_1}\dots\AH_n^{\alpha_n}$, with $\AH_k\in \{A_k,A_k^*\}$. Then we define the Riesz transform of order $\alpha$ by
\begin{align*}
\R^\alpha f(x)=\AH^\alpha \LL^{-|\alpha|/2} f(x)=\sum_{\xi\in\NN_0^n}(2|\xi|+n)^{-|\alpha|/2} \ip{f,h_\xi}\AH^\alpha h_{\xi}(x)
\end{align*}
with kernel
$$ \R^\alpha (x,y)=\sum_{\xi\in\NN_0^n}(2|\xi|+n)^{-|\alpha|/2} \AH^\alpha h_{\xi}(x) h_\xi(y).$$
Note that when $|\alpha|=1$, one has the first-order Riesz transforms 
\begin{align*}
 \R_k^- f(x):=A_k\LL^{-1/2}f(x)=\sum_{\xi\in\NN_0^n}\Big(\f{2\xi_k}{2|\xi|+n}\Big)^\half \ip{f,h_\xi}h_{\xi-e_k}(x), \\
  \R_k^+ f(x):=A_k^*\LL^{-1/2}f(x)=\sum_{\xi\in\NN_0^n}\Big(\f{2\xi_k+2}{2|\xi|+n}\Big)^\half \ip{f,h_\xi}h_{\xi-e_k}(x),
 \end{align*}
 for $1\le k\le n$. 

Both $\R_k^-$ and $\R_k^+$ were studied in \cite{ST03} and shown to be Calderon--Zygmund operators. The higher order Riesz transforms were studied in \cite{BD}; for Riesz Transforms on Lipschitz spaces see  \cite{DT20}.

\begin{Theorem}[Riesz transforms are Hermite--Calder\'on--Zygmund operators]\label{thm:RT-hczo}
Let $\alpha\in\NN_0^n$. Then for any $M\in\NN_0$ and $0<\ve\le 1$, we have $\R^\alpha\in\HCZO(M,\ve)$.
\end{Theorem}

\begin{Theorem}[Riesz transforms on endpoint spaces]\label{thm:RT}
Let $\alpha\in\NN_0^n$. Then the Riesz transforms $\R^\alpha$ are bounded on $h^p_\LL(\RR^n)$ for all $0<p\le 1$ and on $\Lip^s_\LL(\RR^n)$ for all $s\ge 0$.
\end{Theorem}

Before giving the proof of Theorem \ref{thm:RT} we derive some needed kernel estimates. 
For any admissible system $\{\vph_j\}_{j\ge 0}$ we set
\begin{align*}
\R^\alpha_j f(x) :=\R^\alpha \vph_j(\sqrt{\LL})f(x)=\sum_{\xi\in\I_j} (2|\xi|+n)^{-|\alpha|/2} \vph_j(\xi) \AH^\alpha h_\xi(x) \ip{f,h_\xi}
\end{align*}
with kernel
\begin{align*}
\R_j^\alpha(x,y)=\sum_{\xi\in\I_j} (2|\xi|+n)^{-|\alpha|/2} \vph_j(\xi)\AH^\alpha h_\xi(x)h_\xi(y).
\end{align*}
Then we have the following whose proof is found in Appendix \ref{app:technical}.
\begin{Lemma}\label{lem:RTddKj}
For any $\mu\ge0$, $\gamma,\eta\in\NN_0^n$ and $N\ge 0$ we  have
\begin{align}\label{eq:RTKjab}
|x-y|^{N}\big|\partial_x^\gamma\partial_y^\eta \R_j^\alpha(x,y)\big|
\lesssim 2^{j(n+|\gamma|+|\eta|-N)}\big(1+\tfrac{2^{-j}}{\cro(x)}+\tfrac{2^{-j}}{\cro(y)}\big)^{-\mu}
\end{align}
for all $j\in \NN_0$.
\end{Lemma}

\subsubsection{Proof of Theorem \ref{thm:RT-hczo}}
Firstly note that from \eqref{eq:hermitederiv} (see also \eqref{eq:d1a}--\eqref{eq:d1b}), we have
\begin{align}\label{eq:RTddKj0}
\AH^\alpha h_\xi(x) = a_\alpha(\xi)\, h_{\xi+\wt{\alpha}}(x)
\end{align}
where 
\begin{align*}
	\wt{\alpha}_i = \left
\lbrace 
	\begin{array}{ll}
			\alpha_i  \qquad &\text{if}\quad \AH_i=A_i,\\
			-\alpha_i \qquad &\text{if}\quad \AH_i=A_i^*,
	\end{array}
\right.
\end{align*}
and $a_\alpha(\xi)$ are numbers satisfying $|a_\alpha(\xi)|\lesi \ip{|\xi|}^{|\alpha|/2}$. Note that $a_\alpha(\xi)=0$ if for some $i$, $\xi_i+\wt{\alpha}_i\le 0$.

From \eqref{eq:RTddKj0}, and then setting $\sigma_\alpha(\xi):=(2|\xi|+n)^{-|\alpha|/2}$, we have
 \begin{align*}
\R^\alpha f(x)
=\sum_{\xi\in\NN_0^n} \sigma_\alpha(\xi) a_\alpha(\xi) \ip{f,h_\xi} h_{\xi+\wt{\alpha}}(x).
\end{align*}
From this  expression we see that $\R^\alpha$ is a special kind of multiplier with symbol $m(\xi)=\sigma_\alpha(\xi) a_\alpha(\xi)\in L^\infty$. Thus by Parseval's identity $\R^\alpha$ is bounded on $L^2(\RR^n)$, 

Next, it can be readily seen that the size and smoothness estimates (ii),(iii) and (v) of Definition \ref{def:HCZO} follow easily from the following result.
\begin{Proposition}\label{prop:kernelRT}
Suppose that  $\gamma,\eta\in\NN_0^n$ and $\mu\ge 0$. Then 
\begin{align}\label{eq:kernelRT}
\big|\partial_x^\gamma \partial_y^\eta \R^\alpha(x,y)\big| \lesssim  |x-y|^{-n-|\gamma|-|\eta|}\Big(1+\tfrac{|x-y|}{\cro(x)}+\tfrac{|x-y|}{\cro(y)}\Big)^{-\mu}.
\end{align}
\end{Proposition}

\begin{proof}[Proof of Proposition \ref{prop:kernelRT}]
The proof is similar to the proof of Proposition \ref{prop:kernelV1}, and we only highlight the main differences. 

Set $n_{\gamma,\eta}:=n+|\gamma|+|\eta|>0$. 
We write
\begin{align*}
\big|\partial_y^\gamma \partial_x^\eta \R^\alpha(x,y)\big| 
\le\sum_{2^j \le |x-y|^{-1}} \big|\partial_y^\gamma \partial_x^\eta \R_j^\alpha(x,y)\big|+\sum_{2^j>|x-y|^{-1}} \big|\partial_y^\gamma \partial_x^\eta \R_j^\alpha(x,y)\big|
=:I+II.
\end{align*}
Then for term $I$ we can argue exactly as in Proposition \ref{prop:kernelV1}, but invoking Lemma \ref{lem:RTddKj} (again with $N=0$) in place of Lemma \ref{lem:ddKj}.

Term $II$ is also handled similarly to Proposition \ref{prop:kernelV1} but in place of Lemma \ref{lem:ddKj} we use again Lemma \ref{lem:RTddKj}, this time with $N=n_{\gamma,\eta}+\mu+1$. This gives
\begin{align*}
II 
&\lesi |x-y|^{-(n_{\gamma,\eta}+\mu+1)}\sum_{j=\floor{-\log_2(|x-y|)}+1}^\infty 2^{-j(\mu+1)}\big(1+\tfrac{2^{-j}}{\cro(x)}+\tfrac{2^{-j}}{\cro(y)}\big)^{-\mu}\\
&\lesi |x-y|^{-(n_{\gamma,\eta}+1)}\big(1+\tfrac{|x-y|}{\cro(x)}+\tfrac{|x-y|}{\cro(y)}\big)^{-\mu}\sum_{j=\floor{-\log_2(|x-y|)}+1}^\infty 2^{-j}.
\end{align*}
Now noting that 
\begin{align*}
\sum_{j=\floor{-\log_2(|x-y|)}+1}^\infty 2^{-j}
\lesi \big(\tfrac{1}{2}\big)^{\floor{-\log_2(|x-y|)}}
\sim |x-y|,
\end{align*}
 we arrive at the desired estimate for term $II$, and hence obtain \eqref{eq:kernelRT} as required.
\end{proof}
This completes the proof of Proposition \ref{prop:kernelRT} and hence of Theorem \ref{thm:RT-hczo}.

\subsubsection{Proof of Theorem \ref{thm:RT}}
We shall apply Theorems \ref{thm:hardy} and \ref{thm:lip} to $\R^\alpha$. 
The fact that $\R^\alpha$ is a Hermite--Calder\'on--Zygmund operator for any $M\ge 0$ and $\ve \in [0,1]$ was already established in Theorem \ref{thm:RT-hczo}. 

It remains to show the conditions \eqref{eq:hardycond} and \eqref{eq:lipcond}. Recall the definition of $g_{x_0,\beta}$ from \eqref{eq:pdo hardy1}. We first show \eqref{eq:hardycond}. 
From duality, the definition of $\R^\alpha$, and Cauchy--Schwarz, we have
 \begin{align*}
 \big|\vph_j(\sqrt{\LL})(\R^\alpha)^* g_{x_0,\beta}(x)\big| 
 &=\Big|\sum_{\xi\in\I_j}\sigma_\alpha(\xi) \vph_j(\xi) h_\xi(x)\int \AH^\alpha h_\xi(y) g_{x_0,\beta}(y)\,dy\Big|\\
 &\lesi \Big(\sum_{\xi\in\I_j} \sigma_\alpha(\xi) h_\xi(x)^2\Big)^\half \Big(\sum_{\xi\in\I_j} \Big|\int\AH^\alpha h_\xi(y) g_{x_0,\beta}(y)\,dy\Big|^2\Big)^\half.
 \end{align*}
Since $(2|\xi|+n) \sim 2^j$ then by applying \eqref{eq:QQ est} and estimate \eqref{eq:CN3} from Lemma \ref{lem:CN} we have
  \begin{align*}
 \big|\vph_j(\sqrt{\LL})(\R^\alpha)^* g_{x_0,\beta}(x)\big|
 &\lesi 2^{j(n/2-|\alpha|)} \Big(\sum_{\xi\in\I_j} \Big|\int\AH^\alpha h_\xi(y) g_{x_0,\beta}(y)\,dy\Big|^2\Big)^\half \\
 &\lesi 2^{-j(2N-n)} \cro(x_0)^{|\beta|+n-2N}.
 \end{align*}
 Then for any $\om>0$ we have, for any $\beta\in\NN_0^n$, 
 \begin{align}\label{eq:RT1}
2^{j\om}\Vert \vph_j(\sqrt{\LL}) (\R^\alpha)^* g_{x_0,\beta}\Vert_{L^\infty}
\lesi \cro(x_0)^{|\beta|-\om} \Big(\f{1+|x_B|}{2^j}\Big)^{2N-n-\om}
\lesi \cro(x_0)^{|\beta|-\om}.
\end{align}
Now if $1+|x_B|\ge 2^j$ then we choose $N=0$, and $N=\ceil{\f{n+\om}{2}}$ otherwise, yielding \eqref{eq:hardycond}.

Next we turn to \eqref{eq:lipcond}. By orthogonality of the hermite functions,
\begin{align*}
\vph_j(\sqrt{\LL})\R^\alpha g_{x_0,\beta}
&= \sum_{\xi\in\I_j}\sum_{\eta\in\NN_0^n} \vph_j(\xi)\sigma_\alpha(\eta) \ip{g_{x_0,\beta},h_\eta}\ip{\AH^\alpha h_\eta,h_\xi} h_\xi\\
&=\sum_{\xi\in\I_j} \vph_j(\xi) \sigma_\alpha(\xi-\wt{\alpha}) a_\alpha(\xi-\wt{\alpha}) \ip{g_{x_0,\beta},h_{\xi-\wt{\alpha}}}h_\xi.
\end{align*}
From the following estimate
$$ |\sigma_\alpha(\xi-\wt{\alpha}) a_\alpha(\xi-\wt{\alpha}) |\lesi \ip{\xi-\wt{\alpha}}^{-|\alpha|/2} \ip{\xi}^{|\alpha|/2} \lesi 1,$$
and Cauchy--Schwarz, we have
\begin{align*}
\big|\vph_j(\sqrt{\LL})\R^\alpha g_{x_0,\beta}(x)\big|
\lesi \sum_{\xi\in\I_j} \big|\ip{g_{x_0,\beta},h_{\xi-\wt{\alpha}}}\big| |h_\xi(x)|
\le \Big(\sum_{\xi\in\I_j} \big|\ip{g_{x_0,\beta},h_{\xi-\wt{\alpha}}}\big|^2\Big)^\half  \Big(\sum_{\xi\in\I_j} h_\xi(x)^2\Big)^\half.
\end{align*}
We next observe that the following estimate can be obtained much the same way as in \eqref{eq:CN3}: for every $\beta\in\NN_0^n$ and $N\ge 0$ we have 
\begin{align*}
 \Big(\sum_{\xi\in\I_j} \big|\ip{g_{x_0,\beta},h_{\xi-\wt{\alpha}}}\big|^2\Big)^\half
 \lesi 2^{-j(2N-n/2)}\cro(x_0)^{|\beta|+n-2N}.
\end{align*}
Using this estimate, along with \eqref{eq:htoQQ}, we have 
\begin{align*}
\big|\vph_j(\sqrt{\LL})\R^\alpha g_{x_0,\beta}(x)\big|
\lesi 2^{-j(2N-n)}\cro(x_0)^{|\beta|+n-2N}.
\end{align*}
Now we may proceed as in \eqref{eq:RT1} to conclude our proof of Theorem \ref{thm:RT}.

\appendix
\section{Useful Hermite identities and estimates}\label{app:hermite}
In this section we give some key identities and estimates that are important for our proofs. The main results here are Lemma \ref{lem: identities} and Lemma \ref{lem:CN}. 
For a definition of the Hermite derivatives $A_i$, $A^*_i$ and $\AH_i$ see Section \ref{sec:RT}. 

\begin{Lemma}\label{lem: identities}
\begin{enumerate}[\upshape(a)]
\item Suppose that
$$ K(x,y) = \sum_{\xi\in\NN_0^n} k(x,y,\xi) \,h_\xi(x)h_\xi(y), \qquad x,y\in\RR^n.$$
If $N\in\ZZ^+,$ it holds that
\begin{align}\label{eq:identity A}
	2^N(x_i- y_i)^N K(x,y) = \sum_{\f{N}{2}\le \ell\le N}c_{\ell, N}\sum_{\xi\in \NN_0} \diff_i^\ell k(x,y,\xi)\big(A^{(y)}_i-A^{(x)}_i\big)^{2\ell-N} h_\xi(x) h_\xi(y),
\end{align}
where $c_{\ell,N} = (-4)^{N-\ell}(2N-2\ell-1)!!\binom{N}{2\ell-N}$.

\item If $\mu,\alpha\in\NN_0^n,$  $m\in \NN_0$ and $i\in\{1,\dots, n\},$ it holds that for all $x\in\RR^n$,
\begin{align} \label{eq:identity d1}
\big|\big(A^{(x)}_i\big)^m h_\mu(x)\big|& \le \big[2(\mu_i+m)+2\big]^{\f{m}{2}} | h_{\mu+m e_i}(x)| 
\end{align}
and 
\begin{align} \label{eq:identity d2}
\big|\big(A^{(x)}\big)^\alpha h_\mu(x)\big| &\le \big[2(|\mu|+|\alpha|)+2\big]^{\f{|\alpha|}{2}} | h_{\mu+\alpha}(x)|.
\end{align}

\item For each $\beta,\xi\in\NN_0^n$, it holds that for all $x\in\RR^n$,
\begin{align}\label{eq:identity C}
x^\beta \,h_\xi(x)= \sum_{\om\le \beta} b_{\omega,\beta}(\xi)h_{\xi+\beta-2\omega}(x)
\end{align}
where $b_{\om,\beta}(\xi)=\prod_{i=1}^n b_{\om_i,\beta_i}(\xi_i)$ with $ b_{\omega_i,\beta_i}(\mu_i)=0$ if $\mu_i+\beta_i-2\omega_i<0$ and  $ b_{\omega_i,\beta_i}(\mu_i)\sim \mu_i^{\beta_i/2}$ otherwise. 
\item For $m\ge 1$ we have for each $x\in\RR^n$,
\begin{align}\label{eq:d1a}
(A_i^{(x)})^m h_\mu(x)=\prod_{r=0}^{m-1}\sqrt{2(\mu_i+r)+2} \,h_{\mu+me_i}(x),
\end{align}
and if $\mu_i\ge m$,
\begin{align}\label{eq:d1b}
(A_i^{* (x)})^m h_\mu(x)=\prod_{r=0}^{m-1}\sqrt{2(\mu_i-r)} \,h_{\mu-me_i}(x).
\end{align}
\item Let $\AH_i\in \{A_i,A_i^*\}$. If $N, M\in\NN,$ then it holds that for each $x\in\RR^n$,
\begin{align}
	x_i^M \big(\AH^{(x)}_i-\AH^{(y)}_i\big)^N &= \sum_{k=0}^M c_{\AH_i,k}\tbinom{M}{k} \tfrac{N!}{(N-k)!}\big(\AH^{(x)}_i-\AH^{(y)}_i\big)^{N-k}x_i^{M-k} \label{eq:identity b1}
	\end{align}
	and
	\begin{align}
	(x_i-y_i)^N\big(\AH^{(x)}_i\big)^M &=\sum_{k=0}^M c_{\AH_i,k} \tbinom{M}{k}  \tfrac{N!}{(N-k)!}\big(\AH^{(x)}_i\big)^{M-k}(x_i-y_i)^{N-k},\label{eq:identity b2}
\end{align}
where $\f{N!}{(N-k)!}$ is defined to be 0 whenever $N<k$, and $ c_{\AH_i,k} = (-1)^k$ if $\AH_i=A_i^*$, and equal to 1 otherwise. 
\item If $\ell \in\NN_0,$ it holds that
\begin{align}\label{eq:leibniz1}
 \diff_\xi^\ell \big(f(\xi)\,g(\xi)\big) = \sum_{r=0}^\ell \binom{\ell}{r}\diff_\xi^r f(\xi)\,\diff_\xi^{\ell-r}g(\xi+r).
\end{align}
\item If $\alpha\in\NN_0^n$, it holds that
\begin{align}\label{eq:leibniz2}
A^\alpha (fg) =\sum_{\nu\le\alpha} \binom{\alpha}{\nu} (-1)^\nu \partial^\nu f \, A^{\alpha-\nu} g.
\end{align}
\end{enumerate}
\end{Lemma}

\begin{proof}[Proof of Lemma \ref{lem: identities}]
For (a)-(c) and (f)-(g) the reader is directed to \cite{LN21, Ly22}. Part (d) follows readily from the repeated application of \eqref{eq:hermitederiv}. For the commuting identities in part (e), the case $\AH_i=A_i$ is given in \cite[Lemma 9]{PX}. We shall prove (e) in general for $\AH_i \in \{A_i, A_i^*\}$, recovering \cite[Lemma 9]{PX} as a special case. Then we have 
Suppose $i\in\{1,\dots,n\}$ has been fixed. We first obtain some preliminary identities. Consider the commutator of $x$ and $T$ given by  $[x,T](f):=x(Tf)- T(xf)$. 
\begin{align}\label{eq:IDproof1}
\big[x_i,\AH_i^{(y)}\big](1)=\big[y_i,\AH_i^{(x)}\big](1) = 0  
\qquad\text{and}\qquad
\big[x_i,\AH_i^{(x)}\big](1) = 
\left\lbrace \begin{array}{cl}
1 &\quad \text{if}\; \AH_i=A_i\\
-1 &\quad \text{if}\; \AH_i=A_i^*.
 \end{array}\right. 
\end{align}
The first two equalities are easy and follow from the commutativity of the elements. For the third expression in \eqref{eq:IDproof1}, if $\AH_i=A_i$ then by the Leibniz rule $A_i(fg)=gA_if-f\partial_ig$ we have
$$ \big[x_i,\AH_i^{(x)}\big](1) = x_iA_i^{(x)}(1)-A_i^{(x)}(x_i) =  x_iA_i^{(x)}(1) - \big( x_iA_i^{(x)}(1) -\partial_ix_i\big) =1.$$
On the other hand if $\AH_i=A_i^*$ then by the Leibniz rule $A^*_i(fg)=gA^*_if+f\partial_ig$ 
we have
$$ \big[x_i,\AH_i^{(x)}\big](1) = x_iA_i^{* (x)}(1)-A_i^{* (x)}(x_i) =  x_iA_i^{* (x)}(1) - \big( xA_i^{* (x)}(1) +\partial_ix_i\big) =-1.$$
To continue we set $c_{\AH_i,k}=\big(\big[x_i,\AH_i^{(x)}\big](1)\big)^k$ for $k\in\NN_0$. Clearly we have $c_{\AH_i,0}=1$ and $c_{\AH_i,k}c_{\AH_i,\ell}=c_{\AH_i,k+\ell}$. We next wish to show that \eqref{eq:IDproof1} leads to the following identities: for $N\ge 1$ we have
\begin{align}
x_i\Axy^N&=c_{\AH_i,1} N \Axy^{N-1} +\Axy^N x_i  \label{eq:IDproof2.1}\\
(x_i-y_i)^N \AH_i^{(x)} &= c_{\AH_i,1} N (x_i-y_i)^{N-1} +\AH_i^{(x)}(x_i-y_i)^N.  \label{eq:IDproof2.2}
\end{align}
We shall prove these by induction on $N$. Firstly we consider \eqref{eq:IDproof2.1}. For $N=1$ we have by \eqref{eq:IDproof1},
\begin{align*}
x_i \Axy = \big[x_i,\AH_i^{(x)}\big](1)+\AH_i^{(x)}x_i-\AH_i^{(y)}x_i = c_{\AH_i,1} + \Axy x_i.
\end{align*}
Assume now \eqref{eq:IDproof2.1} holds for some $N\ge1$. Then from our induction hypothesis for both $N$ and $N=1$,
\begin{align*}
x_i\Axy^{N+1} 
&=\Big\{c_{\AH_i,1} N \Axy^{N-1} +\Axy^N x_i\Big\} \Axy\\
&=c_{\AH_i,1} N \Axy^{N} +\Axy^N x_i \Axy\\
&=c_{\AH_i,1} N \Axy^{N} +\Axy^N\Big\{ c_{\AH_i,1}  +\Axy x_i\Big\}\\
&=c_{\AH_i,1} (N+1) \Axy^{N} +\Axy^{N+1} x_i.
\end{align*}
For $N=1$ of \eqref{eq:IDproof2.2} we employ \eqref{eq:IDproof1} again to see that
\begin{align*}
(x_i-y_i)\AH_i^{(x)}=\big[x_i,\AH_i^{(x)}\big](1)+\AH_i^{(x)}x_i-\AH_i^{(x)}y_i 
=c_{\AH_i,1}+\AH_i^{(x)}(x_i-y_i).
\end{align*}
Assume now \eqref{eq:IDproof2.2} holds for some $N\ge1$. Then from our induction hypothesis for both $N$ and $N=1$,
\begin{align*}
(x_i-y_i)^{N+1}\AH_i^{(x)}
&=(x_i-y_i)\Big\{c_{\AH_i,1} N (x_i-y_i)^{N-1} +\AH_i^{(x)}(x_i-y_i)^N\Big\}\\
&=c_{\AH_i,1} N (x_i-y_i)^{N} +(x_i-y_i)\AH_i^{(x)}(x_i-y_i)^N\\
&=c_{\AH_i,1} N (x_i-y_i)^{N} + \Big\{c_{\AH_i,1}  +\AH_i^{(x)}(x_i-y_i)\Big\}(x_i-y_i)^N\\
&=c_{\AH_i,1} (N+1) (x_i-y_i)^{N} +\AH_i^{(x)}(x_i-y_i)^{N+1}.
\end{align*}
We are now ready to prove \eqref{eq:identity b1} and \eqref{eq:identity b2}. We shall obtain these by induction on $M$. Firstly we note that the case $M=1$ of \eqref{eq:identity b1} is given by \eqref{eq:IDproof2.1}. Now suppose that \eqref{eq:identity b1} is true for some $M\ge 1$ and all $N\ge 1$. Then from applying our inductive hypothesis followed by \eqref{eq:IDproof2.1}, we obtain
\begin{align*}
&x_i^{M+1}\Axy^N \\
&\qquad=x_i \sum_{k=0}^M\tbinom{M}{k} \f{N!}{(N-k)!} c_{\AH_i,k} \Axy^{N-k} x_i^{M-k}\\
&\qquad=\sum_{k=0}^M\tbinom{M}{k} \f{N!}{(N-k)!} c_{\AH_i,k} \Big\{c_{\AH_i,1} (N-k) \Axy^{N-k-1} +\Axy^{N-k} x_i\Big\}x_i^{M-k}\\
&\qquad=\sum_{k=0}^M\tbinom{M}{k} \f{N!}{(N-k-1)!} c_{\AH_i,k+1} \Axy^{N-k-1} x_i^{M-k} \\
&\qquad\qquad\qquad+\sum_{k=0}^M\tbinom{M}{k} \f{N!}{(N-k)!} c_{\AH_i,k} \Axy^{N-k} x_i^{M+1-k}\\
&\qquad=\sum_{k=0}^{M+1}\Big[\tbinom{M}{k-1}+\tbinom{M}{k}\Big] \f{N!}{(N-k)!} c_{\AH_i,k} \Axy^{N-k} x_i^{M+1-k}
\end{align*}
where  in the final equality it is understood that $\tbinom{M}{j}:=0$ whenever $j>M$ or $j<0$. Now invoking the binomial identity $\tbinom{M}{k-1}+\tbinom{M}{k}=\tbinom{M+1}{k}$ completes the proof of \eqref{eq:identity b1}.

The proof of \eqref{eq:identity b2} is similar. Firstly we note that the case $M=1$ of \eqref{eq:identity b2} is given by \eqref{eq:IDproof2.2}. Now suppose that \eqref{eq:identity b2} is true for some $M\ge 1$ and all $N\ge 1$. Then from applying our inductive hypothesis followed by \eqref{eq:IDproof2.2}, we obtain
\begin{align*}
&(x_i-y_i)^N\big(\AH^{(x)}_i\big)^{M+1}\\
&\qquad=\sum_{k=0}^M c_{\AH_i,k} \tbinom{M}{k}  \tfrac{N!}{(N-k)!}\big(\AH^{(x)}_i\big)^{M-k}(x_i-y_i)^{N-k}\AH^{(x)}_i \\
&\qquad= \sum_{k=0}^M c_{\AH_i,k} \tbinom{M}{k}  \tfrac{N!}{(N-k)!}\big(\AH^{(x)}_i\big)^{M-k} \Big\{c_{\AH_i,1} (N-k) (x_i-y_i)^{N-k-1} +\AH_i^{(x)}(x_i-y_i)^{N-k}\Big\}\\
&\qquad=\sum_{k=0}^M c_{\AH_i,k+1} \tbinom{M}{k}  \tfrac{N!}{(N-k-1)!}\big(\AH^{(x)}_i\big)^{M-k}(x_i-y_i)^{N-k-1} \\
&\qquad\qquad\qquad+\sum_{k=0}^M c_{\AH_i,k} \tbinom{M}{k}  \tfrac{N!}{(N-k)!}\big(\AH^{(x)}_i\big)^{M+1-k}(x_i-y_i)^{N-k}\\
&\qquad=\sum_{k=0}^{M+1} c_{\AH_i,k} \Big[\tbinom{M}{k-1}+\tbinom{M}{k}\Big]  \tfrac{N!}{(N-k)!}\big(\AH^{(x)}_i\big)^{M+1-k}(x_i-y_i)^{N-k}.
\end{align*}
Again invoking the binomial identity $\tbinom{M}{k-1}+\tbinom{M}{k}=\tbinom{M+1}{k}$ we complete the proof of \eqref{eq:identity b2}.
\end{proof}

The next lemma plays an important role  in Section \ref{sec:applications}.
\begin{Lemma}\label{lem:CN}
Let  $N\in\NN_0$. Suppose $x_0\in\RR^n$ and $\chi$ is a smooth cutoff for $B:=B(x_0,\cro(x_0))$.  (That is, $\chi\in C_0^\infty(\RR^n)$, $\chi\equiv 1$ on $B$ and $\Vert \chi^{(\gamma)}\Vert_\infty \lesi \cro(x_0)^{-|\gamma|}$.)
Then we have the following for $\beta\in\NN_0^n$. 
\begin{align}\label{eq:CN0}
\Big(\sum_{\xi\in \I_j} \Big|\int \chi(y)(y-x_0)^\beta h_\xi(y)\,dy\Big|^2\Big)^\half 
\lesi \cro(x_0)^{|\beta|+n-2N}2^{-j(2N-n/2)}.
\end{align}
If $\sigma\in\SM^{0,0,\N}_{1,\delta}$  for some $0\le \delta <1$  and $\N\ge 2N$, then
\begin{align}\label{eq:CN1}
\Big(\sum_{\xi\in \I_j}\Big|\int \chi(y)(y-x_0)^\beta h_\xi(y)\sigma(y,\xi)\,dy\Big|^2\Big)^\half
\lesi \cro(x_0)^{|\beta|+n-2N}2^{-j(2N-n/2)} \max\big\{1,2^j\cro(x_0)\big\}^{2N\delta}.
\end{align}
If $\sigma\in\SM^{0,0,0}_{1,\delta}\cap\CN^{0,M}$  for some $0\le \delta \le1$ and $M\ge 0$ then for $\wt{N}=\floor{\f{n+M}{2}}+1$ we have
\begin{align}\label{eq:CN2}
\Big(\sum_{\xi\in \I_j}\Big|\int \chi(y)(y-x_0)^\beta h_\xi(y)\sigma(y,\xi)\,dy\Big|^2\Big)^\half
\lesi \cro(x_0)^{|\beta|+n-2\wt{N}}2^{-j(2\wt{N}-n/2)}.
\end{align}
For each $\alpha\in\NN_0^n$,
\begin{align}\label{eq:CN3}
\Big(\sum_{\xi\in \I_j} \Big|\int \chi(y)(y-x_0)^\beta \AH^\alpha h_\xi(y)\,dy\Big|^2\Big)^\half 
\lesi \cro(x_0)^{|\beta|+n-2N}2^{-j(2N-n/2-|\alpha|)}.
\end{align}
In addition, all the constants above are independent of $x_0$. 
\end{Lemma}
\noindent It is worth pointing out here that  \eqref{eq:CN0} is needed in Lemma \ref{lem:lip eg}, estimates \eqref{eq:CN1}-\eqref{eq:CN2} are needed for pseudo-mulitpliers in Section \ref{sec:pdo}, and \eqref{eq:CN3} needed for Riesz transforms in Section \ref{sec:RT}.

\begin{proof}[Proof of Lemma \ref{lem:CN}]
We first prove \eqref{eq:CN0}. Here we  follow the method in \cite[Lemma A.1]{LN21}.
For any $N \in \NN_0,$ it holds that
\begin{align*}
\Big| \int_{\RR^n} (y-x_0)^\beta \chi(y) h_\xi(y)\,dy\Big|
&=(2|\xi|+n)^{-N} \Big|\int_{\RR^n} \LL_y^N\big [ (y-x_0)^\beta \chi(y)\big] h_\xi(y)\,dy\Big| \\
&\lesi \ip{\xi}^{-N} \Big\Vert \LL^N \big[(\cdot-x_0)^\beta \chi(\cdot)\big] \Big\Vert_{L^2(2B)} \Vert h_\xi\Vert_{L^2(2B)}.
\end{align*}
Repeated application of the Leibniz' rule gives
\begin{align*}
\LL^N \big[(\cdot-x_0)^\beta \chi(\cdot)\big](y)
= \sum_{a,b,\gamma,\nu,\beta} C_{a,b,\gamma,\nu} \,y^a (y-x_0)^{\beta-\gamma}\chi^{(\nu)}(y),
\end{align*}
with the summation running over $|a|+|b|\le 2N$, $\gamma+\nu=b$ and $|\gamma| \le |\beta|$.
Consequently,
\begin{align*}
 \Big\Vert \LL^N \big[(\cdot-x_0)^\beta \chi(\cdot)\big] \Big\Vert_{L^2(2B)}
 &\sim \sum_{a,b,\gamma,\nu,\beta} \Big(\int_{2B} \big| |y|^{|a|} |y-x_0|^{|\beta|-|\gamma|}|\chi^{(\nu)}(y)| \big|^2\,dy\Big)^{\half} \\
 &\lesssim  \sum_{a,b,\gamma,\nu,\beta}  \cro(x_0)^{|\beta|-|\gamma|-|\nu| +n/2} \sup_{y\in 2B}|y|^{|a|},\\
 &\lesssim  \sum_{\substack{|a|+|b|\le 2N}} \ip{x_0}^{|a|+|b|-|\beta| -n/2},
\end{align*}
and on resolving the summation we thus we arrive at the estimate
\begin{align}\label{eq:CNproof0}
 \Big\Vert \LL^N \big[(\cdot-x_0)^\beta \chi(\cdot)\big] \Big\Vert_{L^2(2B)}
 \lesssim \ip{x_0}^{2N-|\beta|-n/2}.
\end{align}
We then have
\begin{align*}
\Big(\sum_{\xi\in \I_j}\Big|\int \chi(y)(y-x_0)^\beta h_\xi(y)\,dy\Big|^2\Big)^\half
&\lesssim  \ip{x_0}^{2N-|\beta|-n/2} \Big(\sum_{\xi\in I_j}\ip{\xi}^{-2N} \Vert h_\xi\Vert_{L^2(2B)}^2\Big)^{\half}  \\
&\lesssim \ip{x_0}^{2N-|\beta|-n/2} 4^{-jN}\Big(\sum_{\xi\in I_j} \Vert h_\xi\Vert_{L^2(2B)}^2\Big)^{\half}.
\end{align*}
Now note that  $I_j\subset \{\xi\in\NN_0^n: |\xi|\le 4^j\}$. Then by \eqref{eq:QQ est}   we have
$$ \sum_{\xi\in I_j} \Vert h_\xi\Vert_{L^2(2B)}^2 \le \sum_{\xi\le 4^j} \Vert h_\xi\Vert_{L^2(2B)}^2 = \int_{2B} \QQ_{4^j}(y,y)\,dy \lesi 2^{jn}|2B|\sim 2^{jn}\cro(x_0)^n.$$
This gives
\begin{align*}
\Big(\sum_{\xi\in \I_j}\Big|\int \chi(y)(y-x_0)^\beta h_\xi(y)\,dy\Big|^2\Big)^\half
\lesssim \cro(x_0)^{|\beta|+n-2N}2^{-j(2N-n/2)}.
\end{align*}

We next prove \eqref{eq:CN1}. The approach is similar to the proof of \eqref{eq:CN0}; however the estimates require further analysis because of the presence of $\sigma$. Note that one can find a similar  calculation  in \cite[Theorem 4.4(b)]{LN21}.

For $N \in \NN_0$, we have
\begin{align}\label{eq:CNproof1}
&\Big| \int_{\RR^n} \chi(y) (y-x_0)^\beta \,\overline{ \sigma(y,\xi) }\,h_\xi(y)\,dy\Big| \\\notag
&\qquad\qquad=(2|\xi|+n)^{-N} \Big|\int_{\RR^n} \LL_y^N\big [ \chi(y) (y-x_0)^\beta \sigma(y,\xi)\big] h_\xi(y)\,dy\Big| \\\notag
&\qquad\qquad\lesi \ip{\xi}^{-N} \Big\Vert \LL^N \big[\chi(\cdot) (\cdot-x_0)^\beta \sigma(\cdot,\xi)\big] \Big\Vert_{L^2(2B)} \Vert h_\xi\Vert_{L^2(B)}.
\end{align}
Leibniz' rule gives
\begin{align*}
\LL^N \big[\chi(\cdot) (\cdot-x_0)^\beta \sigma(\cdot,\xi)\big](y)
= \sum_{a,b,\gamma,\nu,\eta,\beta} C_{a,b,\gamma,\eta,\nu} \,y^a (y-x_0)^{\beta-\gamma}\chi^{(\eta)}(y) \partial_x^\nu \sigma(y,\xi)
\end{align*}
with the summation running over $|a|+|b|\le 2N$, $\gamma+\nu+\eta=b$ and $|\gamma| \le |\beta|$.
Our condition on $\sigma$ gives, for every $|\nu|\le \N$,
\begin{align*}
\Big(\int_{2B}|\partial_y^\nu \sigma(y,\xi)|^2\,dy\Big)^{\half} 
\lesssim \ip{\xi}^{\dd|\nu|/2} |2B|^\half
\sim \ip{\xi}^{\dd|\nu|/2} \cro(x_0)^{n/2}.
\end{align*}
These two facts give
\begin{align*}
 \Big\Vert \LL^N \big[\chi(\cdot)(\cdot-x_0)^\beta \sigma(\cdot,\xi)\big] \Big\Vert_{L^2(2B)}
 &\lesssim   \sum_{a,b,\gamma,\nu,\eta,\beta}  \cro(x_0)^{|\beta|-|\gamma|-|\eta|} \sup_{y\in 2B}|y|^{|a|} \Big(\int_{2B}|\partial_y^\nu \sigma(y,\xi)|^2\,dy\Big)^{\frac{1}{2}} \\
 &\lesssim  \sum_{a,b,\gamma,\nu,\eta,\beta} \cro(x_0)^{|\beta|-|\gamma| - |\eta|+n/2} \ip{x_0}^{|a|} \ip{\xi}^{\dd|\nu|/2}\\
 &\le  \sum_{a,b,\gamma,\nu,\eta,\beta}\cro(x_0)^{|\beta|-|\gamma| - |\eta|-|\nu|+n/2} \ip{x_0}^{|a|} \Big(\f{\ip{\xi}^\half}{\ip{x_0}}\Big)^{\dd|\nu|} \\
 &\lesi \ip{x_0}^{2N-|\beta|-n/2} \max\Big\{1,\f{\ip{\xi}^\half}{\ip{x_0}}\Big\}^{2N\dd}.
\end{align*}
Inserting the preceding estimate into \eqref{eq:CNproof1} we arrive at 
\begin{align*}
\Big|\int \chi(y)(y-x_0)^\beta\,\overline{\sigma(y,\xi)}\,h_\xi(y)\,dy\Big|
\lesi \ip{\xi}^{-N}  \ip{x_0}^{2N-|\beta|-n/2} \max\Big\{1,\f{\ip{\xi}^\half}{\ip{x_0}}\Big\}^{2N\dd}\Vert h_\xi \Vert_{L^2(2B)}.
\end{align*}
With this latter estimate in hand, we can proceed with
\begin{align*}
&\Big(\sum_{\xi\in \I_j}\Big|\int \chi(y)(y-x_0)^\beta h_\xi(y)\overline{\sigma(y,\xi)}\,dy\Big|^2\Big)^\half \\
&\qquad\lesi   \ip{x_0}^{2N-|\beta|-n/2} \Big(\sum_{\xi\in\I_j}\Big|\ip{\xi}^{-N}\max\big\{1,\ip{\xi}^\half \ip{x_0}^{-1}\big\}^{2N\delta}\Big|^2 \Vert h_\xi\Vert^2_{L^2(B)}\Big)^\half\\
&\qquad\sim 2^{-2Nj} \max\big\{1,2^j\ip{x_0}^{-1}\big\}^{2N\delta}\ip{x_0}^{2N-|\beta|-n/2}\Big(\sum_{\xi\in\I_j}\Vert h_\xi\Vert_{L^2(B)}^2 \Big)^\half \\
&\qquad\sim 2^{-j(2N-n/2)}\max\{1,2^j\cro(x_0)\}^{2N\delta} \cro(x_0)^{|\beta|+n-2N},
\end{align*}
which completes the proof of \eqref{eq:CN1}.

We next prove \eqref{eq:CN2}. We shall proceed as above for \eqref{eq:CN1}; note this is similar to the calculation in \cite[Theorem 4.4(a)]{LN21}. Firstly observe that the condition \eqref{eq:conditionCN} implies
$$\Big(\int_{2B}|\partial_y^\nu \sigma(y,\xi)|^2\,dy\Big)^{\frac{1}{2}}  \lesi \ip{x_0}^{|\nu|}=\cro(x_0)^{-|\nu|}.$$
Set $\wt{N}=\floor{(n+M)/2}+1;$ then we have
\begin{align*}
 \Big\Vert \LL^{\wt{N}} \big[\chi(\cdot)(\cdot-x_0)^\beta \sigma(\cdot,\xi)\big] \Big\Vert_{L^2(2B)}
 &\lesssim \sum_{a,b,\gamma,\nu,\eta,\beta}  \cro(x_0)^{|\beta|-|\gamma|-|\eta|} \sup_{y\in 2B}|y|^{|a|} \Big(\int_{2B}|\partial_y^\nu \sigma(y,\xi)|^2\,dy\Big)^{\frac{1}{2}} \\
 &\lesssim \sum_{a,b,\gamma,\nu,\eta,\beta} \cro(x_0)^{|\beta|-|\gamma| - |\eta|-|\nu|+n/2} \ip{x_0}^{|a|} \\
 &\lesi \ip{x_0}^{2\wt{N}-|\beta|-n/2}.
\end{align*}
Here the summation runs over $|a|+|b|\le 2\wt{N}$, $\gamma+\nu+\eta=b$ and $|\gamma| \le |\beta|$.

Consequently, 
\begin{align*}
\Big(\sum_{\xi\in \I_j}\Big|\int \chi(y)(y-x_0)^\beta h_\xi(y)\overline{\sigma(y,\xi)}\,dy\Big|^2\Big)^\half 
&\lesi   \ip{x_0}^{2\wt{N}-|\beta|-n/2} \Big(\sum_{\xi\in\I_j}\ip{\xi}^{-2N} \Vert h_\xi\Vert^2_{L^2(B)}\Big)^\half\\
&\sim 2^{-2\wt{N}j} \ip{x_0}^{2\wt{N}-|\beta|-n/2}\Big(\sum_{\xi\in\I_j}\Vert h_\xi\Vert_{L^2(B)}^2 \Big)^\half \\
&\sim 2^{-j(2\wt{N}-n/2)}\cro(x_0)^{|\beta|+n-2\wt{N}},
\end{align*}
which yields \eqref{eq:CN2}. 

We now prove \eqref{eq:CN3}. The procedure should now be familiar and we only give the main calculations.
Firstly, 
\begin{align*}
\Big|\int \chi(y)(y-x_0)^\beta \AH^\alpha h_\xi(y)\,dy\Big|
&=|a_\alpha(\xi)|\Big|\int \chi(y)(y-x_0)^\beta  h_{\xi+\wt{\alpha}}(y)\,dy\Big| \\
&= \lambda_{|\xi+\wt{\alpha}|}^{-N}|a_\alpha(\xi)|\Big|\int \LL_y^N\big[\chi(y)(y-x_0)^\beta\big]  h_{\xi+\wt{\alpha}}(y)\,dy\Big| \\
&\le \lambda_{|\xi+\wt{\alpha}|}^{-N}|a_\alpha(\xi)\big\Vert \LL_y^N\big[\chi(\cdot)(\cdot-x_0)^\beta\big] \big\Vert_{L^2(2B)} \Vert h_{\xi+\wt{\alpha}}\Vert_{L^2(2B)}.
\end{align*}
Applying estimate \eqref{eq:CNproof0}, we then have
\begin{align}\label{eq:CNproof2}
\Big|\int \chi(y)(y-x_0)^\beta \AH^\alpha h_\xi(y)\,dy\Big|
\lesi \lambda_{|\xi+\wt{\alpha}|}^{-N}|a_\alpha(\xi)|  \cro(x_0)^{|\beta|+n/2-2N}\Vert h_{\xi+\wt{\alpha}}\Vert_{L^2(2B)}.
\end{align}
Now noting that
\begin{align*}
\sum_{\xi\in\I_j}\Vert h_{\xi+\wt{\alpha}}\Vert_{L^2(2B)}^2
\le \sum_{|\xi|\le 4^j} \Vert h_{\xi+\wt{\alpha}}\Vert_{L^2(2B)}^2
\lesi \sum_{|\xi|\le 4^j+|\alpha|} \Vert h_{\xi}\Vert_{L^2(2B)}^2,
\end{align*}
it follows then that
\begin{align}\label{eq:CNproof3}
\sum_{\xi\in\I_j}\Vert h_{\xi+\wt{\alpha}}\Vert_{L^2(2B)}^2
\lesi \int_{2B} \QQ_{4^j+|\wt{\alpha}|}(y,y)\,dy
\lesi 2^{jn}|2B|
\sim 2^{jn}\cro(x_0)^n.
\end{align}
Estimates \eqref{eq:CNproof2} and \eqref{eq:CNproof3} then allow us to conclude that
\begin{align*}
\Big(\sum_{\xi\in \I_j} \Big|\int \chi(y)(y-x_0)^\beta \AH^\alpha h_\xi(y)\,dy\Big|^2\Big)^\half 
&\lesi \Big(\sum_{\xi\in \I_j}\Big[\lambda_{|\xi+\wt{\alpha}|}^{-N}|a_\alpha(\xi)|   \cro(x_0)^{|\beta|+n/2-2N}\Vert h_{\xi+\wt{\alpha}}\Vert_{L^2(2B)}  \Big]^2\Big)^\half  \\
&\lesi 2^{-j(2N-|\alpha|)}  \cro(x_0)^{|\beta|+n/2-2N} \Big(\sum_{\xi\in \I_j}\Vert h_{\xi+\wt{\alpha}}\Vert_{L^2(2B)}^2\Big)^\half \\
&\lesi 2^{-j(2N-|\alpha|-n/2)}  \cro(x_0)^{|\beta|+n-2N},
\end{align*}
which gives \eqref{eq:CN3} and completes the proof of Lemma \ref{lem:CN}.
\end{proof}

\section{Proofs of some technical results}\label{app:technical}
In this section we give the proofs of some items from Section \ref{sec:applications} that were postponed due to their technicality, namely, Lemmas \ref{lem:ddKj} and \ref{lem:RTddKj}. 

\subsection{Proof of Lemma \ref{lem:ddKj}}
We shall adapt the argument  employed in \cite{LN21}.

We shall consider two cases. 

\underline{Case 1: $N=0$}.
\medskip
By Leibniz' rule,
\begin{align*}
&\D_x^\gamma \D_y^\eta K^\sigma_j(x,y)
= \sum_{\xi\in\I_j} 
\sum_{\nu\le\gamma} \tbinom{\alpha}{\nu} \partial_x^\nu \sigma_j(x,\xi) 
\D_x^{\gamma-\nu} h_{\xi}(x) \, \partial_y^\eta h_\xi(y).
\end{align*}
By the assumption on $\sigma$ and the Cauchy--Schwarz inequality, we get
\begin{align*}
\big|\partial_x^\gamma\partial_y^\eta K^\sigma_j(x,y)\big|
&\le  \sum_{\xi\in\I_j} \sum_{\nu\le \gamma} \tbinom{\alpha}{\nu} |\vph_j(\xi)| |\D_x^\nu \sigma(x,\xi)|\,|\D_x^{\gamma-\nu}h_\xi(x)|\big|\partial_y^\eta h_\xi(y)\big|\\
&\lesi \sum_{\xi\in\I_j} \sum_{\nu\le \gamma}\ip{\xi}^{m/2+|\nu|\dd/2} |\D_x^{\gamma-\nu}h_\xi(x)|\big|\partial_y^\eta h_\xi(y)\big|\\
&\lesi 2^{jm}\sum_{\xi\in\I_j} \sum_{\nu\le \gamma} \ip{\xi}^{\f{\dd}{2}|\nu|}|\D_x^{\gamma-\nu}h_\xi(x)|\big|\partial_y^\eta h_\xi(y)\big|\\
&\le 2^{jm} \sum_{\nu\le \gamma}\Big(\sum_{\xi\in\I_j}  \ip{\xi}^{\dd|\nu|}|\D_x^{\gamma-\nu}h_\xi(x)|^2\Big)^{\half}\Big(\sum_{\xi\in\I_j}\big|\partial_y^\eta h_\xi(y)\big|^2\Big)^\half.
\end{align*}
Now note that 
$$ \big|\partial_y^\eta h_\xi(y)\big| \lesi \ip{\xi}^{|\eta|/2}|h_{\xi+\nu}(y)|.$$
This estimate along with  the fact that $\dd\le 1$ gives
\begin{align*}
\big|\partial_x^\gamma\partial_y^\eta K^\sigma_j(x,y)\big| 
&\lesi 2^{j(m+|\gamma|+|\eta|)}\sum_{\nu\le \gamma}\Big(\sum_{\xi\in\I_j}  h_{\xi+\gamma-\nu}(x)^2\Big)^{\half}\Big(\sum_{\xi\in\I_j} h_{\xi+\eta}(y)^2\Big)^\half \\
&\lesi 2^{j(m+|\gamma|+|\eta|)}\QQ_{4^j+|\gamma|}(x,x)^{\half} \QQ_{4^j+|\eta|}(y,y)^{\half},
\end{align*}
with constants depending on $n,N,\vartheta,\nu,\gamma$. The result for $N=0$ now follows by applying \eqref{eq:QQ est}.

\underline{Case 2: $N\ge 1$}.
\medskip
Note first that we can represent $\partial_x^\gamma$ and $\D_y^\eta$ by
\begin{align*}
\partial_x^\gamma =\sum_{\alpha+\beta\le\gamma} C_{\alpha,\beta} (\A{x})^{\alpha} x^{\beta} &&\text{and}&& \D_y^\eta =\sum_{\wt{\alpha}+\wt{\beta}\le \eta} C_{\wt{\alpha},\wt{\beta}} (\A{y})^{\wt{\alpha}}y^{\wt{\beta}}.
\end{align*}
See (6.15) of \cite{PX}. On combining these expressions  one has
\begin{align}\label{eq:ddrep}
\D_x^\gamma \D_y^\eta = \sum_{\substack{\alpha+\beta\le \gamma \\ \wt{\alpha}+\wt{\beta}\le \eta}} C_{\alpha,\beta,\wt{\alpha},\wt{\beta}} \;(\A{x})^\alpha x^\beta (\A{y})^{\wt{\alpha}} y^{\wt{\beta}}.
\end{align}
This means that in order to prove \eqref{eq:Kjab} it suffices to show
\begin{align}\label{eq:ddKj0}
|x-y|^{N}\big|(\A{x})^{\alpha} x^{\beta}(\A{y})^{\wt{\alpha}}y^{\wt{\beta}} K^\sigma_j(x,y)\big|
\lesssim 2^{j(n+m+|\gamma|+|\eta|+N(1-2\rr))}\big(1+\tfrac{2^{-j}}{\cro(x)}+ \tfrac{2^{-j}}{\cro(y)}\big)^{-\mu}
\end{align}
for any $\alpha,\beta,\wt{\alpha},\wt{\beta}\in \NN_0^n$ such that $0\le |\alpha| + |\beta|\le |\gamma|$, $|\wt{\alpha}|+|\wt{\beta}|\le |\eta|$ and $0\le N\le\K$. 

We first prove bounds for each component $i=1,\dots,n$ by expressing the operators $(\A{x})^\alpha x^\beta$ and $(\A{y})^{\wt{\alpha}}y^{\wt{\beta}}$ in terms of two commuting operators:
\begin{align*}
\big(\A{x}\big)^\alpha x^\beta =\big(\A{x}\big)^{\alpha-\alpha_ie_i}x^{\beta-\beta_ie_i}\big(\A{x}_i\big)^{\alpha_i}x_i^{\beta_i}, \\
\big(\A{y}\big)^{\wt{\alpha}} y^{\wt{\beta}} =\big(\A{y}\big)^{\wt{\alpha}-\wt{\alpha}_ie_i}y^{\wt{\beta}-\wt{\beta}_ie_i}\big(\A{y}_i\big)^{\wt{\alpha}_i}y_i^{\wt{\beta}_i}, 
\end{align*}
where $\{e_i\}_{1\le i\le n}$ is the canonical basis for $\RR^n$. Then by the above expression, identity \eqref{eq:identity b2} and then recombining the commuting operators we have
\begin{align}\label{eq:ddKj1a}
(x_i-y_i)^{N}(\A{x})^{\alpha} x^{\beta}(\A{y})^{\wt{\alpha}}y^{\wt{\beta}} 
= \sum_{s=0}^{\alpha_i}\tbinom{\alpha_i}{s}\tfrac{N!}{(N-s)!}\big(\A{x}\big)^{\alpha-se_i}x^\beta (x_i-y_i)^{N-s}(\A{y})^{\wt{\alpha}}y^{\wt{\beta}}
\end{align}
where $\tfrac{N!}{(N-s)!}:=0$ if $s>N$. Applying the same idea for $(\A{y})^{\wt{\alpha}}y^{\wt{\beta}}$ gives
\begin{align}\label{eq:ddKj1b}
(x_i-y_i)^{N-s}(\A{y})^{\wt{\alpha}}y^{\wt{\beta}}
=\sum_{\wt{s}=0}^{\wt{\alpha}_i} (-1)^{\wt{s}}\tbinom{\wt{\alpha}_i}{\wt{s}}\tfrac{(N-s)!}{(N-s-\wt{s})!}(\A{y})^{\wt{\alpha}-\wt{s}e_i}y^{\wt{\beta}}(x_i-y_i)^{N-s-\wt{s}}.
\end{align}
Combining \eqref{eq:ddKj1a}-\eqref{eq:ddKj1b} we have
\begin{align}\label{eq:ddKj2}
&(x_i-y_i)^{N}(\A{x})^{\alpha} x^{\beta}(\A{y})^{\wt{\alpha}}y^{\wt{\beta}} K^\sigma_j(x,y) \\
&\quad= \sum_{\substack{0\le s\le \alpha_i\\ 0\le \wt{s}\le \wt{\alpha}_i}}(-1)^{\wt{s}}\tbinom{\alpha_i}{s}\tbinom{\wt{\alpha}_i}{\wt{s}}\tfrac{N!}{(N-s-\wt{s})!} \big(\A{x}\big)^{\alpha-se_i}x^\beta (\A{y})^{\wt{\alpha}-\wt{s}e_i}y^{\wt{\beta}}(x_i-y_i)^{N-s-\wt{s}} K^\sigma_j(x,y) 
\notag
\end{align}
with the understanding that the terms in the summation vanish whenever $s+\wt{s}>N$. Then whenever $s+\wt{s}\le N$ we have by \eqref{eq:identity A},
\begin{align*}
 &(x_i-y_i)^{N-s-\wt{s}}K^\sigma_j(x,y)\\
 &\quad=2^{-(N-s-\wt{s})} \sum_{\f{N-s-\wt{s}}{2}\le\ell\le N-s-\wt{s}} c_{\ell, N-s-\wt{s}} \sum_{\xi\in\I_j} \diff_i^\ell \sigma_j(x,\xi) \big(\A{y}_i-\A{x}_i\big)^{2\ell-N+s+\wt{s}}\big[ h_\xi(x)\, h_\xi(y)\big].
\end{align*}
We  next apply the operator $(\A{y})^{\wt{\alpha}-\wt{s}e_i}y^{\wt{\beta}}$ to this expression. Writing $y^{\wt{\beta}}=y^{\wt{\beta}-\wt{\beta} e_i}y_i^{\wt{\beta}_i},$ applying \eqref{eq:identity b1} to $y_i^{\wt{\beta}_i}\big(\A{y}_i-\A{x}_i\big)^{2\ell-N+s+\wt{s}}$, and commuting $y_i^{\wt{\beta}-\wt{\beta}_i e_i}$ with powers of $\A{y}_i-\A{x}_i,$ we obtain
\begin{align*}
&(\A{y})^{\wt{\alpha}-\wt{s}e_i}y^{\wt{\beta}}(x_i-y_i)^{N-s-\wt{s}} K^\sigma_j(x,y)\\
&\quad=2^{-(N-s-\wt{s})} \sum_{\f{N-s-\wt{s}}{2}\le\ell\le N-s-\wt{s}} c_{\ell, N-s-\wt{s}} \sum_{\xi\in\I_j} \diff_i^\ell \sigma_j(x,\xi) \\
&\quad\quad\quad\times (\A{y})^{\wt{\alpha}-\wt{s}e_i}y^{\wt{\beta}} \big(\A{y}_i-\A{x}_i\big)^{2\ell-N+s+\wt{s}}\big[ h_\xi(x)\, h_\xi(y)\big]\\
&\quad=2^{-(N-s-\wt{s})} \sum_{\f{N-s-\wt{s}}{2}\le\ell\le N-s-\wt{s}} c_{\ell, N-s-\wt{s}} \sum_{\xi\in\I_j} \diff_i^\ell \sigma_j(x,\xi) \\
&\quad\quad\quad\times (\A{y})^{\wt{\alpha}-\wt{s}e_i}y^{\wt{\beta}-\wt{\beta} e_i} \sum_{\wt{t}=0}^{\wt{\beta}_i} \tbinom{\wt{\beta}_i}{\wt{t}}\tfrac{(2\ell-N+s+\wt{s})!}{(2\ell-N+s+\wt{s}-\wt{t})!} \big(\A{y}_i-\A{x}_i\big)^{2\ell-N+s+\wt{s}-\wt{t}}y_i^{\wt{\beta}_i-\wt{t}}\big[ h_\xi(x)\, h_\xi(y)\big] \\
&\quad=2^{-(N-s-\wt{s})} \sum_{\f{N-s-\wt{s}}{2}\le\ell\le N-s-\wt{s}} c_{\ell, N-s-\wt{s}} \sum_{\xi\in\I_j} \diff_i^\ell \sigma_j(x,\xi) \\
&\quad\quad\quad\times  \sum_{\wt{t}=0}^{\wt{\beta}_i} \tbinom{\wt{\beta}_i}{\wt{t}}\tfrac{(2\ell-N+s+\wt{s})!}{(2\ell-N+s+\wt{s}-\wt{t})!} \big(\A{y}_i-\A{x}_i\big)^{2\ell-N+s+\wt{s}-\wt{t}}(\A{y})^{\wt{\alpha}-\wt{s}e_i}y^{\wt{\beta}-\wt{t}e_i}\big[ h_\xi(x)\, h_\xi(y)\big].
\end{align*}
We next apply $\big(\A{x}\big)^{\alpha-se_i}x^\beta$ to this expression  in a similar way.  First we invoke  Leibniz' rule for hermite derivatives \eqref{eq:leibniz2}. This produces
\begin{align*}
&\big(\A{x}\big)^{\alpha-se_i}x^\beta(\A{y})^{\wt{\alpha}-\wt{s}e_i}y^{\wt{\beta}}(x_i-y_i)^{N-s-\wt{s}} K^\sigma_j(x,y)\\
&\hspace{20pt}=2^{-(N-s-\wt{s})} \sum_{\f{N-s-\wt{s}}{2}\le\ell\le N-s-\wt{s}} \sum_{\xi\in\I_j} \sum_{\wt{t}=0}^{\wt{\beta}_i}c_{\ell, N-s-\wt{s}} \tbinom{\wt{\beta}_i}{\wt{t}}\tfrac{(2\ell-N+s+\wt{s})!}{(2\ell-N+s+\wt{s}-\wt{t})!}\\
&\hspace{50pt}\times \big(\A{x}\big)^{\alpha-se_i}\Big\{ \diff_i^\ell\sigma_j(x,\xi)  x^\beta \big(\A{y}_i-\A{x}_i\big)^{2\ell-N+s+\wt{s}-\wt{t}}(\A{y})^{\wt{\alpha}-\wt{s}e_i}y^{\wt{\beta}-\wt{t}e_i}\big[ h_\xi(x)\, h_\xi(y)\big] \Big\} \\
&\hspace{20pt}=2^{-(N-s-\wt{s})} \sum_{\f{N-s-\wt{s}}{2}\le\ell\le N-s-\wt{s}} \sum_{\xi\in\I_j} \sum_{\wt{t}=0}^{\wt{\beta}_i}c_{\ell, N-s-\wt{s}}  \tbinom{\wt{\beta}_i}{\wt{t}}\tfrac{(2\ell-N+s+\wt{s})!}{(2\ell-N+s+\wt{s}-\wt{t})!}\\
&\hspace{50pt}\times \sum_{\nu\le \alpha-se_i} (-1)^\nu\tbinom{\alpha-se_i}{\nu} \D_x^\nu\diff_i^\ell\sigma_j(x,\xi)  \\
&\hspace{80pt}\times \big(\A{x}\big)^{\alpha-se_i-\nu}\Big\{x^\beta \big(\A{y}_i-\A{x}_i\big)^{2\ell-N+s+\wt{s}-\wt{t}}(\A{y})^{\wt{\alpha}-\wt{s}e_i}y^{\wt{\beta}-\wt{t}e_i}\big[ h_\xi(x)\, h_\xi(y)\big] \Big\}.
\end{align*}
Now applying \eqref{eq:identity b1}, expanding powers of $(\A{x}_i-\A{y}_i)$ by the binomial theorem, and then absorbing the polynomials $x^{\beta-t e_i}$ and $y^{\wt{\beta}-\wt{t}e_i}$ into $h_\xi(x)h_\xi(y)$ via \eqref{eq:identity C} we obtain
\begin{align}\label{eq:ddKj3}
&\big(\A{x}\big)^{\alpha-se_i}x^\beta(\A{y})^{\wt{\alpha}-\wt{s}e_i}y^{\wt{\beta}}(x_i-y_i)^{N-s-\wt{s}} K^\sigma_j(x,y)\\\notag
&\quad=2^{-(N-s-\wt{s})} \sum_{\f{N-s-\wt{s}}{2}\le\ell\le N-s-\wt{s}} \sum_{\xi\in\I_j} \sum_{\wt{t}=0}^{\wt{\beta}_i}\sum_{t=0}^{\beta_i} \sum_{\nu\le \alpha-se_i}c_{\ell, N-s-\wt{s}}  \tbinom{\wt{\beta}_i}{\wt{t}}\tbinom{\beta_i}{t}\tbinom{\alpha-se_i}{\nu}  \tfrac{(2\ell-N+s+\wt{s})!}{(2\ell-N+s+\wt{s}-t-\wt{t})!} (-1)^{t+\nu}\\\notag
&\quad\quad\times  \D_x^\nu\diff_i^\ell\sigma_j(x,\xi)  \big(\A{y}_i-\A{x}_i\big)^{2\ell-N+s+\wt{s}-t-\wt{t}}\big(\A{x}\big)^{\alpha-se_i-\nu} (\A{y})^{\wt{\alpha}-\wt{s}e_i} x^{\beta-te_i}h_\xi(x)\,y^{\wt{\beta}-\wt{t}e_i} h_\xi(y)\\\notag
&\quad=2^{-(N-s-\wt{s})} \sum_{\f{N-s-\wt{s}}{2}\le\ell\le N-s-\wt{s}} \sum_{\xi\in\I_j} \sum_{\wt{t}=0}^{\wt{\beta}_i}\sum_{t=0}^{\beta_i} \sum_{\nu\le \alpha-se_i}c_{\ell, N-s-\wt{s}}  \tbinom{\wt{\beta}_i}{\wt{t}}\tbinom{\beta_i}{t}\tbinom{\alpha-se_i}{\nu}  \tfrac{(2\ell-N+s+\wt{s})!}{(2\ell-N+s+\wt{s}-t-\wt{t})!} (-1)^{t+\nu}\\\notag
&\quad\quad\times  \D_x^\nu\diff_i^\ell\sigma_j(x,\xi)  
\sum_{\substack{k,\wt{k}\ge 0\\ k+\wt{k}=2\ell-N+s+\wt{s}-t-\wt{t}}} \sum_{\om\le \beta-te_i} \sum_{\wt{\om}\le \wt{\beta}-\wt{t}e_i}
(-1)^k \tbinom{2\ell-N+s+\wt{s}-t-\wt{t}}{k,\;\wt{k}}  b_{\om,\beta-te_i}(\xi) b_{\wt{\om},\wt{\beta}-\wt{t}e_i}(\xi) \\\notag
&\quad\quad\quad\times\big(\A{x}\big)^{\alpha-se_i-\nu+ke_i} h_{\xi+\beta-te_i-2\om}(x) (\A{y})^{\wt{\alpha}-\wt{s}e_i+\wt{k}e_i} h_{\xi+\wt{\beta}-\wt{t}e_i-2\wt{\om}}(y). 
\end{align}
Then from \eqref{eq:ddKj2} and \eqref{eq:ddKj3} we get
\begin{align}\label{eq:ddKj4}
&(x_i-y_i)^{N}(\A{x})^{\alpha} x^{\beta}(\A{y})^{\wt{\alpha}}y^{\wt{\beta}} K^\sigma_j(x,y) \\\notag
&\qquad= \sum_{\substack{s,\wt{s},\xi,\ell,t,\wt{t}\\ \nu,\om,\wt{\om},k,\wt{k}}} \wt{C} \; b_{\om,\beta-te_i}(\xi) \; b_{\wt{\om},\wt{\beta}-\wt{t}e_i}(\xi)  \; \D_x^\nu\diff_i^\ell\sigma_j(x,\xi)  \\\notag
&\qquad\qquad\qquad\times\big(\A{x}\big)^{\alpha-se_i-\nu+ke_i} h_{\xi+\beta-te_i-2\om}(x) \;(\A{y})^{\wt{\alpha}-\wt{s}e_i+\wt{k}e_i} h_{\xi+\wt{\beta}-\wt{t}e_i-2\wt{\om}}(y),
\end{align}
where the summation is over
\begin{align*}
&0\le s\le \alpha_i, && 0\le t\le \beta_i, && \om\le \beta-te_i,  \\
&0\le \wt{s}\le \wt{\alpha}_i, && 0\le \wt{t}\le \wt{\beta}_i, && \wt{\om}\le \wt{\beta}-\wt{t}e_i\\ 
&\xi\in\I_j, && (N-s-\wt{s})/2\le \ell \le N-s-\wt{s}, \\
& k,\wt{k}\ge 0, && k+\wt{k}=2\ell-N+s+\wt{s}-t-\wt{t}
\end{align*}
and 
\begin{align*}
\wt{C}&=(-1)^{t+k+\wt{s}+\nu}2^{-(N-s-\wt{s})}\tbinom{\alpha_i}{s}\tbinom{\wt{\alpha}_i}{\wt{s}} \tbinom{\wt{\beta}_i}{\wt{t}}\tbinom{\beta_i}{t} \\
&\hspace{50pt}\times \tbinom{\alpha-se_i}{\nu}  \tbinom{2\ell-N+s+\wt{s}-t-\wt{t}}{k,\;\wt{k}} \tfrac{N!}{(N-s-\wt{s})!}\tfrac{(2\ell-N+s+\wt{s})!}{(2\ell-N+s+\wt{s}-t-\wt{t})!} c_{\ell, N-s-\wt{s}}. 
\end{align*}
We now proceed to estimate \eqref{eq:ddKj4}. From \eqref{eq:identity C} we have
\begin{align}\label{eq:ddKj5a}
|b_{\om,\beta-te_i}(\xi) \; b_{\wt{\om},\wt{\beta}-\wt{t}e_i}(\xi)|\lesi \ip{\xi}^{\f{|\beta|+|\wt{\beta}|-t-\wt{t}}{2}}.
\end{align}
Firstly, from the Leibniz formula for finite differences \eqref{eq:leibniz1} we have
\begin{align*}
|\D_x^\nu\diff_i^\ell\sigma_j(x,\xi)| = \big|\D_x^\nu\diff_i^\ell \big[\sigma(x,\xi)\vp_j(\xi)\big]\big| 
\le \sum_{r=0}^\ell \tbinom{\ell}{r} \big|\diff_i^r(\vp_j(\xi))\big| \, \big|\D_x^\nu\diff_i^{\ell-r}\sigma(x,\xi+re_i)\big|.
\end{align*}
Next, to estimate the terms in the sum, we shall make use of the following Lemma. 
\begin{Lemma}[Lemma 2.2 in \cite{LN21}]\label{lem: hoppe}
Let $\phi$ be a smooth function defined in $[0,\infty),$ set $\phi_j(x)=\phi(2^{-j}x)$ for $j\in \NN_0$ and consider  $\ell, N\in \NN,$ $N\ge \ell.$
If $\phi^{(m)}(0)=0$ for all $m\in \NN,$ it holds that 
\begin{equation*}
 |\diff^\ell_\xi(\phi_j(\sqrt{2|\xi|+n}))|\lesssim   \|\phi^{(N)}\|_{L^\infty}  2^{-jN} \ip{\xi}^{N/2-\ell}\quad \forall j,\xi\in \NN_0,
 \end{equation*}
 where the implicit constant depends on $N$, $n$ and $\ell.$
\end{Lemma}
With Lemma~\ref{lem: hoppe} in hand, we have
\begin{align*}
\big|\diff_i^r(\vp_j(\xi))\big| \lesssim  2^{-jN}\ip{\xi}^{N/2-r},
\end{align*}
while our assumption on $\sigma$ gives
\begin{align*}
\big|\D_x^\nu\diff_i^{\ell-r}\sigma(x,\xi+re_i)\big| 
\lesssim  \ip{\xi}^{\f{m}{2}+\f{\dd}{2}|\nu|-\rr(\ell-r)}.
\end{align*} 
These last three facts, and the fact that $\dd,\rr\le 1$ and $r\ge 0$ give
\begin{align}\label{eq:ddKj5b}
\big|\D_x^\nu\diff_i^\ell \sigma_j(x,\xi)\big|\lesssim 2^{-jN} \ip{\xi}^{\f{N}{2}+\f{m}{2}+\f{|\nu|}{2} -\rr\ell}.
\end{align}
Observe also that, by applying \eqref{eq:identity d2} we have
\begin{align}\label{eq:ddKj5c}
\big|\big(\A{x}\big)^{\alpha-se_i-\nu+ke_i} h_{\xi+\beta-te_i-2\om}(x)\big|  \lesi \ip{\xi}^{|\alpha-se_i-\nu+ke_i|/2} |h_{\xi+\alpha+\beta-\nu-2\om +(k-t-s)e_i}(x)|
\end{align}
and
\begin{align}\label{eq:ddKj5d}
\big|(\A{y})^{\wt{\alpha}-\wt{s}e_i+\wt{k}e_i} h_{\xi+\wt{\beta}-\wt{t}e_i-2\wt{\om}}(y)\big|
\lesi \ip{\xi}^{|\wt{\alpha}-\wt{s}e_i + \wt{k}e_i|/2} |h_{\xi+\wt{\alpha}+\wt{\beta}-2\wt{\om}+(\wt{k}-\wt{t}-\wt{s})e_i}(y)|.
\end{align}
Inserting \eqref{eq:ddKj5a}-\eqref{eq:ddKj5d} into \eqref{eq:ddKj4}, and noting the binomial bounds $\binom{a}{b}\le 2^a$ and $\f{a!}{(a-b)!}\le a^b$, we obtain 
\begin{align*}
&\big|(x_i-y_i)^{N}(\A{x})^{\alpha} x^{\beta}(\A{y})^{\wt{\alpha}}y^{\wt{\beta}} K^\sigma_j(x,y)\big|\\
&\quad\lesi 2^{-jN} \sum_{\substack{s,\wt{s},\xi,\ell,t,\wt{t}\\ \nu,\om,\wt{\om},k,\wt{k}}} 
\ip{\xi}^{\f{m+|\alpha|+|\wt{\alpha}|+|\beta|+|\wt{\beta}|}{2} +\ell(1-\rr)} 
 |h_{\xi+\alpha+\beta-\nu-2\om +(k-t-s)e_i}(x)|
|h_{\xi+\wt{\alpha}+\wt{\beta}-2\wt{\om}+(\wt{k}-\wt{t}-\wt{s})e_i}(y)|.
\end{align*}
Now since $\rr\le 1$, $\ell\le N$ and recalling that $\ip{\xi}\sim 2^{j}$ we have
\begin{align*}
\ip{\xi}^{\f{m+|\alpha|+|\wt{\alpha}|+|\beta|+|\wt{\beta}|}{2} +\ell(1-\rr)} 
 \le  \ip{\xi}^{\f{m+|\alpha|+|\wt{\alpha}|+|\beta|+|\wt{\beta}|}{2}+N(1-\rr)}
 \lesi 2^{j(m+|\alpha|+|\wt{\alpha}|+|\beta|+|\wt{\beta}|+2N(1-\rr))}.
\end{align*}
Using this estimate along with two applications of the Cauchy--Schwarz inequality, it follows that
\begin{align*}
&\big|(x_i-y_i)^{N}(\A{x})^{\alpha} x^{\beta}(\A{y})^{\wt{\alpha}}y^{\wt{\beta}} K^\sigma_j(x,y)\big|\\
&\quad\lesi 2^{j(m+|\alpha|+|\wt{\alpha}|+|\beta|+|\wt{\beta}|+N(1-2\rr))} \sum_{\substack{s,\wt{s},\xi,\ell,t,\wt{t}\\ \nu,\om,\wt{\om},k,\wt{k}}} 
 |h_{\xi+\alpha+\beta-\nu-2\om +(k-t-s)e_i}(x)|
|h_{\xi+\wt{\alpha}+\wt{\beta}-2\wt{\om}+(\wt{k}-\wt{t}-\wt{s})e_i}(y)| \\
&\quad\lesi 2^{j(m+|\alpha|+|\wt{\alpha}|+|\beta|+|\wt{\beta}|+N(1-2\rr))}
\QQ_{c(4^j + N+ |\alpha|+|\beta|)}(x,x)^\half \QQ_{c(4^j + N+ |\wt{\alpha}|+|\wt{\beta}|)}(y,y)^\half.
\end{align*}
Estimate \eqref{eq:ddKj0} now follows once we invoke \eqref{eq:QQ est}, concluding the proof of Lemma \ref{lem:ddKj}.

\subsection{Proof of Lemma \ref{lem:RTddKj}}
The proof follows along the same lines as that of Lemma \ref{lem:ddKj} but is slightly simpler. We make some preliminary observations.

Secondly note that the symbol $\sigma_\alpha(\xi):=(2|\xi|+n)^{-|\alpha|/2}$ satisfies $\sigma_\alpha\in \SM^{-|\alpha|,\infty,\infty}_{1,0}$. Indeed by the mean value theorem, we have for some $\nu$ satisfying $|\xi| \le |\nu| \le |\xi|+|\kappa|$,
\begin{align*}
|\diff_\xi^\kappa \sigma_\alpha(\xi)| = |\partial_\nu^\kappa \sigma_\alpha(\nu)| \lesi \ip{\nu}^{-\f{|\alpha|}{2}-|\kappa|} \sim \ip{\xi}^{-\f{|\alpha|}{2}-|\kappa|}.
\end{align*}
We now continue with the proof of the lemma.

\underline{Case 1:} $N=0$.  From the expression
\begin{align}\label{eq:RTddKj1}
\partial_x^\gamma \partial_y^\eta \R_j^\alpha(x,y) 
=\sum_{\xi\in\I_j} \vph_j(\xi)\sigma_\alpha(\xi)\, \partial_x^\gamma \AH^\alpha h_\xi(x) \partial_y^\eta h_\xi(y),
\end{align}
and the estimates
\begin{align*}
\big|\partial_y^\eta h_\xi(y)\big|\lesi \ip{\xi}^{\f{|\eta|}{2}} |h_{\xi+\eta}(y)|
\end{align*}
and
\begin{align*}
\big|\partial_x^\gamma\AH^\alpha h_\xi(x)\big|
= |a_\alpha(\xi)| \big|\partial_x^\gamma h_{\xi+\wt{\alpha}}(x)\big|
\lesi \ip{\xi}^{\f{|\alpha|+|\gamma|}{2}} |h_{\xi+\wt{\alpha}+\gamma}(x)|,
\end{align*}
we obtain
\begin{align*}
\big|\partial_x^\gamma \partial_y^\eta \R_j^\alpha(x,y) \big|
\lesi \sum_{\xi\in\I_j} \ip{\xi}^{\f{|\gamma|+|\eta|}{2}}  |h_{\xi+\wt{\alpha}+\gamma}(x)|\, |h_{\xi+\eta}(y)| 
\lesi 2^{j(|\gamma|+|\eta|)}\sum_{\xi\in\I_j}  |h_{\xi+\wt{\alpha}+\gamma}(x)|\, |h_{\xi+\eta}(y)|.
\end{align*}
Now applying the Cauchy--Schwarz inequality and the bounds \eqref{eq:QQ est} we have
\begin{align*}
\big|\partial_x^\gamma \partial_y^\eta \R_j^\alpha(x,y) \big|
&\lesi 2^{j(|\gamma|+|\eta|)}\QQ_{4^j+|\wt{\alpha}|+|\gamma|}(x,x)^\half \QQ_{4^j+|\eta|}(y,y)^\half \\
&\lesi 2^{j(n+|\gamma|+|\eta|)}\big(1+\tfrac{2^{-j}}{\cro(x)}+\tfrac{2^{-j}}{\cro(y)}\big)^{-\mu}.
\end{align*}

\underline{Case 2:} $N\ge 1$. As in \eqref{eq:ddrep}, one can represent $\D_x^\gamma \D_y^\eta$ via the expression
\begin{align*}
\D_x^\gamma \D_y^\eta = \sum_{\substack{\gamma'+\gamma''\le \gamma \\ \eta'+\eta''\le \eta}} C_{\gamma',\gamma'',\eta',\eta''} \;(A^{(x)})^{\gamma'} x^{\gamma''} (A^{(y)})^{\eta'} y^{\eta''}.
\end{align*}
Thus, in order to prove \eqref{eq:RTKjab} it suffices to show
\begin{align}\label{eq:ddKj0}
|x-y|^{N}\big|(A^{(x)})^{\gamma'} x^{\gamma''} (A^{(y)})^{\eta'} y^{\eta''} \R_j^\alpha(x,y)\big|
\lesssim 2^{j(n+|\gamma|+|\eta|-N)}\big(1+\tfrac{2^{-j}}{\cro(x)}+ \tfrac{2^{-j}}{\cro(y)}\big)^{-\mu}
\end{align}
for any $\gamma',\gamma'',\eta',\eta''\in \NN_0^n$ such that $0\le |\gamma'|+|\gamma''|\le |\gamma|$, $|\eta'|+|\eta''|\le |\eta|$. 
As in \cite{Ly22} we first prove bounds for each component $i=1,\dots,n$ by expressing the operators $(A^{(x)})^{\gamma'} x^{\gamma''}$ and $(A^{(y)})^{\eta'}y^{\eta''}$ in terms of two commuting operators:
\begin{align*}
\big(A^{(x)}\big)^{\gamma'} x^{\gamma''} =\big(A^{(x)}\big)^{\gamma'-\gamma'_ie_i}x^{\gamma''-\gamma''_ie_i}\big(A^{(x)}_i\big)^{\gamma'_i}x_i^{\gamma''_i}, \\
\big(A^{(y)}\big)^{\eta'} y^{\eta''} =\big(A^{(y)}\big)^{\eta'-\eta'_ie_i}y^{\eta''-\eta''_ie_i}\big(A^{(y)}_i\big)^{\eta'_i}y_i^{\eta''_i}, 
\end{align*}
where $\{e_i\}_{1\le i\le n}$ is the canonical basis for $\RR^n$. Then using \eqref{eq:identity b2} with $\AH_i=A_i$ and proceeding as in \cite{Ly22} we have
\begin{align}\label{eq:ddKj2}
&(x_i-y_i)^{N}(A^{(x)})^{\gamma'} x^{\gamma''}(A^{(y)})^{\eta'}y^{\eta''} \R_j^\alpha(x,y) \\
&\quad= \sum_{\substack{0\le s\le \gamma'_i\\ 0\le \wt{s}\le \eta'_i}}(-1)^{\wt{s}}\tbinom{\gamma'_i}{s}\tbinom{\eta'_i}{\wt{s}}\tfrac{N!}{(N-s-\wt{s})!} \big(A^{(x)}\big)^{\gamma'-se_i}x^{\gamma''}(A^{(y)})^{\eta'-\wt{s}e_i}y^{\eta''}(x_i-y_i)^{N-s-\wt{s}} \R_j^\alpha(x,y) 
\notag
\end{align}
Next we use the same idea to commute $(x_i-y_i)^{N-s-\wt{s}}$ with $(\AH^{(x)})^\alpha$. That is, we write $
\AH^{\alpha} = \AH^{\alpha-\alpha_ie_i} \AH_i^{\alpha_i}$ and then apply \eqref{eq:identity b2} to obtain
\begin{align*}
(x_i-y_i)^{N-s-\wt{s}}(\AH^{(x)})^\alpha =\sum_{t=0}^{\alpha_i} c_{\AH_i,t} \tbinom{\alpha_i}{t} \tfrac{(N-s-\wt{s})!}{(N-s-\wt{s}-t)!} (\AH^{(x)})^{\alpha-te_i} (x_i-y_i)^{N-s-\wt{s}-t}
\end{align*}
Inserting this expression into \eqref{eq:ddKj2} we have
\begin{align}\label{eq:ddKj3}
&(x_i-y_i)^{N}(A^{(x)})^{\gamma'} x^{\gamma''}(A^{(y)})^{\eta'}y^{\eta''} \R_j^\alpha(x,y) \\
&\qquad\qquad=
\sum_{s,\wt{s},t}(-1)^{\wt{s}}c_{\AH_i,t}\tbinom{\gamma'_i}{s}\tbinom{\eta'_i}{\wt{s}}\tbinom{\alpha_i}{t}\tfrac{N!}{(N-s-\wt{s}-t)!} \big(A^{(x)}\big)^{\gamma'-se_i}x^{\gamma''}(A^{(y)})^{\eta'-\wt{s}e_i}y^{\eta''}  \notag\\
&\qquad\qquad\qquad\qquad \times (A^{(x)})^{\alpha-te_i} 
 (x_i-y_i)^{N-s-\wt{s}-t} \sum_{\xi\in\I_j} \sigma_\alpha(\xi)\vph_j(\xi)h_\xi(x)h_\xi(y),
\notag
\end{align}
with the understanding that the terms in the summation vanish whenever $s+\wt{s}+t>N$; note also the summation here runs over $0\le s\le \gamma'_i$,  $0\le \wt{s}\le \eta'_i$, and $0\le t\le \alpha_i$.
Applying \eqref{eq:identity A} to above expression and then using the binomial theorem, we have
\begin{align*}
&(x_i-y_i)^{N-s-\wt{s}-t} \sum_{\xi\in\I_j} \sigma_\alpha(\xi)\vph_j(\xi)h_\xi(x)h_\xi(y)\\
&\qquad\qquad= 2^{-\wt{N}}\sum_{\wt{N}/2\le \ell\le \wt{N}} c_{\ell,\wt{N}}\sum_{\xi\in\I_j} \diff_i^\ell\big[\vph_j(\xi)\sigma_\alpha(\xi)\big] \big(A_i^{(y)}-A_i^{(x)}\big)^{2\ell-\wt{N}}h_\xi(x)h_\xi(y) \\
&\qquad\qquad= 2^{-\wt{N}}\sum_{\wt{N}/2\le \ell\le \wt{N}} c_{\ell,\wt{N}}\sum_{\xi\in\I_j} \diff_i^\ell\big[\vph_j(\xi)\sigma_\alpha(\xi)\big]  \sum_{\substack{k,\wt{k}\ge 0\\ k+\wt{k}=2\ell-\wt{N}}} \tbinom{2\ell-\wt{N}}{k,\wt{k}} (A_i^{(x)})^k h_\xi(x) (A_i^{(y)})^{\wt{k}}h_\xi(y)
\end{align*}
where $\wt{N}:=N-s-\wt{s}-t$. Now inserting this into \eqref{eq:ddKj3} we arrive at
\begin{align}\label{eq:ddKj4}
&(x_i-y_i)^{N}(A^{(x)})^{\gamma'} x^{\gamma''}(A^{(y)})^{\eta'}y^{\eta''} \R_j^\alpha(x,y) \\\notag
&\quad= \sum_{\substack{s,\wt{s},t,\xi,\\ \ell,k,\wt{k}}} \wt{C} \diff_i^\ell\big[\vph_j(\xi)\sigma_\alpha(\xi)  \big] 
\big(A^{(x)}\big)^{\gamma'-se_i} x^{\gamma''} \big(\AH^{(x)}\big)^{\alpha-te_i} \big(A_i^{(x)}\big)^{k} h_{\xi}(x) 
\;(A^{(y)})^{\eta'-\wt{s}e_i} y^{\eta''} \big(A_i^{(y)}\big)^{\wt{k}} h_{\xi}(y)
\end{align}
where the summation is over
\begin{align*}
&0\le s\le \gamma'_i, &&0\le \wt{s}\le \eta'_i,&&  0\le t\le \alpha_i, &&  \xi\in\I_j, \\
& \wt{N}/2\le \ell \le \wt{N}, 
&& k,\wt{k}\ge 0, && k+\wt{k}=2\ell-\wt{N}
\end{align*}
and 
\begin{align*}
\wt{C}=(-1)^{\wt{s}}c_{\AH_i,t} 2^{-\wt{N}}\tbinom{\gamma'_i}{s}\tbinom{\eta'_i}{\wt{s}}\tbinom{\alpha_i}{t}  \tbinom{2\ell-\wt{N}}{k,\;\wt{k}}\tfrac{N!}{\wt{N}!} c_{\ell, \wt{N}}. 
\end{align*}
We now proceed to estimate \eqref{eq:ddKj4}. From \eqref{eq:d1a}, \eqref{eq:RTddKj0} and \eqref{eq:identity C} we have
\begin{align*}
&\big|\big(A^{(x)}\big)^{\gamma'-se_i} x^{\gamma''} \big(\AH^{(x)}\big)^{\alpha-te_i} \big(A_i^{(x)}\big)^{k} h_{\xi}(x)\big| \\
&\qquad\qquad\qquad\lesi \sum_{\om \le \gamma''} \ip{\xi}^{(k+|\gamma'-se_i|+|\gamma''| +|\alpha-te_i|)/2} |h_{\xi +\gamma'+\gamma''-2\om+\wt{\alpha} + (k-s)e_i}(x)|;
\end{align*}
from \eqref{eq:d1a} and \eqref{eq:identity C} we have
\begin{align*}
\big|(A^{(y)})^{\eta'-\wt{s}e_i} y^{\eta''} \big(A_i^{(y)}\big)^{\wt{k}} h_{\xi}(y)\big|
\lesi \sum_{\wt{\om}\le \eta''} \ip{\xi}^{(|\eta'-\wt{s}e_i|+|\eta''|+\wt{k})/2} |h_{\xi+\eta'+\eta''-2\wt{\om} +(\wt{k}-\wt{s})e_i}(y)|.
\end{align*}
From \cite[Lemma 2.2]{LN21} we have the estimate
\begin{align*}
\big|\diff_i^r(\vp_j(\xi))\big| \lesssim 2^{-j\wt{N}}\ip{\xi}^{\wt{N}/2-r},
\end{align*}
and since $\sigma_\alpha\in \SM^{-|\alpha|,\infty,\infty}_{1,0}$ we also have
\begin{align*}
\big|\diff_i^{\ell-r}\sigma_\alpha(\xi)\big| 
\lesi \ip{\xi}^{-|\alpha|/2-\ell+r}.
\end{align*} 
These two estimates along with the Leibniz formula for finite differences gives
\begin{align}\label{eq:RTddKj3}
\big|\diff_i^\ell \big[\sigma_\alpha(\xi)\vph_j(\xi)\big]\big| 
\le \sum_{r=0}^\ell \tbinom{\ell}{r} \big|\diff_i^r\vp_j(\xi)\big| \, \big|\diff_i^{\ell-r}\sigma_\alpha(\xi+re_i)\big|
\lesi 2^{-j\wt{N}}\ip{\xi}^{\wt{N}/2-|\alpha|/2-\ell}.
\end{align}
Now taking absolute values in \eqref{eq:ddKj4} and then employing the above estimates we obtain
\begin{align*}
&\big|(x_i-y_i)^{N}(A^{(x)})^{\gamma'} x^{\gamma''}(A^{(y)})^{\eta'}y^{\eta''} \R_j^\alpha(x,y)\big|\\
&\qquad\qquad \lesi 2^{-j\wt{N}}\sum_{\substack{s,\wt{s},t,\xi,\\\ell,k,\wt{k}}} \sum_{\substack{\om\le \gamma''\\\wt{\om}\le \eta''}}  \ip{\xi}^\Theta \,
 |h_{\xi +\gamma'+\gamma''-2\om+\wt{\alpha} + (k-s)e_i}(x)| |h_{\xi+\eta'+\eta''-2\wt{\om} +(\wt{k}-\wt{s})e_i}(y)|, 
\end{align*}
where $\Theta = \f{\wt{N}+|\gamma'|+|\gamma''|+|\eta'|+|\eta''|+k+\wt{k}-2\ell - s-\wt{s}-t}{2}$.

Now recall that $\ip{\xi}\sim 4^j$, $k+\wt{k}=2\ell-\wt{N}$ and that $\wt{N}=N-s-\wt{s}-t$. Then we have
\begin{align*}
&\big|(x_i-y_i)^{N}(A^{(x)})^{\gamma'} x^{\gamma''}(A^{(y)})^{\eta'}y^{\eta''} \R_j^\alpha(x,y)\big|\\
&\qquad \lesi 2^{j(|\gamma'|+|\gamma''|+|\eta'|+|\eta''|-N)}\sum_{\substack{s,\wt{s},t,\xi,\\\ell,k,\wt{k}}} \sum_{\substack{\om\le \gamma''\\\wt{\om}\le \eta''}} |h_{\xi +\gamma'+\gamma''-2\om+\wt{\alpha} + (k-s)e_i}(x)| |h_{\xi+\eta'+\eta''-2\wt{\om} +(\wt{k}-\wt{s})e_i}(y)| \\
&\qquad\lesi 2^{j(|\gamma|+|\eta|-N)}\QQ_{4^j +|\gamma|+|\wt{\alpha}|+N}(x,x)^\half \QQ_{4^j+|\eta|+N}(y,y)^\half,
\end{align*}
since $|\gamma'|+|\gamma''|\le |\gamma|$, $|\eta'|+|\eta''|\le |\eta|$ and $\wt{N}\le N$. We may now invoke \eqref{eq:QQ est} to arrive at \eqref{eq:ddKj0}, completing the proof of the lemma.

%%%%%%%%%%%%Bibliography  %%%%%%%%%%%%%%%%%%%%%%%%%%%

\end{document}